\documentclass[11pt,english,10pt,onecolumn,draftclsnofoot]{IEEEtran}
\usepackage[T1]{fontenc}
\usepackage[latin9]{inputenc}
\usepackage{array}
\usepackage{dsfont}
\usepackage{amsmath}
\usepackage{amssymb}
\usepackage{graphicx}

\makeatletter

\providecommand{\tabularnewline}{\\}

\newtheorem{assumption}{Assumption}
\newtheorem{remrk}{Remark}
\newtheorem{definitn}{Definition}
\newtheorem{lemma}{Lemma}
\newtheorem{thm}{Theorem}
\newtheorem{cor}{Corollary}
\newtheorem{example}{Example}

\usepackage{cite}
\usepackage{times}

\usepackage{fancyhdr} 
\usepackage{amsmath}

\@ifundefined{showcaptionsetup}{}{%
 \PassOptionsToPackage{caption=false}{subfig}}
\usepackage{subfig}
\makeatother

\usepackage{babel}
\begin{document}
\title{On the Value of Linear Quadratic Zero-sum Difference Games with Multiplicative
Randomness: Existence and Achievability}
\author{Songfu Cai, \textsl{Member, IEEE}, and Xuanyu Cao, \textsl{Senior
Member, IEEE}{\normalsize{}}\thanks{Songfu Cai is with the College of ISEE, Zhejiang University (email:
sfcai@zju.edu.cn) and Xuanyu Cao is with the School of EECS, Washington
State University (email: xuanyu.cao@wsu.edu).}\\
}
\maketitle
\begin{abstract}
We consider a wireless networked control system (WNCS) with multiple
controllers and multiple attackers. The dynamic interaction between
the controllers and the attackers is modeled as a linear quadratic
(LQ) zero-sum difference game with multiplicative randomness induced
by the the multiple-input and multiple-output (MIMO) wireless fading.
The existence of the game value is of significant importance to the
LQ zero-sum game. If it exists, the value characterizes the minimum
loss that the controllers can secure against the attackers' best possible
strategies. However, if the value does not exist, then there will
be no Nash equilibrium policies for the controllers and attackers.
Therefore, we focus on analyzing the existence and achievability of
the value of the LQ zero-sum game. In  stark contrast to existing
literature, where the existence of the game value depends heavily
on the controllability of the closed-loop systems, the multiplicative
randomness induced by the MIMO wireless channel fading may destroy
closed-loop controllability and introduce intermittent controllability
or almost sure uncontrollability. We first establish a general sufficient
and necessary condition for the existence of the game value. This
condition relies on the solvability of a modified game algebraic Riccati
equation (MGARE) under an implicit concavity constraint, which is
generally difficult to verify due to the multiplicative randomness.
We next introduce a novel positive semidefinite (PSD) kernel decomposition
method induced by multiplicative randomness. A verifiable tight sufficient
condition is then obtained by applying the proposed PSD kernel decomposition
to the constrained MGARE. Under the existence condition, we finally
construct a closed-form saddle-point policy based on the minimal solution
to the constrained MGARE, which is able to achieve the game value
in a certain class of admissible policies. We demonstrate that the
proposed saddle-point policy is backward compatible to the existing
strictly feedback stabilizing saddle-point policy. Our developed theory
is also verified via numerical simulations.

\medskip{}
\medskip{}
\medskip{}

\providecommand{\keywords}[1]{\textbf{\textit{Index terms---}}#1} 

\end{abstract}

\section{Introduction}

Wireless networked control systems (WNCSs) have drawn great attention
in both academia and industry in recent years. A wide spectrum of
WNCS applications can be found in areas such as aerospace control,
factory automation, mobile sensor networks, remote surgery, intelligent
transportation, unmanned aerial vehicles (UAV), etc. {[}1{]}, {[}2{]},
{[}3{]}. Compared to classical feedback control systems, networked
control systems have many advantages such as flexibility in communication
architectures, low cost in maintenance, and fast installation and
deployment {[}4{]}. Despite the great advantages offered by WNCS,
the wireless communication network with shared common spectrum in
between the dynamic plant and the remote controllers is often unreliable
and vulnerable to potential threats and malicious attacks. The attackers
can intrude the shared wireless communication network and deteriorate
the closed-loop control performance, as illustrated in Fig. 1. The
interaction between the controllers and the attackers is typically
modeled as a linear quadratic (LQ) zero-sum game {[}5{]}, {[}6{]},
{[}7{]}, where the controllers try to minimize the negative impacts
of adversarial attacks and the attackers aim to maximize the disruption
to the system. From a practical point of view, the value of the LQ
zero sum game, if it exists, is equal to the loss that the controllers
can secure against any strategy of the attackers, and can prevent
the attackers from doing any better to cause more than that loss.
Therefore, analyzing the existence conditions of the game value is
of paramount importance for the controllers to comprehend the worst-case
scenario and to prepare defenses accordingly, while also providing
insight into how the attackers might strategize. 

\begin{figure}[t!]
\centering{}\includegraphics[width=0.7\columnwidth]{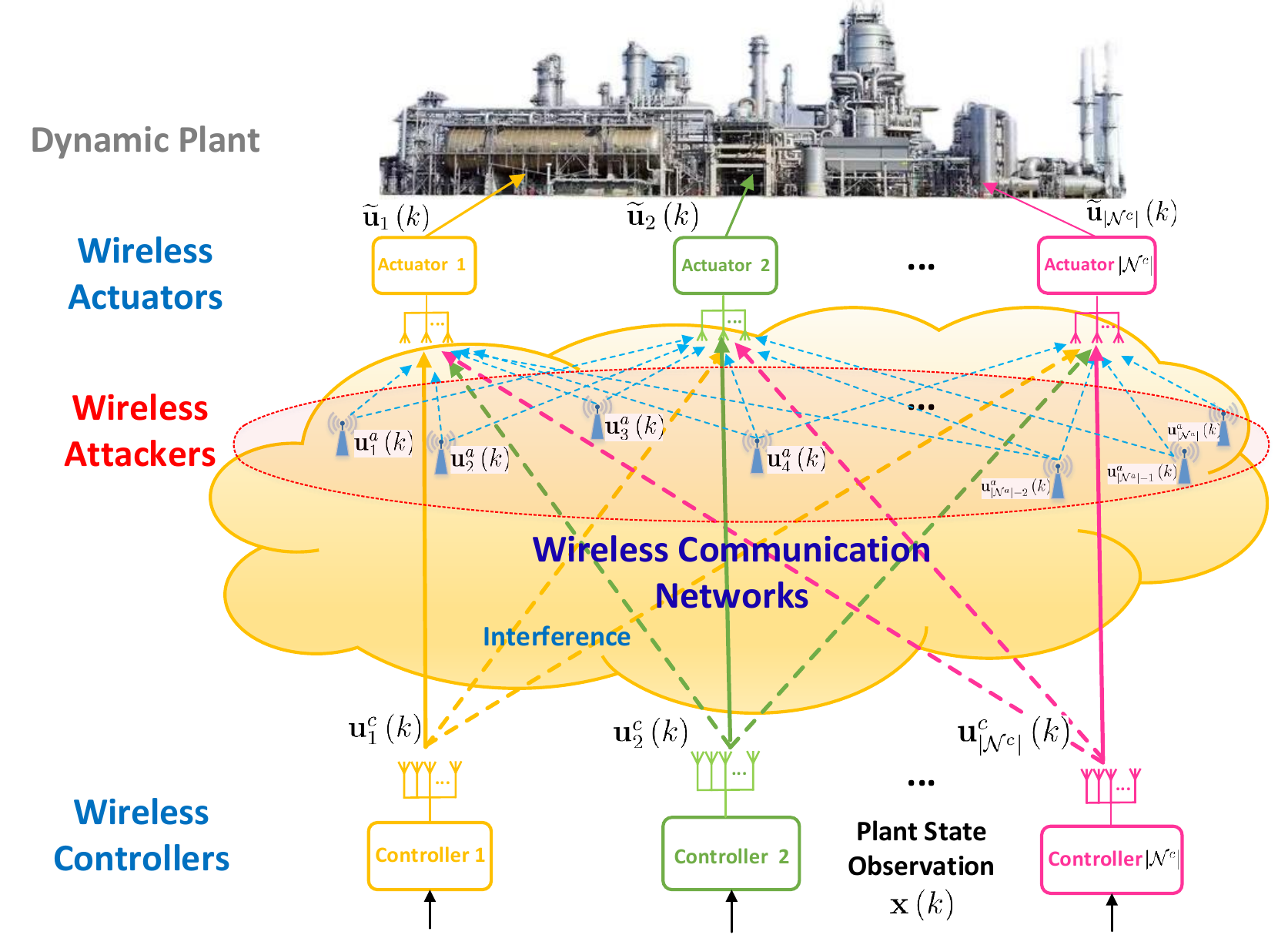} \caption{\label{fig-system-topo}Illustration of a WNCS with multiple wireless
controllers and multiple wireless attackers. The remote controllers
transmit the control actions to the actuators over a shared MIMO wireless
communication network. The wireless attackers can intrude the wireless
communication network and inject malicious attacks to the actuators.}
\end{figure}

Several existing works have studied LQ zero-sum differential games.
The seminal work by Ho, Bryson, and Baron {[}8{]} developed the zero-sum
framework for the deterministic LQ pursuit-evasion differential games.
Based on Isaacs' Verification Theorem {[}9{]}, they derived a sufficient
condition for achieving the value of the game through closed-loop
saddle-point policies. Ba\c{s}ar and Bernhard {[}10{]}, {[}11{]},
{[}12{]} established the connections between the LQ zero-sum differential
games and the disturbance attenuation in the $H_{\infty}$ optimal
control, where the controller and the disturbance are considered as
minimizing and maximizing players, respectively. Mou and Yong {[}13{]}
considered a class of two-person LQ zero-sum differential games. They
identified the connection between the existence of the game values
and the solvability of a generalized differential Riccati equation
(GDRE). Goldys, Yang and Zhou {[}14{]} further investigated a class
of dynamical systems with coupled slow and fast state processes and
showed that, when the slow-fast state coupling is weak enough, the
solvability of the associated GDRE could be determined by the solvability
of a decoupled pair of differential algebraic Riccati equations. Moon
{[}15{]} studied the LQ zero-sum differential games for Markov jump
systems. It was shown that the existence of a feedback saddle point
is determined by the solvability of a set of coupled Riccati differential
equations (CRDEs). Other notable works, such as those by Schmitendorf
{[}16{]}, Zhang {[}17{]}, Sun {[}18{]}, and Yu {[}19{]}, provided
comprehensive studies on closed-loop saddle-points, open-loop saddle-points,
and the values of LQ zero-sum differential games.

In the studies of LQ zero-sum difference games, the game algebraic
Riccati equation (GARE) is of paramount importance. For linear time-invariant
(LTI) systems, pioneering works by Ba\c{s}ar {[}20{]} revealed that
finding the saddle-point solution to the LQ zero-sum difference game
involves solving a GARE. Tremendous effort has been dedicated to developing
iterative methods for solving the GARE, assuming the existence of
the game value. Luo, Yang, and Liu {[}21{]} proposed a data-based
policy iteration Q-learning algorithm, while Kiumarsi, Lewis, and
Jiang {[}22{]} used off-policy reinforcement learning (RL) and demonstrated
its robustness to probing noises. Al-Tamimi, Lewis and Abu-Khalaf
{[}23{]} developed a novel actor-critic-based online learning algorithm,
where the critic converges to the game value function and the actor
networks converge to the Nash equilibrium (NE) of the game. Zhang,
Yang and Ba\c{s}ar {[}24{]} adopted the policy optimization framework
to optimize the state feedback gain matrices of the players directly.
They investigated the optimization landscape and proved the global
convergence of policy optimization methods for finding the NE. Several
other works have considered the LQ zero-sum difference games for linear
time-varying (LTV) systems with model uncertainties. Specifically,
Wu et al. {[}25{]} considered a stochastic LQ zero-sum difference
game for WNCS with one controller and one attacker, where the attacker's
malicious injection is governed by a Bernoulli random variable. Under
concavity assumption, a sufficient condition for the existence of
game value is obtained in terms of the attack injection probability.
Gravell, Ganapathy and Summers {[}26{]} studied an LTV system, where
the state transition matrix, the control input matrix and disturbance
input matrix are perturbed by additive random noises. Under the assumption
that the associated GARE is solvable, two policy iteration algorithms
are proposed to compute the equilibrium strategies and the value of
the game. All the aforementioned existing works {[}20{]}-{[}26{]}
require a common assumption that the closed-loop system is controllable,
i.e., the pair $\left(\mathbf{A},\mathbf{B}^{c}\right)$ is controllable.
This controllability assumption guarantees the existence of a positive
definite solution to the GARE associated with the LQ zero-sum difference
game considered in {[}20{]}-{[}26{]}. Furthermore, in order to guarantee
the existence of the value of the zero-sum game, the attacker's weight
matrix, e.g., $\gamma\mathbf{I}$ in {[}20{]}-{[}23{]}, $\mathbf{R}^{\upsilon}$
in {[}24{]}, $\mathbf{R}_{2}$ in {[}25{]} and $\mathbf{S}$ in {[}26{]},
is directly assumed to satisfy an implicit concavity constraint. However,
the explicit closed-form characterization on the feasible domain of
the attacker's weight matrix has never been investigated.

\begin{figure}[t!]
\centering{}\includegraphics[width=0.9\columnwidth]{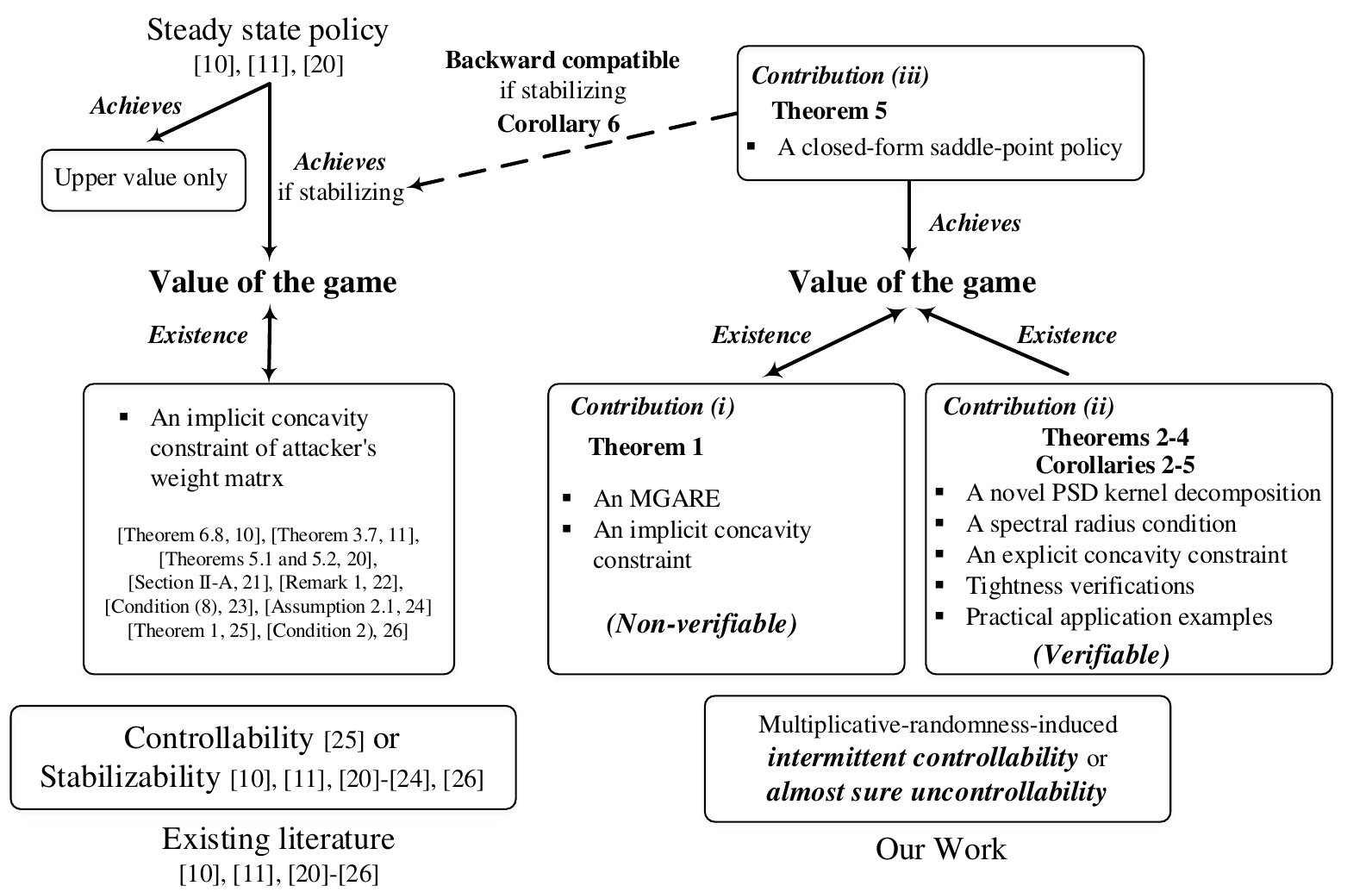}
\caption{\label{fig-summary-contribution} Summary of the key contributions.}
\end{figure}

In this paper, we adopt the LQ zero-sum difference game framework
to analyze the WNCS with multiple remote controllers and multiple
attackers over MIMO wireless fading channels. In stark contrast to
the existing literature {[}20{]}-{[}26{]} where closed-loop controllability
is a prerequisite, the multiplicative randomness induced by the MIMO
wireless channel fading may destroy controllability and introduce
intermittent controllability or almost sure uncontrollability. We
investigate various conditions, including the explicit closed-form
characterizations on the feasible region of the attackers' weight
matrices, that guarantee the existence of the game value. Our key
contributions are outlined below and are briefly depicted in Fig.
\ref{fig-summary-contribution}.

\begin{itemize}
\item \textbf{General Sufficient and Necessary Condition for Existence of
Game Value: }In existing literature, sufficient and necessary conditions
for the existence of the value of the zero-sum game rely heavily on
the controllability of the closed-loop system {[}20{]}-{[}25{]}. However,
in our work, the multiplicative randomness induced by the MIMO wireless
channel fading may introduce intermittent controllability or almost
sure uncontrollability, rendering the existing conditions inapplicable.
To address this issue, we establish a general sufficient and necessary
condition in terms of the solvability of a modified game algebraic
Riccati equation (MGARE) subject to an implicit concavity constraint.
We have also proposed an equivalent analytic condition based on the
modified game algebraic Riccati recursions. 
\item \textbf{Closed-form Tight Verifiable Sufficient Condition:} Finding
a closed-form verifiable sufficient condition for the existence of
the value of a zero-sum game is generally difficult, even for deterministic
zero-sum difference games. The current literature {[}21{]}, {[}23{]}-{[}26{]}
presumes the game value exists and only concentrates on iterative
methods to solve the associated GARE. To the best of our knowledge,
closed-form tight verifiable conditions for the existence of the value
of stochastic zero-sum difference games have not been studied in the
literature. In this paper, we introduce a new PSD kernel decomposition
induced by the multiplicative randomness. A verifiable sufficient
condition is then established by applying this decomposition to both
the MGARE and the implicit concavity constraint. We discover that
the solvability of the MGARE is closely related to the existence of
a positive definite solution to a generalized Lyapunov equation (GLE).
Besides, the implicit concavity constraint is also transformed into
a verifiable explicit constraint. Additionally, we demonstrate the
tightness of the proposed verifiable sufficient condition using a
specific application example.
\item \textbf{Closed-form Game Value Achieving Saddle-point Policy:} It
is well-known that, for deterministic controllable systems, the steady
state policy associated with the minimal solution to the GARE is usually
not a saddle-point policy because it can only achieve the upper value
of the zero-sum game {[}11{]}, {[}20{]}. However, finding a closed-form
saddle-point policy is often difficult. To address this challenge,
a novel timestamp based on the MAGRE and the growth of the plant state's
Euclidean norm is defined. A closed-form policy is then constructed
such that it adopts the steady state policy associated with the minimal
solution to the MGARE for time slots that are less than the timestamp,
and another closed-form time-varying policy for all time slots that
are larger than the timestamp. It is proven that this proposed policy
is a game value achieving saddle-point policy in the class of admissible
policies with a prescribed plant state growth. 
\end{itemize}

\begin{figure}[t!]
\centering{}\includegraphics[width=0.4\columnwidth]{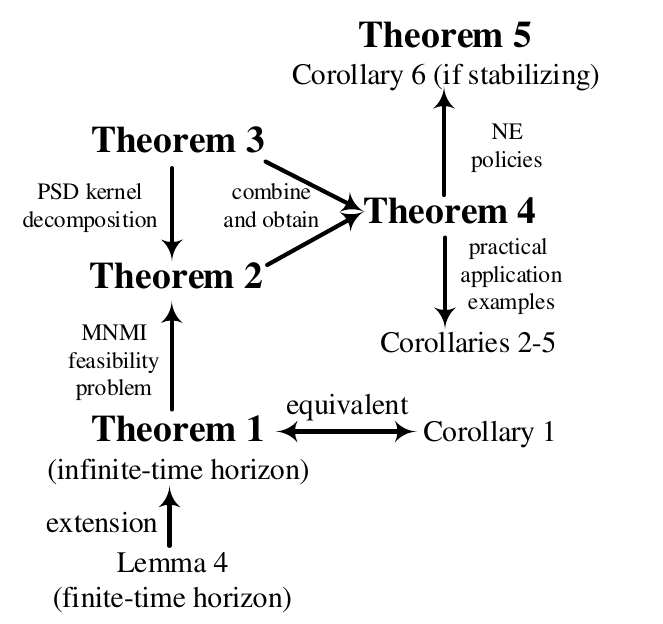}
\caption{\label{fig-summary-connection} Connections between the key theorems,
lemmas and corollaries.}
\end{figure}

A flowchart that elucidates the connections between the key lemmas,
theorems and corollaries in this work is provided in Fig. \ref{fig-summary-connection}.
The notations that will be used throughout this paper are summarized
below.

\begin{itemize}
\item[\textendash{}]  $\mathbb{S}_{+}^{S}$ and $\hat{\mathbb{S}}_{+}^{S}$ denote the
set of $S\times S$ dimensional positive definite matrices and positive
semidefinite matrices, respectively. $\mathbb{Z}^{+}$ and $\mathbb{Z}_{0}^{+}$
denote the set of positive integers and nonnegative integers, respectively.
\textrm{Ø} is the null set. $\left|\mathcal{I}\right|$ denotes the
cardinality of an integer set $\mathcal{I}$. $\mathcal{I}\setminus i$
denotes an integer set $\mathcal{I}$ with element $i$ removed. 
\item[\textendash{}]  $\left(\cdot\right)^{T}$ and $\mathrm{diag}\left(\cdot\right)$
denote the matrix transpose and the diag operator, respectively. $\mathrm{vec}\left(\cdot\right)$
is the vectorization operator, and $\mathrm{vec}^{-1}\left(\cdot\right)$
is the inverse vectorization operator of a vector to a square matrix
with appropriate dimension. $\mathbf{0}_{m}$ and $\mathbf{I}_{m}$
denote the $m\times m$ dimensional null matrix and identity matrix,
respectively.
\item[\textendash{}] $\Vert\mathbf{a}\Vert_{2}$ denotes the Euclidean norm of vector
$\mathbf{a}$. $\Vert\mathbf{a}\Vert_{\mathbf{Q}}=\Vert\mathbf{Q}^{\frac{1}{2}}\mathbf{a}\Vert_{2}$
denotes the Euclidean norm weighted by a positive definite matrix
$\mathbf{Q}$. 
\item[\textendash{}] $\mathrm{rank}\left(\mathbf{A}\right),$ $\mathrm{ker}\left(\mathbf{A}\right)$,
$\Vert\mathbf{A}\Vert$, $\sigma_{\mathrm{max}}\left(\mathbf{A}\right)$,
$\sigma_{\mathrm{min}}\left(\mathbf{A}\right)$, $\sigma_{i}\left(\mathbf{A}\right)$
and $\mathrm{\rho}\left(\mathbf{A}\right)$ denote the rank, left
kernel, spectral norm, the maximum singular value, the minimum singular
value, the $i$-th largest singular value and the spectral radius
of matrix $\mathbf{A}$. $\mathrm{cov}\left(\mathbf{A}\right)$ denotes
the covariance matrix of a random matrix $\mathbf{A}$.
\item[\textendash{}] $\mathbf{A}_{\left(i,j\right)}$ denotes the element in the $i$-th
row and $j$-th column of matrix $\mathbf{A}$. $\left(\mathbf{A}\right)_{i}$
denotes the $i$-th order leading principal sub-matrix of $\mathbf{A}$.
$\mathbf{\left(A\right)}_{\left(i:j;l:m\right)}$ denotes the $(j-i+1)\times(m-l+1)$
dimensional block sub-matrix of $\mathbf{A}$ with the first element
being $\mathbf{A}_{\left(i,l\right)}$. $\mathbf{\left(A\right)}_{\left(i:j;:\right)}$
denotes the block sub-matrix of $\mathbf{A}$ that contains a total
number of $\left(j-i+1\right)$ consecutive rows of $\mathbf{A}$
starting from the $i$-th row until the $j$-th row of $\mathbf{A}$.
\item[\textendash{}] $\left[\mathbf{A}_{i}\right]_{i\in\mathcal{I}}\in\mathbb{R}^{m\times\left(n\left|\mathcal{I}\right|\right)}$
denotes a matrix constructed by the row concatenation of each $\mathbf{A}_{i}\in\mathbb{R}^{m\times n}$,
whose subscript $i$ lies in the integer set $\mathcal{I}$, in subscript
ascending order, i.e., $\left[\mathbf{A}_{i}\right]_{i\in\mathcal{I}}=[\mathbf{A}_{i_{1}},\mathbf{A}_{i_{2}},$
$\cdots,\mathbf{A}_{i_{\left|\mathcal{I}\right|}}]$ with $i_{1}<i_{2}<\cdots<i_{\left|\mathcal{I}\right|}.$
$\left[\mathbf{A}_{i,j}\right]_{\left(i,j\right)\in\mathcal{I}\times\mathcal{J}}\in\mathbb{R}^{m\times\left(n\left|\mathcal{I}\right|\left|\mathcal{J}\right|\right)}$
denotes a matrix constructed by the row concatenation of each $\mathbf{A}_{i,j}\in\mathbb{R}^{m\times n}$,
whose subscript $\left(i,j\right)$ lies in the Cartesian product
$\mathcal{I}\times\mathcal{J}$ of the integer sets $\mathcal{I}$
and $\mathcal{J}$, in subscript ascending order, i.e., $\left[\mathbf{A}_{i,j}\right]_{\left(i,j\right)\in\mathcal{I}\times\mathcal{J}}=\big[\mathbf{A}_{i_{1},j_{1}},\cdots,\mathbf{A}_{i_{1},j_{\left|\mathcal{J}\right|}},\cdots,\mathbf{A}_{j_{\left|\mathcal{I}\right|},j_{1}},$
$\cdots,\mathbf{A}_{j_{\left|\mathcal{I}\right|},j_{\left|\mathcal{J}\right|}}\big]$
with $i_{1}<i_{2}<\cdots<i_{\left|\mathcal{I}\right|}$ and $j_{1}<j_{2}<\cdots<j_{\left|\mathcal{J}\right|}$.
$\left\{ \mathbf{A}\left(i\right)\right\} _{i\in\mathcal{I}}$ (resp.,
$\left\{ \mathbf{A}_{i}\right\} _{i\in\mathcal{I}}$) denotes a set
that containing all the $\mathbf{A}\left(i\right)$'s (resp., $\mathbf{A}_{i}$'s)
with $i\in\mathcal{I}$. $\left\{ \mathbf{A}_{i,j}\right\} _{\left(i,j\right)\in\mathcal{I}\times\mathcal{J}}$
denotes 
\item[\textendash{}] $\mathbf{A}\prec\mathbf{B}$ and $\mathbf{A}\preceq\mathbf{B}$ mean
$\mathbf{B}-\mathbf{A}$ is positive definite and positive semidefinite,
respectively.
\item[\textendash{}]  For a block matrix $\mathbf{A}=\left[\begin{array}{cc}
\mathbf{A}_{1} & \mathbf{A}_{2}\\
\mathbf{A}_{3} & \mathbf{A}_{4}
\end{array}\right]$, if $\mathbf{A}_{1}$ is invertible, the Schur complement of the
block $\mathbf{A}_{1}$ of the matrix $\mathbf{A}$ is defined as
$\mathbf{A}\big/\mathbf{A}_{1}=\mathbf{\mathbf{A}}_{4}-\mathbf{A}_{3}\mathbf{A}_{1}^{-1}\mathbf{A}_{2}$.
Similarly, if $\mathbf{A}_{4}$ is invertible, the Schur complement
of the block $\mathbf{A}_{4}$ of the matrix $\mathbf{A}$ is defined
as $\mathbf{A}\big/\mathbf{A}_{4}=\mathbf{\mathbf{A}}_{1}-\mathbf{A}_{2}\mathbf{A}_{4}^{-1}\mathbf{A}_{3}$.
\item[\textendash{}]  $\mathrm{Pr}\left(A\right)$ and $\mathds{1}_{A}$ denote the probability
and the indicator function of event $A$, respectively. 
\item[\textendash{}] $\Omega\left(\cdot\right)$ is the big omega notation and $f\left(x\right)=\Omega\left(g\left(x\right)\right)$
if $\limsup_{x\rightarrow\infty}\left|\frac{f\left(x\right)}{g\left(x\right)}\right|>0.$
\end{itemize}

To enhance clarity, a summary of the key technical symbols is provided
in Table 1.

\begin{table}
\begin{centering}
\begin{tabular}{|>{\centering}m{0.3\columnwidth}|>{\centering}m{0.55\columnwidth}|}
\hline 
$\mathcal{N}^{c}$ (resp., $\mathcal{N}^{a}$) & the index set of all the controllers (resp., attackers)\tabularnewline
\hline 
$\mathbf{H}_{j,i}^{c}\left(k\right)$ (resp., $\mathbf{H}_{j,i}^{a}\left(k\right)$
) & the wireless MIMO fading channel matrix

between the $j$-th controller (resp., $l$-th attacker) and the $i$-th
actuator\tabularnewline
\hline 
$\mathbf{B}^{c}\left(k\right)$ (resp., $\mathbf{B}^{a}\left(k\right)$) & the concatenated aggregated CSI of the controllers (resp., attackers)\tabularnewline
\hline 
$\mathcal{T}_{k}^{\mathrm{ker}}$(resp., $\mathcal{T}_{k}^{\mathrm{\mathbf{0}}}$) & the controllable (resp. uncontrollable) cone associated with $\mathbf{B}^{c}\left(k\right)$\tabularnewline
\hline 
$\mathbf{T}=\mathbf{T}_{k}^{\mathrm{ker}}+\mathbf{T}_{k}^{\mathrm{\mathbf{0}}}$ & the PSD kernel decomposition for any positive definite matrix $\mathbf{T}$
such that $\mathbf{T}_{k}^{\mathrm{ker}}\in\mathcal{T}_{k}^{\mathrm{ker}}$
and $\mathbf{T}_{k}^{\mathrm{\mathbf{0}}}\in\mathcal{T}_{k}^{\mathrm{\mathbf{0}}}$\tabularnewline
\hline 
$\mathbf{P}=f\left(\mathbf{P}\right)$ & the MGARE\tabularnewline
\hline 
$f^{k}\left(\mathbf{P}\right)$ & the $k$-step modified game algebraic Riccati recursion\tabularnewline
\hline 
$\mathcal{P}_{\mathbf{R}^{a}}$ & the positive definite solution set of the MGARE\tabularnewline
\hline 
$\mathcal{R}_{\mathbf{R}^{a}}$ & the concavity constraint set of the attackers' weight matrix $\mathbf{R}^{a}$\tabularnewline
\hline 
$\left(\mathbf{u}^{c}\right)^{*}\left(\mathbf{P};k\right)$ (resp.,
$\left(\mathbf{u}^{a}\right)^{*}\left(\mathbf{P};k\right)$) & the steady state control (resp., attacker) actions\tabularnewline
\hline 
$\overline{J}_{K}$ (resp., $\underline{J}_{K}$) & the upper (resp., lower) value representation of the $K$-length finite-time
horizon average cost\tabularnewline
\hline 
$\overline{J}$ (resp., $\underline{J}$) & the upper (resp., lower) value of the infinite-time horizon LQ zero-sum
game\tabularnewline
\hline 
$J^{*}$ & the value of the infinite-time horizon LQ zero-sum game\tabularnewline
\hline 
$\alpha\left(T_{0}\right),$ $\beta\left(T_{0}\right)$ & constants quantifying the expected growth of the plant state's norm
in the first $T_{0}$ time-slots\tabularnewline
\hline 
$\left\{ \overline{\boldsymbol{\mu}}_{\beta\left(\infty\right)}^{c},\overline{\boldsymbol{\mu}}_{\beta\left(\infty\right)}^{a}\right\} $ & a Nash equilibrium solution that achieves the value of the infinite-time
horizon zero-sum game\tabularnewline
\hline 
\end{tabular}
\par\end{centering}
\caption{Summary of key technical symbols.}
\end{table}

\section{System Model}

In this section, we introduce the key components of the wireless networked
control system with multi-controllers and multi-attackers and present
the LQ zero-sum game formulation.

\subsection{Dynamic Plant Model}

We consider a time-slotted wireless networked control system, which
contains a potentially unstable dynamic plant, $\left|\mathcal{N}^{c}\right|$
actuators, $\left|\mathcal{N}^{c}\right|$ remote controllers, and
$\left|\mathcal{N}^{a}\right|$ remote attackers. The actuators are
colocated with the dynamic plant, whereas the controllers and attackers
are geographically distributed and are connected to the actuators
through a wireless communication network. 

The dynamic evolution of the plant state $\mathbf{x}_{k}$ is summarized
below.

\begin{assumption}
\emph{(Plant State Model)} The plant state $\mathbf{x}_{k}$ follows
the dynamics of
\begin{align}
 & \mathbf{x}\left(k+1\right)=\mathbf{Ax}\left(k\right)+\sum_{i\in\mathcal{N}^{c}}\mathbf{B}_{i}\widetilde{\mathbf{u}}_{i}\left(k\right)+\mathbf{w}\left(k\right),k\in\mathbb{Z}_{0}^{+},\label{eq:sys-dynamics}
\end{align}
where $\mathbf{x}\left(k\right)\in\mathbb{R}^{S\times1}$ is the state
vector, $S$ is the state dimension; $\mathbf{A}\in\mathbb{R}^{S\times S}$
is the state transition matrix; $\mathcal{N}^{c}=\left\{ 1,2,\cdots,\left|\mathcal{N}^{c}\right|\right\} $,
$\mathbf{B}_{i}\in\mathbb{R}^{S\times N_{r}}$ and $\widetilde{\mathbf{u}}_{i}\left(k\right)\in\mathbb{R}^{N_{r}\times1}$
is the $i$-th actuator's control input gain matrix and control action,
respectively; $\mathbf{w}\left(k\right)\in\mathbb{R}^{S\times1}$
is the additive plant noise. The initial state $\mathbf{x}\left(0\right)$
is a constant vector.
\end{assumption}

\subsection{Wireless Communication Model}

The wireless communication channel between each controller-actuator
pair and each attacker-actuator pair is modeled as a wireless MIMO
fading channel. We assume each actuator, remote controller, remote
attacker has $N_{r}$ receive antennas, $N_{t}^{c}$ and $N_{t}^{a}$
transmit antennas, respectively.At the $k$-th time-slot, the received
signals $\widetilde{\mathbf{u}}_{i}\left(k\right)\in\mathbb{R}^{N_{r}\times1}$
at the $i$-th actuator is given by
\begin{align}
\widetilde{\mathbf{u}}_{i}\left(k\right)= & \mathbf{H}_{i,i}^{c}\left(k\right)\mathbf{u}_{i}^{c}\left(k\right)+\sum_{j\in\mathcal{N}^{c},j\neq i}\mathbf{H}_{j,i}^{c}\left(k\right)\mathbf{u}_{j}^{c}\left(k\right)\label{eq:u_i(k)}\\
 & +\sum_{l\in\mathcal{N}^{a}}\mathbf{H}_{l,i}^{a}\left(k\right)\mathbf{u}_{l}^{a}\left(k\right)+\mathbf{v}_{i}\left(k\right),i\in\mathcal{N}^{c},k\in\mathbb{Z}_{0}^{+},\nonumber 
\end{align}
where $\mathbf{u}_{j}^{c}\left(k\right)$ and $\mathbf{u}_{l}^{a}\left(k\right)$
is the transmitted control action of the $j$-th controller and the
malicious attack injection of the $l$-th attacker, respectively;
$\mathcal{N}^{a}=\left\{ 1,2,\cdots,\left|\mathcal{N}^{a}\right|\right\} $
is the index set of the attackers; $\mathbf{H}_{j,i}^{c}\left(k\right)\in\mathbb{R}^{N_{r}\times N_{t}^{c}}$
and $\mathbf{H}_{l,i}^{a}\left(k\right)\in\mathbb{R}^{N_{r}\times N_{t}^{a}}$
are random matrices that denote the wireless MIMO fading channel between
the $j$-th controller and the $i$-th actuator and between the $l$-th
attacker and the $i$-th actuator, respectively; and $\mathbf{v}_{i}\left(k\right)$
is the additive wireless channel noise at the $i$-th actuator. 

We make the following assumption on the additive randomness of the
system.

\begin{assumption}
\emph{(}\textit{Additive Randomness}\emph{)} \label{Assumption-Additive-Randomness}The
additive plant noise and the channel noise process satisfy the following
conditions: 
\end{assumption}

\begin{itemize}
\item $\mathbf{w}\left(k\right)$ is i.i.d. over the time-slots, i.e., \textbf{$\mathbf{w}\left(k_{1}\right)$}
and $\mathbf{w}\left(k_{2}\right)$ are i.i.d., $\forall k_{1},k_{2}\in\mathbb{Z}_{0}^{+},k_{1}\neq k_{2}$,
with zero mean and finite covariance matrix $\mathbf{W}$.
\item $\mathbf{v}_{i}\left(k\right),\forall i\in\mathcal{N}^{c},$ is i.i.d.
across time and actuators with zero mean and finite covariance matrix
$\mathbf{V}$.
\item The realizations of $\mathbf{w}\left(k\right)$ and $\mathbf{v}_{i}\left(k\right),\forall i\in\mathcal{N}^{c},$
remain constant within each time-slot $k$, $\forall k\in\mathbb{Z}_{0}^{+}$.
\end{itemize}

The i.i.d. additive dynamic plant noise model in Assumption \ref{Assumption-Additive-Randomness}
characterizes the environmental factors of dynamic systems {[}27{]},
e.g., random changes in temperature or humidity, the random unpredictable
disturbances from external sources {[}28{]}, and the kinetic approximation
errors in system modeling {[}29{]}, etc. It has been extensively adopted
in many practical applications of stochastic networked control {[}30{]},
{[}31{]}, as well as remote state estimation {[}32{]}, {[}33{]}. The
i.i.d. additive channel model is the most commonly used channel model
in wireless communications {[}34{]}. It serves as an important foundation
for the classical Shannon information theory and has wide applications
in channel coding design, channel capacity characterization, and rate
distortion analysis. It is also an ideal model for many practical
applications, including satellite {[}35{]}, deep space {[}36{]} and
terrestrial communication {[}37{]}. 

We make the following assumptions on the multiplicative randomness
of the system.

\begin{assumption}
\emph{(}\textit{MIMO Wireless Fading Channel Model}\emph{)\label{MIMO-Wireless-Fading}}
The random MIMO fading processes of the remote controllers $\left\{ \left[\mathbf{H}_{j,i}^{c}\left(k\right)\right]_{\left(j,i\right)\in\mathcal{N}^{c}\times\mathcal{N}^{c}}\right\} _{k\in\mathbb{Z}_{0}^{+}}$
and the remote attackers $\left\{ \left[\mathbf{H}_{l,i}^{a}\left(k\right)\right]_{\left(l,i\right)\in\mathcal{N}^{a}\times\mathcal{N}^{c}}\right\} _{k\in\mathbb{Z}_{0}^{+}}$
satisfy the following properties:
\end{assumption}

\begin{itemize}
\item For each time-slot $k\in\mathbb{Z}_{0}^{+}$, the concatenated wireless
MIMO channel matrices of the controllers $\left[\mathbf{H}_{j,i}^{c}\left(k\right)\right]_{\left(j,i\right)\in\mathcal{N}^{c}\times\mathcal{N}^{c}}$
follows some general distribution. Besides, $\left[\mathbf{H}_{j,i}^{c}\left(k\right)\right]_{\left(j,i\right)\in\mathcal{N}^{c}\times\mathcal{N}^{c}}$
is i.i.d. over the time-slots. 
\item For each time-slot $k\in\mathbb{Z}_{0}^{+}$, the concatenated wireless
MIMO channel matrices of the attackers $\left[\mathbf{H}_{l,i}^{a}\left(k\right)\right]_{\left(l,i\right)\in\mathcal{N}^{a}\times\mathcal{N}^{c}}$
follows some general distribution. Besides, $\left[\mathbf{H}_{j,i}^{c}\left(k\right)\right]_{\left(j,i\right)\in\mathcal{N}^{c}\times\mathcal{N}^{c}}$
is i.i.d. over the time-slots. 
\item The realizations of $\left[\mathbf{H}_{j,i}^{c}\left(k\right)\right]_{\left(j,i\right)\in\mathcal{N}^{c}\times\mathcal{N}^{c}}$
and $\left[\mathbf{H}_{l,i}^{a}\left(k\right)\right]_{\left(l,i\right)\in\mathcal{N}^{a}\times\mathcal{N}^{c}}$
remain constant within each time-slot $k$, $\forall k\in\mathbb{Z}_{0}^{+}$.
\end{itemize}

Assumption \ref{MIMO-Wireless-Fading} is a commonly used model for
wireless communications, known as the block-fading channel model {[}38{]},
{[}39{]}. In this model, the wireless channel gains remain constant
within a fading block, which is represented by the time index $k$,
but is i.i.d. from one block to another. Each block is referred to
as the coherence interval of the wireless channel and is equal to
the duration of a time-slot. 

We also have the following independence assumption.

\begin{assumption}
\emph{(}\textit{Mutual Independence}\emph{)\label{Indenpendence}}
The multiplicative randomness processes $\left\{ \left[\mathbf{H}_{j,i}^{c}\left(k\right)\right]_{\left(j,i\right)\in\mathcal{N}^{c}\times\mathcal{N}^{c}}\right\} _{k\in\mathbb{Z}_{0}^{+}}$
and $\left\{ \left[\mathbf{H}_{l,i}^{a}\left(k\right)\right]_{\left(l,i\right)\in\mathcal{N}^{a}\times\mathcal{N}^{c}}\right\} _{k\in\mathbb{Z}_{0}^{+}},$
and the additive randomness processes $\left\{ \mathbf{w}\left(k\right)\right\} _{k\in\mathbb{Z}_{0}^{+}}$
and $\left\{ \left[\mathbf{v}_{i}\left(k\right)\right]_{i\in\mathcal{N}^{c}}\right\} _{k\in\mathbb{Z}_{0}^{+}}$
are mutually independent.
\end{assumption}

For Assumption \ref{Indenpendence}, the multiplicative randomness
and additive randomness are independent of each other as they originate
from different sources of randomness. 

\subsection{Zero-sum Game Formulation}

Substitute (\ref{eq:u_i(k)}) into (\ref{eq:sys-dynamics}), the closed-loop
system dynamics can be represented as
\begin{align}
\mathbf{x}\left(k+1\right)= & \mathbf{Ax}\left(k\right)+\mathbf{B}^{c}\left(k\right)\mathbf{u}^{c}\left(k\right)+\mathbf{B}^{a}\left(k\right)\mathbf{u}^{a}\left(k\right)\label{eq:closed-loop dynamics}\\
 & +\sum_{i\in\mathcal{N}^{c}}\mathbf{B}_{i}\mathbf{v}_{i}\left(k\right)+\mathbf{w}\left(k\right),\nonumber 
\end{align}
where
\begin{align}
 & \mathbf{B}^{c}\left(k\right)=\left[\sum_{i\in\mathcal{N}^{c}}\mathbf{B}_{i}\mathbf{H}_{j,i}^{c}\left(k\right)\right]_{j\in\mathcal{N}^{c}}\in\mathbb{R}^{S\times N_{t}^{c}\left|\mathcal{N}^{c}\right|},\label{eq:agg_CSI_B^c}\\
 & \mathbf{B}^{a}\left(k\right)=\left[\sum_{i\in\mathcal{N}^{c}}\mathbf{B}_{i}\mathbf{H}_{l,i}^{a}\left(k\right)\right]_{l\in\mathcal{N}^{a}}\in\mathbb{R}^{S\times N_{t}^{a}\left|\mathcal{N}^{a}\right|},\label{eq:agg_CSI_B^a}
\end{align}
{\small{}$\mathbf{u}^{c}\left(k\right)=\left[\left(\mathbf{u}_{j}^{c}\left(k\right)\right)^{T}\right]_{j\in\mathcal{N}^{c}}^{T}\in\mathbb{R}^{N_{t}^{c}\left|\mathcal{N}^{c}\right|\times1}$}
and {\small{}$\mathbf{u}^{a}\left(k\right)=\left[\left(\mathbf{u}_{l}^{a}\left(k\right)\right)^{T}\right]_{l\in\mathcal{N}^{a}}^{T}\in\mathbb{R}^{N_{t}^{a}\left|\mathcal{N}^{a}\right|\times1}$}.
The coefficient $\sum_{j\in\mathcal{N}^{c}}$ $\mathbf{B}_{i}\mathbf{H}_{j,i}^{c}\left(k\right)$
in (\ref{eq:agg_CSI_B^c}) is the aggregated channel state information
(CSI) of the $j$-th controller and $\sum_{i\in\mathcal{N}^{c}}\mathbf{B}_{i}\mathbf{H}_{l,i}^{c}\left(k\right)$
in (\ref{eq:agg_CSI_B^a}) is the aggregated CSI of the $l$-th attacker. 

\begin{remrk}
\textit{(Intermittent Controllability and Uncontrollability)}\emph{\label{Remark: Unctrl+intermitent Ctrl-1}}
The intermittent controllability refers to the case that, for any
time slot $k\in\mathbb{Z}_{0}^{+},$ the pair $\left(\mathbf{A},\mathbf{B}^{c}\left(k\right)\right)$
can be either controllable or uncontrollable for different realizations
of $\mathbf{B}^{c}\left(k\right)$ {[}3{]}, {[}40{]}. For example,
let $\mathbf{A}=\left[\begin{array}{cc}
1 & 0\\
0 & 1
\end{array}\right]$, $\mathrm{Pr}\left(\mathbf{B}^{c}\left(k\right)=\left[\begin{array}{cc}
1 & 0\\
0 & 1
\end{array}\right]\right)=\frac{1}{2}$ and $\mathrm{Pr}\left(\mathbf{B}^{c}\left(k\right)=\mathbf{0}\right)=\frac{1}{2},\forall k\in\mathbb{Z}_{0}^{+},$
the pair $\left(\mathbf{A},\mathbf{B}^{c}\left(k\right)\right)$ is
intermittent controllable. The almost surely uncontrollability means
that, for every time slot $k\in\mathbb{Z}_{0}^{+},$ the pair $\left(\mathbf{A},\mathbf{B}^{c}\left(k\right)\right)$
is uncontrollable with probability 1 {[}41{]}. Consider another toy
example that $\mathbf{A}=\left[\begin{array}{cc}
1 & 0\\
0 & 1
\end{array}\right]$, $\mathrm{Pr}\left(\mathbf{B}^{c}\left(k\right)=\left[\begin{array}{cc}
1 & 0\\
0 & 0
\end{array}\right]\right)=\frac{1}{2}$ and $\mathrm{Pr}\left(\mathbf{B}^{c}\left(k\right)=\left[\begin{array}{cc}
0 & 0\\
0 & 1
\end{array}\right]\right)=\frac{1}{2},\forall k\in\mathbb{Z}_{0}^{+}.$ In this case, the system suffers from almost surely uncontrollability.
Due to the controllers' multiplicative randomness of $\left\{ \mathbf{H}_{j,i}^{c}\left(k\right)\right\} _{\left(j,i\right)\in\mathcal{N}^{c}\times\mathcal{N}^{c}}$,
the pair $\left(\mathbf{A},\mathbf{B}^{c}\left(k\right)\right)$ may
not be controllable for every time slot $k\in\mathbb{Z}_{0}^{+}$.
In fact, the controllability of $\left(\mathbf{A},\mathbf{B}^{c}\left(k\right)\right)$
will depend on the realizations of $\left\{ \mathbf{H}_{j,i}^{c}\left(k\right)\right\} _{\left(j,i\right)\in\mathcal{N}^{c}\times\mathcal{N}^{c}}$,
and $\left(\mathbf{A},\mathbf{B}^{c}\left(k\right)\right)$ can be
either intermittent controllable, i.e., controllable at some time
slots and uncontrollable at some other slots, or even almost surely
uncontrollable at every time slots.
\end{remrk}

We consider the following closed-loop information structure at the
controllers and attackers.

\begin{assumption}
\textit{(Closed-loop Information Structure)}\emph{\label{Assumption: Info-Structure}}
We assume the controllers and attackers have perfect knowledge of
the current and the entire past realizations of the controllers' aggregated
CSI and the plant states. The CSI of the attackers is assumed unknown.
Specifically, the information set available to each controller and
each attacker at time $k\in\mathbb{Z}_{0}^{+}$ is $\mathcal{I}\left(k\right)=\left\{ \mathbf{x}\left(\widetilde{k}\right),\left[\sum_{i\in\mathcal{N}^{c}}\mathbf{B}_{i}\mathbf{H}_{j,i}^{c}\left(\widetilde{k}\right)\right]_{j\in\mathcal{N}^{c}}\right\} _{\widetilde{k}\in\left[0,k\right]}.$
\end{assumption}

In practice, Assumption \ref{Assumption: Info-Structure} can be achieved
by existing pilot-based channel estimation approach in wireless communications
{[}34{]}, {[}39{]}. Specifically, at the beginning of each time-slot
$k\in\mathbb{Z}_{0}^{+}$, each actuator first broadcasts a weighted
common pilot $\mathbf{B}_{i}^{T}\mathbf{C}\in\mathbb{R}^{N_{r}\times C}$,
where $\mathbf{C}\in\mathbb{R}^{S\times C}$ is a full row-rank common
pilot symbol that is known only to the controllers and is unknown
to the attackers. The received aggregated pilot symbols at the $j$-th
controller is $\mathbf{Y}_{j}^{c}\left(k\right)=\sum_{i\in\mathcal{N}^{c}}\left(\mathbf{H}_{j,i}^{c}\left(k\right)\right)^{T}\mathbf{B}_{i}^{T}\mathbf{C}+\mathbf{N}_{j}^{c}\left(k\right)$,
where $\mathbf{N}_{j}^{c}\left(k\right)\in\mathbb{R}^{N_{t}^{c}\times C}$
is the additive channel noise. The $j$-th controller obtains its
aggregated CSI $\sum_{i\in\mathcal{N}^{c}}\mathbf{B}_{i}\mathbf{H}_{j,i}^{c}\left(k\right)$
based on $\mathbf{Y}_{j}^{c}\left(k\right)$ and $\mathbf{C}$ using
the least square channel estimation approach {[}39{]}. Each controller
then shares its aggregated CSI with the other controllers and obtains
the global aggregated CSI $\left[\sum_{i\in\mathcal{N}^{c}}\mathbf{B}_{i}\mathbf{H}_{j,i}^{c}\left(k\right)\right]_{j\in\mathcal{N}^{c}}$
in Assumption \ref{Assumption: Info-Structure}. The attackers can
also receive the broadcasted pilots from the actuators. However, each
attacker cannot obtain its aggregated CSI because the common pilot
$\mathbf{C}$ is unknown to the attackers. Nevertheless, the attackers
can still monitor the sharing of aggregated CSI by the controllers
and gain knowledge about the global aggregated CSI $\left[\sum_{i\in\mathcal{N}^{c}}\mathbf{B}_{i}\mathbf{H}_{j,i}^{c}\left(k\right)\right]_{j\in\mathcal{N}^{c}}$. 

For each controller $j\in\mathcal{N}^{c}$, a control policy $\boldsymbol{\mu}_{j}^{c}$
is a sequence of mappings $\boldsymbol{\mu}_{j}^{c}=\left\{ \boldsymbol{\mu}_{j}^{c}\left(\cdot,k\right)\right\} _{k\in\mathbb{Z}_{0}^{+}},$
where $\boldsymbol{\mu}_{j}^{c}\left(\cdot,k\right):\mathcal{I}\left(k\right)\rightarrow\mathbb{R}^{N_{t}^{c}\times1}$
and $\mathbf{u}_{j}^{c}\left(k\right)=\boldsymbol{\mu}_{j}^{c}\left(\mathcal{I}\left(k\right),k\right)$.
Similarly, for each attacker $l\in\mathcal{N}^{a}$, an attack policy
$\boldsymbol{\mu}_{l}^{a}$ is $\boldsymbol{\mu}_{l}^{a}=\left\{ \boldsymbol{\mu}_{l}^{a}\left(\cdot,k\right)\right\} _{k\in\mathbb{Z}_{0}^{+}}$
with $\boldsymbol{\mu}_{l}^{a}\left(\cdot,k\right):\mathcal{I}\left(k\right)\rightarrow\mathbb{R}^{N_{t}^{a}\times1}$
and $\mathbf{u}_{l}^{a}\left(k\right)=\boldsymbol{\mu}_{l}^{a}\left(\mathcal{I}\left(k\right),k\right)$.
Let $\mathcal{U}^{c}$ and $\mathcal{U}^{a}$ denote the set of all
such Borel measurable control policies and attack policies, respectively.

We define an infinite-horizon average quadratic cost $J$ as
\begin{align}
 & J\left(\left\{ \boldsymbol{\mu}_{j}^{c}\right\} _{j\in\mathcal{N}^{c}},\text{\ensuremath{\left\{  \boldsymbol{\mu}_{l}^{a}\right\}  _{l\in\mathcal{N}^{a}}}}\right)\\
 & =\lim_{K\rightarrow\infty}\frac{1}{K}\sum_{k=0}^{K-1}\mathbb{E}\Bigg[\left\Vert \mathbf{x}\left(k+1\right)\right\Vert _{\mathbf{Q}}^{2}+\sum_{j\in\mathcal{N}^{c}}\left\Vert \mathbf{u}_{j}^{c}\left(k\right)\right\Vert _{\mathbf{R}_{j}^{c}}^{2}\nonumber \\
 & -\sum_{l\in\mathcal{N}^{a}}\left\Vert \mathbf{u}_{l}^{a}\left(k\right)\right\Vert _{\mathbf{R}_{l}^{a}}^{2}\Bigg]=\lim_{K\rightarrow\infty}\frac{1}{K}\sum_{k=0}^{K-1}\mathbb{E}\Big[\left\Vert \mathbf{x}\left(k+1\right)\right\Vert _{\mathbf{Q}}^{2}\nonumber \\
 & +\left\Vert \mathbf{u}^{c}\left(k\right)\right\Vert _{\mathbf{R}^{c}}^{2}-\left\Vert \mathbf{u}^{a}\left(k\right)\right\Vert _{\mathbf{R}^{a}}^{2}\Big],\nonumber 
\end{align}
where $\mathbf{Q}\in\mathbb{S}_{+}^{S}$, $\mathbf{R}_{j}^{c}\in\mathbb{\mathbb{S}}_{+}^{N_{t}^{c}},\forall j\in\mathcal{N}^{c},$
and $\mathbf{R}_{l}^{a}\in\mathbb{\mathbb{S}}_{+}^{N_{t}^{a}},\forall l\in\mathcal{N}^{a},$
are the positive definite weight matrices; and $\mathbf{R}^{c}=\mathrm{diag}\left(\left[\mathbf{R}_{j}^{c}\right]_{j\in\mathcal{N}^{c}}\right),\ \mathbf{R}^{a}=\mathrm{diag}\left(\left[\mathbf{R}_{l}^{a}\right]_{l\in\mathcal{N}^{a}}\right).$ 

In a zero-sum stochastic linear-quadratic difference game, the objective
of the controllers is to minimize $J$, while the attackers aim to
maximize it. There are two important notions in zero-sum game theory,
namely the upper and lower values, which are summarized in Definition
\ref{definition-game-value}.

\begin{definitn}
\textsl{(}\textit{Upper Value, Lower Value and Value of the game}\textsl{)}\label{definition-game-value}
The upper value $\overline{J}$ and the lower value $\underline{J}$
of the LQ zero-sum game are defined as {\small{}
\begin{align}
 & \overline{J}=\label{eq:upper-value}\\
 & \inf_{\left\{ \boldsymbol{\mu}_{l}^{c}\right\} _{j\in\mathcal{N}^{c}}\in\left(\mathcal{U}^{c}\right)^{\left|\mathcal{N}^{c}\right|}}\sup_{\left\{ \boldsymbol{\mu}_{l}^{a}\right\} _{l\in\mathcal{N}^{a}}\in\left(\mathcal{U}^{a}\right)^{\left|\mathcal{N}^{a}\right|}}J\left(\left\{ \boldsymbol{\mu}_{l}^{c}\right\} _{j\in\mathcal{N}^{c}},\text{\ensuremath{\left\{  \boldsymbol{\mu}_{l}^{a}\right\}  _{l\in\mathcal{N}^{a}}}}\right),\nonumber 
\end{align}
}and{\small{}
\begin{align}
 & \underline{J}=\label{eq:lower-value}\\
 & \sup_{\left\{ \boldsymbol{\mu}_{l}^{a}\right\} _{l\in\mathcal{N}^{a}}\in\left(\mathcal{U}^{a}\right)^{\left|\mathcal{N}^{a}\right|}}\inf_{\left\{ \boldsymbol{\mu}_{l}^{c}\right\} _{j\in\mathcal{N}^{c}}\in\left(\mathcal{U}^{c}\right)^{\left|\mathcal{N}^{c}\right|}}J\left(\left\{ \boldsymbol{\mu}_{l}^{c}\right\} _{j\in\mathcal{N}^{c}},\text{\ensuremath{\left\{  \boldsymbol{\mu}_{l}^{a}\right\}  _{l\in\mathcal{N}^{a}}}}\right),\nonumber 
\end{align}
}respectively. The value $J^{*}$ of the LQ zero-sum game exists if
and only if the upper value $\overline{J}$ and the lower value $\underline{J}$
are equal and is given by 
\begin{align}
 & J^{*}=\overline{J}=\underline{J}.
\end{align}

\end{definitn}

The Nash equilibrium is another notion of paramount importance to
an LQ zero-sum game. If the LQ zero-sum game has a Nash equilibrium,
then the value of the game is the Nash equilibrium cost, i.e., the
value of the quadratic cost function in such a Nash equilibrium. It
is important to note that the value of the LQ zero-sum game, if it
exists, is unique. This is becasue the value of the LQ zero-sum game
characterizes the global minimum loss that the controllers can secure
against the attackers\textquoteright{} best possible strategies. Therefore,
if there exist multiple Nash equilibria, the value of the quadratic
cost function in each Nash equilibrium is identical and equal to the
value of the game. The formal definition of the Nash equilibrium (i.e.,
the saddle point policies) as summarized in the following Definition
\ref{definition-NE-1}. Based on the previous discussion, it is clear
that $J\left(\left\{ \left(\boldsymbol{\mu}_{j}^{c}\right)^{*}\right\} _{j\in\mathcal{N}^{c}},\left\{ \left(\boldsymbol{\mu}_{l}^{a}\right)^{*}\right\} _{l\in\mathcal{N}^{c}}\right)$
in Definition \ref{definition-NE-1} is the value of the LQ zero-sum
difference game. 

\begin{definitn}
\textsl{(}\textit{Nash Equilibrium}\textsl{)}\label{definition-NE-1}\textsl{
}A policy pair {\small{}$\left\{ \left\{ \left(\boldsymbol{\mu}_{j}^{c}\right)^{*}\right\} _{j\in\mathcal{N}^{c}},\left\{ \left(\boldsymbol{\mu}_{l}^{a}\right)^{*}\right\} _{l\in\mathcal{N}^{a}}\right\} $,
$\left(\boldsymbol{\mu}_{j}^{c}\right)^{*}$ $\in\mathcal{U}^{c},\forall j\in\mathcal{N}^{c},\left(\boldsymbol{\mu}_{l}^{a}\right)^{*}$
$\in\mathcal{U}^{a},\forall l\in\mathcal{N}^{a},$} constitutes a
Nash equilibrium, or equivalently a saddle-point policy, in the class
of admissible policy pairs $\left(\mathcal{U}^{c}\right)^{\left|\mathcal{N}^{c}\right|}\times\left(\mathcal{U}^{a}\right)^{\left|\mathcal{N}^{a}\right|}$
if the following two conditions are satisfied simultaneously.
\end{definitn}

\begin{itemize}
\item $J\left(\left\{ \left(\boldsymbol{\mu}_{j}^{c}\right)^{*}\right\} _{j\in\mathcal{N}^{c}},\left\{ \left(\boldsymbol{\mu}_{l}^{a}\right)^{*}\right\} _{l\in\mathcal{N}^{c}}\right)\leq$\\$J\left(\left\{ \boldsymbol{\mu}_{\widetilde{j}}^{c}\right\} _{\widetilde{j}\in\mathcal{N}_{1}},\left\{ \left(\boldsymbol{\mu}_{j}^{c}\right)^{*}\right\} _{j\in\left(\mathcal{N}^{c}\setminus\mathcal{N}_{1}\right)},\left\{ \left(\boldsymbol{\mu}_{l}^{a}\right)^{*}\right\} _{l\in\mathcal{N}^{a}}\right)$,
$\forall\widetilde{j}\in\mathcal{N}_{1},\forall\mathcal{N}_{1}\subseteq\mathcal{N}^{c},\forall\boldsymbol{\mu}_{\widetilde{j}}^{c}\in\mathcal{U}^{c}$;
\item $J\left(\left\{ \left(\boldsymbol{\mu}_{j}^{c}\right)^{*}\right\} _{j\in\mathcal{N}^{c}},\left\{ \left(\boldsymbol{\mu}_{l}^{a}\right)^{*}\right\} _{l\in\mathcal{N}^{c}}\right)\geq$\\$J\left(\left\{ \left(\boldsymbol{\mu}_{j}^{c}\right)^{*}\right\} _{j\in\mathcal{N}^{c}},\left\{ \boldsymbol{\mu}_{\widetilde{l}}^{a}\right\} _{\widetilde{l}\in\mathcal{N}_{2}},\left\{ \left(\boldsymbol{\mu}_{l}^{a}\right)^{*}\right\} _{l\in\left(\mathcal{N}^{a}\setminus\mathcal{N}_{2}\right)}\right),$
$\forall\widetilde{l}\in\mathcal{N}_{2},\forall\mathcal{N}_{2}\subseteq\mathcal{N}^{a},\forall\boldsymbol{\mu}_{\widetilde{l}}^{a}\in\mathcal{U}^{a}.$
\end{itemize}

For the finite-time horizon zero-sum game, Lemma \ref{Lemma: J_K representations}
provides a sufficient and necessary condition for the existence and
uniqueness of the Nash equilibrium. For the infinite-time horizon
case, Theorem \ref{Thm: General case- suff and necessary} provides
a sufficient and necessary condition for the existence of a Nash equilibrium,
which may not be unique.

\section{Conditions for the Existence of Game Value}

In this section, we aim to derive various conditions that can guarantee
the existence of the game value. 

\subsection{Modified GARE}

We first define an MGARE $\mathbf{P}=f\left(\mathbf{P}\right)$ for
$\mathbf{P}\in\mathbb{S}_{+}^{S}$, where $f\left(\mathbf{P}\right)$
is given by
\begin{align}
f\left(\mathbf{P}\right)= & \mathbf{A}^{T}\mathbb{E}\Big[\mathbf{P}-\mathbf{P}\mathbf{B}\left(k\right)\Big[\mathbf{R}\left(\mathbf{P}\right)+\mathbf{B}^{T}\left(k\right)\mathbf{P}\mathbf{B}\left(k\right)\Big]^{-1}\label{eq: general-GARE}\\
 & \cdot\mathbf{B}^{T}\left(k\right)\mathbf{P}\Big]\mathbf{A}+\mathbf{Q},\nonumber 
\end{align}
with $\mathbf{B}\left(k\right)=\left[\begin{array}{cc}
\mathbf{B}^{c}\left(k\right) & \mathbb{E}\left[\mathbf{B}^{a}\left(k\right)\right]\end{array}\right],$ $\mathbf{R}\left(\mathbf{P}\right)=\mathrm{diag}(\mathbf{R}^{c},-\mathbf{\widetilde{R}}^{a}\left(\mathbf{P}\right)),$
and $\mathbf{\widetilde{R}}^{a}\left(\mathbf{P}\right)=\mathbf{R}^{a}-\mathrm{cov}(\mathbf{B}^{a}\left(k\right)\mathbf{P}^{\frac{1}{2}}).$ 

The expectation in (\ref{eq: general-GARE}) is taken w.r.t. controllers'
multiplicative randomness $\left\{ \mathbf{H}_{j,i}^{c}\left(k\right)\right\} _{\left(j,i\right)\in\mathcal{N}^{c}\times\mathcal{N}^{c}}$.
The expectations in $\mathbf{B}\left(k\right)$ and $\mathbf{\widetilde{R}}^{a}\left(\mathbf{P}\right)$
are taken w.r.t. attackers' multiplicative randomness $\left\{ \mathbf{H}_{l,i}^{a}\left(k\right)\right\} _{\left(l,i\right)\in\mathcal{N}^{a}\times\mathcal{N}^{c}}$. 

For each $\mathbf{R}^{a}\in\mathbb{S}_{+}^{N_{t}^{a}\left|\mathcal{N}^{a}\right|\times N_{t}^{a}\left|\mathcal{N}^{a}\right|},$
we further define two subsets of positive definite matrices as follows.
\begin{alignat}{1}
 & \mathcal{P}_{\mathbf{R}^{a}}\triangleq\left\{ \mathbf{P}\in\mathbb{S}_{+}^{S}:\mathbf{P}=f\left(\mathbf{P}\right)\right\} ;\\
 & \mathcal{R}_{\mathbf{R}^{a}}\triangleq\bigg\{\mathbf{P}\in\mathbb{S}_{+}^{S}:\mathbf{R}^{a}\succ\mathbb{E}\left[\left(\mathbf{B}^{a}\left(k\right)\right)^{T}\mathbf{P}\mathbf{B}^{a}\left(k\right)\right]\bigg\}.\label{eq: Definition R^a-set}
\end{alignat}
The set $\mathcal{P}_{\mathbf{R}^{a}}$ contains all the positive
definite solutions to the MGARE $\mathbf{P}=f\left(\mathbf{P}\right)$,
and $\mathcal{P}_{\mathbf{R}^{a}}=\textrm{Ø}$ if the MGARE has no
positive definite solution.

We denote $f^{k+1}\left(\mathbf{Q}\right)=f\left(f^{k}\left(\mathbf{Q}\right)\right),\forall k\in\mathbb{Z}_{0}^{+},$
with $f^{0}\left(\mathbf{Q}\right)=\mathbf{Q}.$ Several key properties
of the sets $\mathcal{P}_{\mathbf{R}^{a}}$ and $\mathcal{R}_{\mathbf{R}^{a}}$
are summarized below, which will be frequently used in our later analysis
on the existence of the game value. 

\begin{lemma}
\textit{(Monotonicity)}\emph{\label{Lemma: Monotonicity f(.)}} $\forall\mathbf{P}_{1},\mathbf{P}_{2}\in\mathcal{R}_{\mathbf{R}^{a}}$,
$\mathbf{P}_{1}\preceq\mathbf{P}_{2}$ implies $f\left(\mathbf{P}_{1}\right)\preceq f\left(\mathbf{P}_{2}\right)$.
\end{lemma}

\begin{lemma}
\textit{(Existence of $f^{\infty}\left(\mathbf{Q}\right)$)}\emph{\label{Lemma: Minimal-Solution-P*}}
If $\mathcal{P}_{\mathbf{R}^{a}}\cap\mathcal{R}_{\mathbf{R}^{a}}\neq\textrm{Ø}$,
$f^{\infty}\left(\mathbf{Q}\right)=\lim_{k\rightarrow\infty}f^{k}\left(\mathbf{Q}\right)$
exists and $f^{\infty}\left(\mathbf{Q}\right)\in\mathcal{P}_{\mathbf{R}^{a}}\cap\mathcal{R}_{\mathbf{R}^{a}}$.
\end{lemma}

\begin{lemma}
\textit{(Minimum Element)}\emph{\label{Lemma: P*=00003Df^infty(Q)}}
If $\mathcal{P}_{\mathbf{R}^{a}}\cap\mathcal{R}_{\mathbf{R}^{a}}\neq\textrm{Ø}$,
then $f^{\infty}\left(\mathbf{Q}\right)$ is the minimum element of
$\mathcal{P}_{\mathbf{R}^{a}}\cap\mathcal{R}_{\mathbf{R}^{a}}$, i.e.,
$f^{\infty}\left(\mathbf{Q}\right)\preceq\mathbf{\mathbf{P}}$ for
all $\mathbf{\mathbf{P}}\in\mathcal{P}_{\mathbf{R}^{a}}\cap\mathcal{R}_{\mathbf{R}^{a}}.$
\end{lemma}

\subsection{Finite-Time Horizon Average Cost}

In this subsection, we study the structural properties of the finite-time
horizon average cost, which facilitates our later analysis of the
infinite-time horizon problem by letting the length of the horizon
approach infinity. 

The average quadratic cost $J_{K}$ over a finite time horizon of
length $K\in\mathbb{Z}^{+}$ is defined as
\begin{align}
 & J_{K}=\frac{1}{K}\sum_{k=0}^{K-1}\mathbb{E}\left[\left\Vert \mathbf{x}\left(k+1\right)\right\Vert _{\mathbf{Q}}^{2}+\left\Vert \mathbf{u}^{c}\left(k\right)\right\Vert _{\mathbf{R}^{c}}^{2}-\left\Vert \mathbf{u}^{a}\left(k\right)\right\Vert _{\mathbf{R}^{a}}^{2}\right].\label{eq:def-J_k}
\end{align}

We introduce a list of time varying matrix functions as follows, which
would facilitate the later analysis of $J_{K}$. Specifically, $\forall\mathbf{P}\in\mathcal{R}_{\mathbf{R}^{a}},\forall0\leq k\leq K-1$,
we define
\begin{align}
 & \Phi_{1}\left(\mathbf{\mathbf{P}};k\right)=\left(\mathbf{B}^{c}\left(k\right)\right)^{T}\mathbf{\mathbf{P}}\mathbf{B}^{c}\left(k\right)+\mathbf{R}^{c},\\
 & \Phi_{2}\left(\mathbf{\mathbf{P}};k\right)=\left(\mathbf{B}^{c}\left(k\right)\right)^{T}\mathbf{\mathbf{P}}\mathbb{E}\left[\mathbf{B}^{a}\left(k\right)\right],\\
 & \Phi_{3}\left(\mathbf{\mathbf{P}}\right)=\mathbf{R}^{a}-\mathbb{E}\left[\left(\mathbf{B}^{a}\left(k\right)\right)^{T}\mathbf{P}\mathbf{B}^{a}\left(k\right)\right],\\
 & \Phi\left(\mathbf{P};k\right)=\left[\begin{array}{cc}
\Phi_{1}\left(\mathbf{\mathbf{P}};k\right) & \Phi_{2}\left(\mathbf{\mathbf{P}};k\right)\\
\Phi_{2}^{T}\left(\mathbf{\mathbf{P}};k\right) & -\Phi_{3}\left(\mathbf{\mathbf{P}}\right)
\end{array}\right],\label{eq:big-phi-MGARE}\\
 & \left(\mathbf{u}^{c}\right)^{*}\left(\mathbf{P};k\right)=-\left(\Phi^{-1}\left(\mathbf{P};k\right)\mathbf{B}^{T}\left(k\right)\mathbf{P}\mathbf{A}\mathbf{x}\left(k\right)\right)_{\left(1:N_{t}^{c}\left|\mathcal{N}^{c}\right|;:\right)},\label{eq:def-u_c*}\\
 & \left(\mathbf{u}^{a}\right)^{*}\left(\mathbf{P};k\right)=-\big(\Phi^{-1}\left(\mathbf{P};k\right)\mathbf{B}^{T}\left(k\right)\mathbf{P}\label{eq:def-u_a*}\\
 & \cdot\mathbf{A}\mathbf{x}\left(k\right)\big)_{\left(\left(N_{t}^{c}\left|\mathcal{N}^{c}\right|+1\right):\left(N_{t}^{c}\left|\mathcal{N}^{c}\right|+N_{t}^{a}\left|\mathcal{N}^{a}\right|\right);:\right)},\nonumber \\
 & \mathbf{u}_{\Delta}^{c}\left(\mathbf{P};k\right)=\begin{array}{c}
\mathbf{u}^{c}\left(k\right)\end{array}-\left(\mathbf{u}^{c}\right)^{*}\left(\mathbf{P};k\right),\label{eq:delta-u-c}\\
 & \mathbf{u}_{\Delta}^{a}\left(\mathbf{\mathbf{P}};k\right)=\mathbf{u}^{a}\left(k\right)-\left(\mathbf{u}^{a}\right)^{*}\left(\mathbf{P};k\right).\label{eq:delta-u-a}
\end{align}

It is clear that $\Phi\left(\mathbf{P};k\right)$ in (\ref{eq:big-phi-MGARE})
is a representation of the inverse term in the MGARE (\ref{eq: general-GARE}),
i.e., $\Phi\left(\mathbf{P};k\right)=\mathbf{R}\left(\mathbf{P}\right)+\mathbf{B}^{T}\left(k\right)\mathbf{P}\mathbf{B}\left(k\right).$
Besides, $\left(\mathbf{u}^{c}\right)^{*}\left(\mathbf{P};k\right)$
in (\ref{eq:def-u_c*}) and $\left(\mathbf{u}^{a}\right)^{*}\left(\mathbf{P};k\right)$
in (\ref{eq:def-u_a*}) have the same form as the steady state controls
and attacks defined in {[}11{]}, {[}20{]}, respectively. 

The inverse $\Phi^{-1}\left(\mathbf{P};k\right)$ in (\ref{eq:def-u_c*})
and (\ref{eq:def-u_a*}) is well-defined because we restrict $\mathbf{P}\in\mathcal{R}_{\mathbf{R}^{a}}$.
We now have the following lemma  on the representations of $J_{K}.$

\begin{lemma}
\textit{\emph{(}}\textit{Representations of }\textit{\emph{$J_{K}$)}}\emph{\label{Lemma: J_K representations}}
For any given $K\in\mathbb{Z}^{+}$, if $\mathbf{R}^{a}$ is chosen
such that the following condition is satisfied, 
\begin{align}
 & \mathbf{R}^{a}\succ\mathbb{E}\left[\left(\mathbf{B}^{a}\left(k\right)\right)^{T}f^{K-1}\left(\mathbf{Q}\right)\mathbf{B}^{a}\left(k\right)\right],\label{eq:Ra-Jk-Condi}
\end{align}
then the finite time horizon average cost $J_{K}$ admits the following
two different representations:
\end{lemma}

\begin{itemize}
\item \textit{Upper Value Representation}: $J_{K}=\overline{J}_{K}$, where
$\overline{J}_{K}$ is given by (\ref{eq:J_K-Upper}).
\item \textit{Lower Value Representation}: $J_{K}=\underline{J}_{K}$, where
$\underline{J}_{K}$ is given by (\ref{eq:J_K-Lower}).
\end{itemize}
\begin{figure*}
{\footnotesize{}
\begin{align}
 & \overline{J}_{K}=\frac{1}{K}\left\Vert \mathbf{x}\left(0\right)\right\Vert _{\left(f^{K}\left(\mathbf{Q}\right)-\mathbf{Q}\right)}^{2}+\mathrm{Tr}\left(\frac{1}{K}\sum_{k=0}^{K-1}\left(f^{k}\left(\mathbf{Q}\right)\right)\left(\mathbf{\mathbf{W}}+\sum_{i\in\mathcal{N}^{c}}\mathbf{B}_{i}^{T}\mathbf{V}\mathbf{B}_{i}\right)\right)\label{eq:J_K-Upper}\\
 & +\frac{1}{K}\sum_{k=0}^{K-1}\mathbb{E}\left[\left\Vert \mathbf{u}_{\Delta}^{c}\left(f^{K-k-1}\left(\mathbf{Q}\right);k\right)\right\Vert _{\Phi\left(f^{K-k-1}\left(\mathbf{Q}\right);k\right)\big/\left(-\Phi_{3}\left(f^{K-k-1}\left(\mathbf{Q}\right)\right)\right)}^{2}\right]\nonumber \\
 & -\frac{1}{K}\sum_{k=0}^{K-1}\mathbb{E}\left[\left\Vert \mathbf{u}_{\Delta}^{a}\left(f^{K-k-1}\left(\mathbf{Q}\right);k\right)+\Phi_{3}^{-1}\left(f^{K-k-1}\left(\mathbf{Q}\right)\right)\Phi_{2}\left(f^{K-k-1}\left(\mathbf{Q}\right);k\right)\mathbf{u}_{\Delta}^{c}\left(f^{K-k-1}\left(\mathbf{Q}\right);k\right)\right\Vert _{\Phi_{3}\left(f^{K-k-1}\left(\mathbf{Q}\right)\right)}^{2}\right];\nonumber \\
 & \underline{J}_{K}=\frac{1}{K}\left\Vert \mathbf{x}\left(0\right)\right\Vert _{\left(f^{K}\left(\mathbf{Q}\right)-\mathbf{Q}\right)}^{2}+\mathrm{Tr}\left(\frac{1}{K}\sum_{k=0}^{K-1}\left(f^{k}\left(\mathbf{Q}\right)\right)\left(\mathbf{\mathbf{W}}+\sum_{i\in\mathcal{N}^{c}}\mathbf{B}_{i}^{T}\mathbf{V}\mathbf{B}_{i}\right)\right)\label{eq:J_K-Lower}\\
 & -\frac{1}{K}\sum_{k=0}^{K-1}\mathbb{E}\left[\left\Vert \mathbf{u}_{\Delta}^{a}\left(f^{K-k-1}\left(\mathbf{Q}\right);k\right)\right\Vert _{-\left(\Phi\left(f^{K-k-1}\left(\mathbf{Q}\right);k\right)\big/\left(\Phi_{1}\left(f^{K-k-1}\left(\mathbf{Q}\right);k\right)\right)\right)}^{2}\right]\nonumber \\
 & +\frac{1}{K}\sum_{k=0}^{K-1}\mathbb{E}\left[\left\Vert \mathbf{u}_{\Delta}^{c}\left(f^{K-k-1}\left(\mathbf{Q}\right);k\right)+\Phi_{1}^{-1}\left(f^{K-k-1}\left(\mathbf{Q}\right);k\right)\Phi_{2}\left(f^{K-k-1}\left(\mathbf{Q}\right);k\right)\mathbf{u}_{\Delta}^{a}\left(f^{K-k-1}\left(\mathbf{Q}\right);k\right)\right\Vert _{\Phi_{1}\left(f^{K-k-1}\left(\mathbf{Q}\right);k\right)}^{2}\right]\nonumber 
\end{align}
}{\footnotesize\par}

\hrulefill
\end{figure*}
 
\begin{IEEEproof}
Please see Appendix \ref{subsec:Proof-of-Lemma-J_k-representations}.
\end{IEEEproof}

Lemma \ref{Lemma: J_K representations} delivers several key insights
on the structural properties of the finite time horizon cost $J_{K}$.
Based on $\overline{J}_{K}$ in (\ref{eq:J_K-Upper}), we have
\begin{align}
 & \min_{\left\{ \mathbf{u}_{\Delta}^{c}\left(f^{K-k-1}\left(\mathbf{Q}\right);k\right)\right\} _{k\in\left[0,K-1\right]}}\max_{\left\{ \left(\mathbf{u}_{\Delta}^{a}\right)\left(f^{K-k-1}\left(\mathbf{Q}\right);k\right)\right\} _{k\in\left[0,K-1\right]}}\overline{J}_{K}\label{eq:min-max-J_K}\\
 & =\frac{1}{K}\left\Vert \mathbf{x}\left(0\right)\right\Vert _{\left(f^{K}\left(\mathbf{Q}\right)-\mathbf{Q}\right)}^{2}\nonumber \\
 & +\mathrm{Tr}\left(\frac{1}{K}\sum_{k=0}^{K-1}\left(f^{k}\left(\mathbf{Q}\right)\right)\left(\mathbf{\mathbf{W}}+\sum_{i\in\mathcal{N}^{c}}\mathbf{B}_{i}^{T}\mathbf{V}\mathbf{B}_{i}\right)\right),\nonumber 
\end{align}
where the minimax value of $\overline{J}_{K}$ in (\ref{eq:min-max-J_K})
is achieved if $\mathbf{u}_{\Delta}^{c}\left(f^{K-k-1}\left(\mathbf{Q}\right);k\right)=\mathbf{0}$
and 
\begin{align}
 & \left(\mathbf{u}_{\Delta}^{a}\right)\left(f^{K-k-1}\left(\mathbf{Q}\right);k\right)=-\Phi_{3}^{-1}\left(f^{K-k-1}\left(\mathbf{Q}\right)\right)\\
 & \cdot\Phi_{2}\left(f^{K-k-1}\left(\mathbf{Q}\right);k\right)\mathbf{u}_{\Delta}^{c}\left(f^{K-k-1}\left(\mathbf{Q}\right);k\right),\nonumber 
\end{align}
$\forall0\leq k\leq K-1.$ Based on (\ref{eq:delta-u-c}) and (\ref{eq:delta-u-a}),
this is equivalent to $\mathbf{u}^{c}\left(k\right)=\left(\mathbf{u}^{c}\right)^{*}\left(\mathbf{P};k\right)$
and $\mathbf{u}^{a}\left(k\right)=\left(\mathbf{u}^{a}\right)^{*}\left(\mathbf{P};k\right),$
$\forall0\leq k\leq K-1.$ Similarly, from $\underline{J}_{K}$ in
(\ref{eq:J_K-Lower}), the maximin value of $\underline{J}_{K}$ is
achieved if $\mathbf{u}^{a}\left(k\right)=\left(\mathbf{u}^{a}\right)^{*}\left(\mathbf{P};k\right)$
and $\mathbf{u}^{c}\left(k\right)=\left(\mathbf{u}^{c}\right)^{*}\left(\mathbf{P};k\right),$
$\forall0\leq k\leq K-1.$ This means that, under condition (\ref{eq:Ra-Jk-Condi})
the control and attack policy pair $\left\{ \left(\mathbf{u}^{c}\right)^{*}\left(\mathbf{P};k\right),\left(\mathbf{u}^{a}\right)^{*}\left(\mathbf{P};k\right)\right\} $
constitutes a Nash equilibrium for the finite time horizon zero-sum
game. 

\subsection{General Sufficient and Necessary Condition for Existence of Game
Value}

We are now ready to give a general sufficient and necessary condition
for the existence of the value of the zero-sum game, which is summarized
in the following theorem.

\begin{thm}
\textit{(General Sufficient and Necessary Condition)}\emph{\label{Thm: General case- suff and necessary}}
The value of the infinite-time horizon zero-sum game exists if and
only if 
\begin{align}
 & \mathcal{P}_{\mathbf{R}^{a}}\cap\mathcal{R}_{\mathbf{R}^{a}}\neq\textrm{Ø}.\label{eq: general-GARE condition}
\end{align}
Furthermore, under condition (\ref{eq: general-GARE condition}), 
\end{thm}

\begin{itemize}
\item $\mathbf{P}^{*}\triangleq f^{\infty}\left(\mathbf{Q}\right)$ is the
minimal positive definite solution to the MGARE (\ref{eq: general-GARE}),
i.e., every positive definite solution $\mathbf{P}$ to the MGARE
(\ref{eq: general-GARE}) satisfies $\mathbf{P}\succeq\mathbf{P}^{*}.$
\item the value of the game is given by
\begin{align}
 & J^{*}=\mathrm{Tr}\left(\mathbf{P}^{*}\left(\mathbf{\mathbf{W}}+\sum_{i\in\mathcal{N}^{c}}\mathbf{B}_{i}^{T}\mathbf{V}\mathbf{B}_{i}\right)\right).\label{eq:game-value-static-1-1}
\end{align}
\end{itemize}
\begin{IEEEproof}
Please see Appendix \ref{subsec:Proof-of-Thm General Suff and Necessary-1}.
\end{IEEEproof}

Condition (\ref{eq: general-GARE condition}) in Theorem \ref{Thm: General case- suff and necessary}
guarantees the existence of a Nash equilibrium for the infinite-time
horizon case. There may exist multiple Nash equilibria, but the value
of the game is unique and is given by (\ref{eq:game-value-static-1-1}).
This is consistent with the fact, pointed out in {[}42{]}-{[}46{]},
that in general infinite-time horizon linear quadratic dynamic games,
the Nash equilibrium may not be unique. It is worth noting that condition
(\ref{eq:Ra-Jk-Condi}) in Lemma \ref{Lemma: J_K representations}
do imply that the finite-time horizon game has a unique Nash equilibrium. 
\begin{remrk}
\textsl{(}\textit{Verifiability Challenge: Intermittent Controllability
and Uncontrollability}\textsl{)} \emph{\label{Remark: Unctrl+intermitent Ctrl}}The
general sufficient and necessary condition (\ref{eq: general-GARE condition})
conveys two essential messages. First, the MGARE (\ref{eq: general-GARE})
should have a positive definite solution. And then, such a positive
definite solution should be an element of the set $\mathcal{R}_{\mathbf{R}^{a}}$.
In existing literature, where linear time-invariant systems without
randomness are considered, the existence of solution to the GARE associated
with the zero-sum game relies heavily on the controllability of the
pair $\left(\mathbf{A},\mathbf{B}^{c}\right)$. As shown in {[}11{]},
{[}20{]}, when $\mathbf{R}^{a}$ is chosen sufficiently large, if
$\left(\mathbf{A},\mathbf{B}^{c}\right)$ is controllable, then the
GARE admits a positive definite solution. However, in this paper,
we are considering a time-varying stochastic system, where the pair
$\left(\mathbf{A},\mathbf{B}^{c}\left(k\right)\right)$ may not be
controllable at every time slot. In this case, even if $\mathbf{R}^{a}$
is chosen sufficiently large, the MAGRE (\ref{eq: general-GARE})
may still not possess a positive definite solution. Existing literature
has not investigated the conditions for the existence of a positive
definite solution to the MAGRE (\ref{eq: general-GARE}) in the presence
of intermittent controllability and uncontrollability.
\end{remrk}

\begin{remrk}
\textsl{(}\textit{Verifiability Challenge: Tight Coupling Between
$\mathcal{P}_{\mathbf{R}^{a}}$ and $\mathcal{R}_{\mathbf{R}^{a}}$}\textsl{)}
\emph{\label{Remark: Tight Coupling}}The value of the positive definite
solution $\mathbf{P}^{*}$ to the MAGRE (\ref{eq: general-GARE})
depends on the specific choice of $\mathbf{R}^{a}$. However, the
choice of $\mathbf{R}^{a}$ also depends on the specific value of
$\mathbf{P}^{*}$, i.e., $\mathbf{R}^{a}\succ\mathbb{E}\left[\left(\mathbf{B}^{a}\left(k\right)\right)^{T}\mathbf{P}^{*}\mathbf{B}^{a}\left(k\right)\right]$.
As a result, the sets $\mathcal{P}_{\mathbf{R}^{a}}$ and $\mathcal{R}_{\mathbf{R}^{a}}$
are intrinsically coupled together. Moreover, due to the expectation
w.r.t. the multiplicative randomness of the controllers in the MAGRE
(\ref{eq: general-GARE}), the positive definite solution $\mathbf{P}^{*}\left(\mathbf{R}^{a}\right)$
does not have closed-form characterization. As a result, finding an
$\mathbf{R}^{a}$ that satisfies the implicit concavity constraint
$\mathbf{R}^{a}\succ\mathbb{E}\left[\left(\mathbf{B}^{a}\left(k\right)\right)^{T}\mathbf{P}^{*}\left(\mathbf{R}^{a}\right)\mathbf{B}^{a}\left(k\right)\right]$
is generally difficult. 
\end{remrk}

As outlined in Remarks \ref{Remark: Unctrl+intermitent Ctrl} and
\ref{Remark: Tight Coupling}, the verification of condition (\ref{eq: general-GARE condition})
in Theorem \ref{Thm: General case- suff and necessary} encounters
two significant challenges: intermittent controllability and uncontrollability,
and tight coupling between $\mathcal{P}_{\mathbf{R}^{a}}$ and $\mathcal{R}_{\mathbf{R}^{a}}$.
Nevertheless, we have established a variety of methods to verify condition
(\ref{eq: general-GARE condition}) in Theorem \ref{Thm: General case- suff and necessary}.
Specifically, the following Corollary \ref{Corollary: Refined Suff=000026Necessary}
provides two conditions that are equivalent to condition (\ref{eq: general-GARE condition}).
Lemma \ref{Lemma: P-tilde bound} and Theorem \ref{Thm: g(T)<T} provide
several sufficient conditions for condition (\ref{eq: general-GARE condition})
to hold. Lemma \ref{Lemma: Suff-Condition-Existence of T} and Theorem
\ref{Thm: Verifiable Suff Condition} provide a verifiable sufficient
condition for condition (\ref{eq: general-GARE condition}) to hold.
Corollaries \ref{Coro: Example 1 Suff}, \ref{Coro: Example 2 Suff}
and \ref{Coro: Example 3 Suff} provide verifiable sufficient conditions
for condition (\ref{eq: general-GARE condition}) to hold under the
application setups in Examples 1, 2 and 3, respectively. 

\begin{cor}
\textsl{(}\textit{Equivalent Sufficient and Necessary Condition}\textsl{)}
\emph{\label{Corollary: Refined Suff=000026Necessary}}The general
sufficient and necessary condition (\ref{eq: general-GARE condition})
holds if and only if the following two conditions are satisfied simultaneously.
\end{cor}

\begin{itemize}
\item $f^{k}\left(\mathbf{Q}\right)$ is bounded, $\forall k\in\mathbb{Z}_{0}^{+}.$
\item $f^{k}\left(\mathbf{Q}\right)\in\mathcal{R}_{\mathbf{R}^{a}},\forall k\in\mathbb{Z}_{0}^{+}.$
\end{itemize}
\begin{IEEEproof}
Please see Appendix \ref{subsec:Proof-of-Corollary: Refined Suff=000026Necessary}.
\end{IEEEproof}

\begin{remrk}
\textsl{(}\textit{Verifiability Challenge: Infinite Computational
Complexity}\textsl{)} \emph{\label{Remark: Verifiability Issues}}Although
Corollary \ref{Corollary: Refined Suff=000026Necessary} provides
a closed-form analytical sufficient and necessary condition, it involves
the calculation of $f^{k}\left(\mathbf{Q}\right)$ for all $k\in\mathbb{Z}_{0}^{+}$,
which has infinite computational complexity. As a result, verifying
Lemma \ref{Corollary: Refined Suff=000026Necessary} for practical
systems is generally infeasible. This prompts us to pursue the derivation
of alternative verifiable conditions.
\end{remrk}

\subsection{Verifiable Sufficient Condition for Existence of Game Value}

We first propose a sufficient condition for the non-emptiness of the
set intersection $\mathcal{P}_{\mathbf{R}^{a}}\cap\mathcal{R}_{\mathbf{R}^{a}}.$

\begin{lemma}
\textit{(Sufficient Condition for }$\mathcal{P}_{\mathbf{R}^{a}}\cap\mathcal{R}_{\mathbf{R}^{a}}\neq\textrm{Ø}$\textit{)}
\emph{\label{Lemma: P-tilde bound}} If there exists a matrix pair
$\left(\widetilde{\mathbf{P}},\mathbf{R}^{a}\right)$ that satisfies
the following multivariate nonlinear matrix inequality (MNMI) 
\begin{align}
 & f\left(\widetilde{\mathbf{P}}\right)\prec\widetilde{\mathbf{P}},\ \widetilde{\mathbf{P}}\in\mathcal{R}_{\mathbf{R}^{a}},\label{eq:condition P-tilde}
\end{align}
then\textit{ }$\mathcal{P}_{\mathbf{R}^{a}}\cap\mathcal{R}_{\mathbf{R}^{a}}\neq\textrm{Ø}.$
\end{lemma}

Note that if condition (\ref{eq:condition P-tilde}) is satisfied,
based on the monotonicity of the operator $f$ in Lemma \ref{Lemma: Monotonicity f(.)},
we can continuously apply the operator $f$ to such a $\widetilde{\mathbf{P}}$
and obtain 
\begin{align}
 & \mathbf{Q}\prec\cdots\prec f^{k+1}\left(\widetilde{\mathbf{P}}\right)\prec f^{k}\left(\widetilde{\mathbf{P}}\right)\prec\cdots\prec f\left(\widetilde{\mathbf{P}}\right)\prec\widetilde{\mathbf{P}},\\
 & f^{k}\left(\widetilde{\mathbf{P}}\right)\in\mathcal{R}_{\mathbf{R}^{a}},\forall k\in\mathbb{Z}_{0}^{+}.
\end{align}
As a result, the matrix sequence $\left\{ f^{k}\left(\widetilde{\mathbf{P}}\right)\right\} _{k\in\mathbb{Z}_{0}^{+}}$
is monotonically decreasing and bounded from below. By the monotonic
convergence theorem, we have $f^{k}\left(\widetilde{\mathbf{P}}\right)$
converges as $k\rightarrow\infty$. Therefore, we have 
\begin{align*}
 & \lim_{k\rightarrow\infty}f^{k}\left(\widetilde{\mathbf{P}}\right)=\widetilde{\mathbf{P}}^{*},f\left(\widetilde{\mathbf{P}}^{*}\right)=\widetilde{\mathbf{P}}^{*},\widetilde{\mathbf{P}}^{*}\in\mathcal{R}_{\mathbf{R}^{a}}.
\end{align*}
This means that $\mathcal{P}_{\mathbf{R}^{a}}\cap\mathcal{R}_{\mathbf{R}^{a}}$
is not empty and has en element of $\widetilde{\mathbf{P}}^{*}$,
which in turn guarantees the existence of the game value.

It is worth noting that the existence of the matrix pair $\left(\widetilde{\mathbf{P}},\mathbf{R}^{a}\right)$
in Lemma \ref{Lemma: P-tilde bound} is not always guaranteed. In
fact, we have derived various conditions to guarantee the existence
of such a matrix pair $\left(\widetilde{\mathbf{P}},\mathbf{R}^{a}\right)$
in Lemma \ref{Lemma: P-tilde bound}. Specifically, the next Theorem
\ref{Thm: g(T)<T} shows that if there exists a $\mathbf{T}\in\mathbb{S}_{+}^{S}$
satisfying condition (\ref{eq: g(T)<T}), then the matrix pair $\big(\widetilde{\mathbf{P}},\mathbf{R}^{a}\big)$
in Lemma \ref{Lemma: P-tilde bound} is guaranteed to exist. Moreover,
Lemma \ref{Lemma: Suff-Condition-Existence of T} further provides
a sufficient condition to guarantee the existence of a $\mathbf{T}\in\mathbb{S}_{+}^{S}$
that satisfies (\ref{eq: g(T)<T}), which in turn guarantees the existence
of the matrix pair $\big(\widetilde{\mathbf{P}},\mathbf{R}^{a}\big)$
in Lemma \ref{Lemma: P-tilde bound}. We next discuss how to explicitly
construct such a pair of $\left(\widetilde{\mathbf{P}},\mathbf{R}^{a}\right)$
that satisfies condition (\ref{eq:condition P-tilde}) in the following
theorem.

\begin{thm}
\textit{(Construction of $\left(\widetilde{\mathbf{P}},\mathbf{R}^{a}\right)$)}
\emph{\label{Thm: g(T)<T}}Suppose there exists a $\mathbf{T}\in\mathbb{S}_{+}^{S}$
such that
\begin{align}
 & \mathbb{E}\left[\mathbf{A}^{T}\left(\mathbf{B}^{c}\left(k\right)\left(\mathbf{R}^{c}\right)^{-1}\left(\mathbf{B}^{c}\left(k\right)\right)^{T}+\mathbf{T}^{-1}\right)^{-1}\mathbf{A}\right]+\mathbf{Q}\prec\mathbf{T}.\label{eq: g(T)<T}
\end{align}
Let the matrix pair $\big(\widetilde{\mathbf{P}},\mathbf{R}^{a}\big)$
be constructed based on $\mathbf{T}$ as {\small{}
\begin{align}
 & \widetilde{\mathbf{P}}=\label{eq:P-tilde}\\
 & \left(\mathbb{E}\left[\mathbf{B}^{a}\left(k\right)\right]\left(\mathbf{R}^{a}-\mathrm{cov}\left(\mathbf{B}^{a}\left(k\right)\widetilde{\mathbf{P}}^{\frac{1}{2}}\right)\right)^{-1}\mathbb{E}\left[\left(\mathbf{B}^{a}\left(k\right)\right)^{T}\right]+\mathbf{T}^{-1}\right)^{-1},\nonumber \\
 & \mathbf{R}^{a}\succ\mathbb{E}\left[\left(\mathbf{B}^{a}\left(k\right)\right)^{T}\right]\label{eq:Ra}\\
 & \cdot\left(\vphantom{\left(\mathbb{E}\left[\mathbf{A}^{T}\left(\mathbf{B}^{c}\left(k\right)\left(\mathbf{R}^{c}\right)^{-1}\left(\mathbf{B}^{c}\left(k\right)\right)^{T}+\mathbf{T}^{-1}\right)^{-1}\mathbf{A}\right]+\mathbf{Q}\right)^{-1}}\right.\left(\mathbb{E}\left[\mathbf{A}^{T}\left(\mathbf{B}^{c}\left(k\right)\left(\mathbf{R}^{c}\right)^{-1}\left(\mathbf{B}^{c}\left(k\right)\right)^{T}+\mathbf{T}^{-1}\right)^{-1}\mathbf{A}\right]+\mathbf{Q}\right)^{-1}\nonumber \\
 & -\mathbf{T}^{-1}\left.\vphantom{\vphantom{\left(\mathbb{E}\left[\mathbf{A}^{T}\left(\mathbf{B}^{c}\left(k\right)\left(\mathbf{R}^{c}\right)^{-1}\left(\mathbf{B}^{c}\left(k\right)\right)^{T}+\mathbf{T}^{-1}\right)^{-1}\mathbf{A}\right]+\mathbf{Q}\right)^{-1}}}\right)^{-1}\mathbb{E}\left[\mathbf{B}^{a}\left(k\right)\right]+\sigma_{\mathbf{B}^{a}}^{2}\mathrm{Tr}\left(\mathbf{T}\right)\mathbf{I},\nonumber 
\end{align}
}where $\sigma_{\mathbf{B}^{a}}^{2}=\max\Big\{\mathrm{cov}\left(\left(\mathbf{B}^{a}\left(k\right)\right)_{\left(i,j\right)}\right),1\leq i\leq S,1\leq j\leq\left(N_{t}^{a}\left|\mathcal{N}^{a}\right|\right)$$\Big\},$
then the pair $\left(\widetilde{\mathbf{P}},\mathbf{R}^{a}\right)$
given by (\ref{eq:P-tilde}) and (\ref{eq:Ra}) satisfies condition
(\ref{eq:condition P-tilde}) in Lemma \ref{Lemma: P-tilde bound},
which in turn guarantees the existence of the game value.
\end{thm}

\begin{IEEEproof}
Please see Appendix \ref{subsec:Proof-of-Thm: g(T)<T}.
\end{IEEEproof}

Note that the construction of the matrix pair $\left(\widetilde{\mathbf{P}},\mathbf{R}^{a}\right)$
in Theorem \ref{Thm: g(T)<T} is not based on time $k$ because of
the expectations taken w.r.t. the randomness of $\mathbf{B}^{c}\left(k\right)$
and $\mathbf{B}^{a}\left(k\right).$ As such, the matrix pair $\left(\widetilde{\mathbf{P}},\mathbf{R}^{a}\right)$
in Theorem \ref{Thm: g(T)<T} is not time-varying. Besides, under
the existence of a $\mathbf{T}\in\mathbb{S}_{+}^{S}$ satisfying the
condition in Equation (\ref{eq: g(T)<T}), all the inverse matrices
in Equation (\ref{eq:P-tilde}) and (\ref{eq:Ra}) are well-defined. 

The verification of condition (\ref{eq: g(T)<T}) in Theorem \ref{Thm: g(T)<T}
contains two major parts: the expectation on the L.H.S. of equation
(\ref{eq: g(T)<T}) and the conditions for the existence of a matrix
$\mathbf{T}$ that satisfying the nonlinear matrix inequality in (\ref{eq: g(T)<T}).
We want to highlight that, even if the expected value in the L.H.S.
of (\ref{eq: g(T)<T}) is known, finding the existence conditions
for $\mathbf{T}$ is still very difficult due to the nonlinearities
in (\ref{eq: g(T)<T}). To address this challenge, we will further
show that, under condition (\ref{eq: =00005CrhomaxA_k}) in Lemma
\ref{Lemma: Suff-Condition-Existence of T}, there exists a $\mathbf{T}=\mathbf{T}^{*}$
given by (\ref{eq:T*-Closed-form}) that satisfies equation (\ref{eq: g(T)<T}).
As a result, the verification of condition (\ref{eq: g(T)<T}) in
Theorem \ref{Thm: g(T)<T} is reduced to the verification of condition
(\ref{eq: =00005CrhomaxA_k}) in Lemma \ref{Lemma: Suff-Condition-Existence of T}.
Further note that the random variables $\mathbf{U}_{k}$ and $r_{k}$
within the expectation in (\ref{eq: =00005CrhomaxA_k}) depend only
on the multiplicative randomness $\left\{ \mathbf{H}_{j,i}^{c}\left(k\right)\right\} _{\left(j,i\right)\in\mathcal{N}^{c}\times\mathcal{N}^{c}}$
of the controllers that satisfies Assumption \ref{MIMO-Wireless-Fading},
as shown in equation (\ref{eq: g(T)<T}) and (\ref{eq:agg_CSI_B^c}).
Therefore, we can verify condition (\ref{eq: =00005CrhomaxA_k}) based
on the statistics of $\left\{ \mathbf{H}_{j,i}^{c}\left(k\right)\right\} _{\left(j,i\right)\in\mathcal{N}^{c}\times\mathcal{N}^{c}}$,
which in turn leads to the verification of condition (\ref{eq: g(T)<T}). 

We further propose a novel PSD kernel decomposition technique to derive
a sufficient condition for the existence of such a positive definite
matrix $\mathbf{T}$ that satisfies condition (\ref{eq: g(T)<T}),
which in turn is sufficient for the existence of game value. 

We associate each realization of $\mathbf{B}^{c}\left(k\right),\forall k\in\mathbb{Z}_{0}^{+}$,
with two sets of PSD matrices as follows:
\begin{eqnarray}
\mathcal{T}_{k}^{\mathrm{ker}} & \triangleq & \left\{ \left.\widetilde{\mathbf{T}}\in\hat{\mathbb{S}}_{+}^{S}\right|\mathrm{\text{ker}}\left(\widetilde{\mathbf{T}}\mathbf{B}^{c}\left(k\right)\right)=\mathrm{\text{ker}}\left(\widetilde{\mathbf{T}}\right)\right\} ;\\
\mathcal{T}_{k}^{\mathrm{\mathbf{0}}} & \triangleq & \left\{ \left.\widetilde{\mathbf{T}}\in\hat{\mathbb{S}}_{+}^{S}\right|\widetilde{\mathbf{T}}\mathbf{B}^{c}\left(k\right)=\mathbf{0}\right\} .
\end{eqnarray}

For each realization of $\mathbf{B}^{c}\left(k\right),\forall k\in\mathbb{Z}_{0}^{+},$
we let $\mathrm{rank}\left(\mathbf{B}^{c}\left(k\right)\right)=r_{k}$
and denote the SVD 
\begin{align}
 & \mathbf{B}^{c}\left(k\right)\left(\mathbf{R}^{c}\right)^{-1}\left(\mathbf{B}^{c}\left(k\right)\right)^{T}=\mathbf{U}_{k}^{T}\boldsymbol{\Sigma}_{k}\mathbf{U}_{k}\label{eq:SVD B^c}\\
 & =\mathbf{U}_{k}^{T}\mathrm{diag}\left(\left[\left(\boldsymbol{\Sigma}_{k}\right)_{r_{k}},\boldsymbol{0}_{S-r_{k}}\right]\right)\mathbf{U}_{k},\nonumber 
\end{align}
where $\mathbf{U}_{k}\in\mathbb{R}^{S}$ is an orthogonal matrix and
$\boldsymbol{\Sigma}_{k}$ is a diagonal matrix with the diagonal
elements in descending order. 

We have the following theorem on the closed-form PSD kernel decomposition.

\begin{thm}
\textit{(PSD Kernel Decomposition)} \label{Thm: Ker-Decomposition}For
any given $\mathbf{B}^{c}\left(k\right)\neq\mathbf{0},\forall k\in\mathbb{Z}_{0}^{+},$
every $\mathbf{T}\in\mathbb{S}_{+}^{S}$ can be decomposed into a
sum of two PSD matrices as $\mathbf{T}=\mathbf{T}_{k}^{\mathrm{ker}}+\mathbf{T}_{k}^{\mathrm{\mathbf{0}}},$
where $\mathbf{T}_{k}^{\mathrm{ker}}\in\mathcal{T}_{k}^{\mathrm{ker}}$
and $\mathbf{T}_{k}^{\mathrm{\mathbf{0}}}\in\mathcal{T}_{k}^{\mathrm{\mathbf{0}}}.$
The closed-form expressions of $\mathbf{T}_{k}^{\mathrm{ker}}$ and
$\mathbf{T}_{k}^{\mathrm{\mathbf{0}}}$ are given by $\mathbf{T}_{k}^{\mathrm{ker}}=\left(\widetilde{\mathbf{T}}_{k}^{\mathrm{ker}}\right)^{T}\widetilde{\mathbf{T}}_{k}^{\mathrm{ker}}$
and $\mathbf{T}_{k}^{\mathrm{\mathbf{0}}}=\left(\widetilde{\mathbf{T}}_{k}^{\mathrm{\mathbf{0}}}\right)^{T}\widetilde{\mathbf{T}}_{k}^{\mathrm{\mathbf{0}}},$
respectively, where
\begin{align}
\widetilde{\mathbf{T}}_{k}^{\mathrm{ker}} & =\left[\begin{array}{c}
\left(\mathbf{U}_{k}\mathbf{T}\mathbf{U}_{k}^{T}\right)_{r_{k}}^{\frac{1}{2}}\\
\mathbf{L}_{k}\left(\mathbf{U}_{k}\mathbf{T}\mathbf{U}_{k}^{T}\right)_{r_{k}}^{\frac{1}{2}}
\end{array}\right]^{T}\mathbf{U}_{k},\label{eq:T-ker-closed-form}\\
\widetilde{\mathbf{T}}_{k}^{\mathrm{\mathbf{0}}} & =\left[\begin{array}{c}
\mathbf{0}_{r_{s}}\\
\left[-\mathbf{L}_{k},\mathbf{I}_{S-r_{k}}\right]\left(\mathbf{U}_{k}\mathbf{T}\mathbf{U}_{k}^{T}\right)^{\frac{1}{2}}
\end{array}\right]^{T}\mathbf{U}_{k},\label{eq:T-0-closed-form}
\end{align}
and $\mathbf{L}_{k}=\left(\mathbf{U}_{k}\mathbf{T}\mathbf{U}_{k}^{T}\right)_{\left(1:r_{k};r_{k}+1:S\right)}^{T}\left(\mathbf{U}_{k}\mathbf{T}\mathbf{U}_{k}^{T}\right)_{r_{k}}^{-1}$.
\end{thm}

\begin{IEEEproof}
Since $\mathbf{B}^{c}\left(k\right)\neq\mathbf{0}$, we have $r_{k}>0$
and the inverse $\left(\mathbf{U}_{k}\mathbf{T}\mathbf{U}_{k}^{T}\right)_{r_{k}}^{-1}$
in $\mathbf{L}_{k}$ is well-defined. It is clear that $\mathbf{T}_{k}^{\mathrm{ker}}$
and $\mathbf{T}_{k}^{\mathrm{\mathbf{0}}}$ are PSD matrices due to
the construction of a matrix multiplies its transpose. We can then
directly verify that $\mathbf{T}_{k}^{\mathrm{ker}}+\mathbf{T}_{k}^{\mathrm{\mathbf{0}}}=\mathbf{T}$,
$\mathbf{T}_{k}^{\mathrm{ker}}\in\mathcal{T}_{k}^{\mathrm{ker}}$
and $\mathbf{T}_{k}^{\mathrm{\mathbf{0}}}\in\mathcal{T}_{k}^{\mathrm{\mathbf{0}}}.$
\end{IEEEproof}

Based on the proposed PSD kernel decomposition in Theorem \ref{Thm: Ker-Decomposition},
we can readily obtain the following lemma on sufficient condition
for the existence of a positive definite matrix $\mathbf{T}$ that
satisfies condition (\ref{eq: g(T)<T}).

\begin{lemma}
\textit{(Sufficient Condition for Existence of $\mathbf{T}$)} \emph{\label{Lemma: Suff-Condition-Existence of T}}If
the following condition is satisfied
\begin{align}
\rho\Big( & \mathbb{E}\big[\left(\mathbf{U}_{k}^{T}\mathrm{diag}\left(\left[\boldsymbol{0}_{r_{k}},\mathbf{I}_{S-r_{k}}\right]\right)\mathbf{U}_{k}\mathbf{A}\right)\label{eq: =00005CrhomaxA_k}\\
 & \otimes\big(\mathbf{U}_{k}^{T}\mathrm{diag}\left(\left[\boldsymbol{0}_{r_{k}},\mathbf{I}_{S-r_{k}}\right]\right)\mathbf{U}_{k}\mathbf{A}\big)\big]\Big)<1,\nonumber 
\end{align}
 then $\mathbf{T}=\mathbf{T}^{*}$ satisfies condition (\ref{eq: g(T)<T}),
where $\mathbf{T}^{*}$ is given by
\begin{align}
 & \mathbf{T}^{*}=\mathrm{vec}^{-1}\bigg(\Big(\mathbf{I}-\mathbb{E}\big[\left(\mathbf{U}_{k}^{T}\mathrm{diag}\left(\left[\boldsymbol{0}_{r_{k}},\mathbf{I}_{S-r_{k}}\right]\right)\mathbf{U}_{k}\mathbf{A}\right)\label{eq:T*-Closed-form}\\
 & \otimes\left(\mathbf{U}_{k}^{T}\mathrm{diag}\left(\left[\boldsymbol{0}_{r_{k}},\mathbf{I}_{S-r_{k}}\right]\right)\mathbf{U}_{k}\mathbf{A}\right)\big]\Big)^{-1}\mathrm{vec}\Big(\left\Vert \mathbf{A}\right\Vert ^{2}\nonumber \\
 & \cdot\left\Vert \mathbb{E}\left[\mathds{1}_{\left\{ r_{k}\neq0,\left(\boldsymbol{\Sigma}_{k}\right)_{r_{k}}\succ\xi\mathbf{I}_{r_{k}}\right\} }\mathrm{Tr}\left(\left(\boldsymbol{\Sigma}_{k}\right)_{r_{k}}^{-1}\right)\right]\mathbf{I}+\mathbf{Q}\right\Vert \mathbf{I}\Big)\bigg)\nonumber 
\end{align}
and $\xi$ is some positive constant.
\end{lemma}

\begin{IEEEproof}
Please see Appendix \ref{subsec:Proof-of-Lemma: Suff-Condition-Existence of T}.
\end{IEEEproof}

Based on Theorem \ref{Thm: g(T)<T} and Lemma \ref{Lemma: Suff-Condition-Existence of T},
a verifiable sufficient condition for existence of game value can
be readily obtained in the following theorem.

\begin{thm}
\textit{(Verifiable Sufficient Condition)} \emph{\label{Thm: Verifiable Suff Condition}}The
value of the zero-sum game exists if the following two conditions
are satisfied simultaneously: $\left(i\right)$ Condition (\ref{eq: =00005CrhomaxA_k})
is satisfied, and $\left(ii\right)$ $\mathbf{R}^{a}$ is chosen to
satisfy (\ref{eq:Ra}) with $\mathbf{T}=\mathbf{T}^{*}$ in (\ref{eq:T*-Closed-form}). 
\end{thm}

Theorem \ref{Thm: Verifiable Suff Condition} is a direct result of
Lemma \ref{Lemma: P-tilde bound}, Theorem \ref{Thm: g(T)<T}, and
Lemma \ref{Lemma: Suff-Condition-Existence of T}. Specifically, based
on Lemma \ref{Lemma: Suff-Condition-Existence of T}, we know that
under condition (\ref{eq: =00005CrhomaxA_k}), $\mathbf{T}^{*}$ in
(\ref{eq:T*-Closed-form}) satisfies condition (\ref{eq: g(T)<T}).
Based on Theorem \ref{Thm: g(T)<T}, we can construct a pair of $\left(\widetilde{\mathbf{P}},\mathbf{R}^{a}\left(\mathbf{T}^{*}\right)\right)$
to satisfy condition (\ref{eq:condition P-tilde}) in Lemma \ref{Lemma: P-tilde bound},
which is sufficient for the existence of the game value.

The verification of conditions in Theorem \ref{Thm: Verifiable Suff Condition}
also contains two parts: the verification of condition (\ref{eq: =00005CrhomaxA_k});
and choose an $\mathbf{R}^{a}$ to satisfy inequality (\ref{eq:Ra})
with $\mathbf{T}=\mathbf{T}^{*}$ given by equality (\ref{eq:T*-Closed-form}).
We can verify these two parts in a similar way as the verifications
of Theorem \ref{Thm: g(T)<T}. Additionally, in the following Section
IV, we will use various practical application examples to illustrate
the verifications of the conditions in Theorem \ref{Thm: Verifiable Suff Condition}.

\subsection{Comparison and Differentiation from Existing Literature}

The existing or classical results on infinite time-horizon LQ zero-sum
difference game in {[}10{]}, {[}11{]}, {[}20{]}-{[}26{]} has the following
two key restrictions. Specifically, to ensure the existence of the
game value, the GARE associated with the infinite time-horizon LQ
zero-sum difference game in {[}10{]}, {[}11{]}, {[}20{]}-{[}26{]}
must possess a positive definite solution. To guarantee such solvability
of the GARE, the controllability or stabilizability is an indispensable
prerequisite in {[}10{]}, {[}11{]}, {[}20{]}-{[}26{]}. To be more
specific, if GARE has a positive definite solution, then the closed-loop
system, or equivalently the pair $\left(\mathbf{A},\mathbf{B}^{c}\right)$,
where $\mathbf{B}^{c}$ is the controller input gain matrix, must
be either controllable or stabilizable. As a result, the existing
literature {[}10{]}, {[}11{]}, {[}20{]}-{[}26{]} all restrict the
closed-loop system to be either controllable or stabilizable. Moreover,
the existing results in {[}10{]}, {[}11{]}, {[}20{]}-{[}26{]} also
require the\textbf{ }attacker's weight matrix $\mathbf{R}^{a}$ to
satisfy an implicit concavity constraint of 
\begin{align}
 & \mathbf{R}^{a}\succ\left(\mathbf{B}^{a}\right)^{T}\mathbf{P}^{*}\left(\mathbf{R}^{a}\right)\mathbf{B}^{a},\label{eq:implicit cons gamma-weakness}
\end{align}
where $\mathbf{P}^{*}\left(\mathbf{R}^{a}\right)$ is the minimal
positive definite solution to the GARE. For more details of the implicit
constraint (\ref{eq:implicit cons gamma-weakness}), please refer
to Condition (5.3) in {[}20{]}, Paragraph 1 of Section II-A in {[}21{]},
Remark 1 in {[}22{]}, Condition (8) in {[}23{]}, Assumption 2.1 in
{[}24{]}, Theorem 1 in {[}25{]} and Condition 2) in Section III-B
in {[}26{]}. We want to highlight that both sides of (\ref{eq:implicit cons gamma-weakness})
are functions of the variable $\mathbf{R}^{a}$. However, the $\mathbf{P}^{*}\left(\mathbf{R}^{a}\right)$
in R.H.S. of (\ref{eq:implicit cons gamma-weakness}) is dependent
of $\mathbf{R}^{a}$ but has no closed-form. As such, condition (\ref{eq:implicit cons gamma-weakness})
in existing results cannot provide us with an explicit closed-form
characterization or any practical insight on the feasible domain of
$\mathbf{R}^{a}$, wherein the game value is guaranteed to exist.

Compared to the existing works {[}10{]}, {[}11{]}, {[}20{]}-{[}26{]},
the key advantages of our results are two-fold: (i) our main theorems,
i.e., Theorems \ref{Thm: General case- suff and necessary}- \ref{Thm: Verifiable Suff Condition},
are not subject to the restriction of controllability or stabilizability
of the closed-loop systems and are therefore applicable to a more
extensive range of scenarios, including intermittent controllable
or almost surely uncontrollable closed-loop systems; (ii) our main
theorems offer explicit characterizations of the feasible domain of
the attackers' weight matrices $\mathbf{R}^{a}$. Specifically, our
Theorems \ref{Thm: g(T)<T}, \ref{Thm: Ker-Decomposition}, and \ref{Thm: Verifiable Suff Condition},
together with Lemma \ref{Lemma: Suff-Condition-Existence of T}, provide
various conditions under which the the MGARE is still solvable in
the context of intermittent controllability or almost surely uncontrollability.
In addition, our results also provide an explicit closed-form characterization
of the feasible domain of $\mathbf{R}^{a}$. We want to highlight
that the R.H.S. of (\ref{eq:Ra}) is independent of $\mathbf{R}^{a}$.
Therefore, compared to the existing literature {[}10{]}, {[}11{]},
{[}20{]}-{[}26{]} (i.e., condition (\ref{eq:implicit cons gamma-weakness})),
our work offers an insightful and explicit closed-form picture of
the feasible domain of $\mathbf{R}^{a},$ which guarantees the existence
of the game value. 

\section{Application Examples and Tightness Analysis }

To emphasize the impact of our results in practical scenarios, in
this section we consider several concrete application examples where
the existing literature {[}10{]}, {[}11{]}, {[}20{]}-{[}26{]} fails
to analyze because the closed-loop systems are not controllable or
stabilizable at every time slots. We discuss the tightness of our
proposed verifiable sufficient conditions and provide various practical
application insights.

\begin{example}
\emph{\label{Exmp: E.g.1}}The state transition matrix $\mathbf{A}$
is symmetric. The total numbers of transmit antennas of the controllers
and the attackers satisfy $N_{t}^{c}\left|\mathcal{N}^{c}\right|\geq S$
and $N_{t}^{a}\left|\mathcal{N}^{a}\right|\geq S$, respectively.
The transmissions of the controllers are highly correlated such that
the aggregated CSI $\mathbf{B}^{c}\left(k\right)$ of the controllers
is i.i.d. over the time-slots with $\mathrm{Pr}\left(\mathrm{rank}\left(\mathbf{B}^{c}\left(k\right)\right)=S\right)=\delta,\ \mathrm{Pr}\left(\mathbf{B}^{c}\left(k\right)=\mathbf{0}\right)=1-\delta$.
Besides, the minimum non-zero singular value of $\mathbf{B}^{c}\left(k\right)$
is bounded away from 0 w. p. 1.
\end{example}

Applying Theorem \ref{Thm: Verifiable Suff Condition} to Example
\ref{Exmp: E.g.1}, we obtain the following verifiable sufficient
condition.

\begin{cor}
\textit{(Verifiable Sufficient Condition for Example \ref{Exmp: E.g.1})}\emph{\label{Coro: Example 1 Suff}}
The game value of Example \ref{Exmp: E.g.1} exists if $\delta$ and
$\mathbf{R}^{a}$ satisfy the following conditions 
\begin{align}
 & \delta>1-\rho^{-2}\left(\mathbf{A}\right),\label{eq: delta condition-suff-example-3}\\
 & \mathbf{R}^{a}\succ\gamma_{1}\left(1-\left(1-\delta\right)\rho^{2}\left(\mathbf{A}\right)\right)^{-2}\mathbf{I},\label{eq:Ra condition-suff-example-3}
\end{align}
 where $\gamma_{1}$ is some positive constant scalar.
\end{cor}

\begin{IEEEproof}
Please see Appendix \ref{subsec:Proof-of-Coro: Example 1 Suff}.
\end{IEEEproof}

The following corollary characterizes necessary conditions for existence
of game value of Example \ref{Exmp: E.g.1}.

\begin{cor}
\textit{(Verifiable Necessary Condition for Example \ref{Exmp: E.g.1})}\emph{
\label{Coro: Example 1 Necessary}}Consider the same system configurations
as in Example \ref{Exmp: E.g.1}. We further assume that $\mathbf{A}$
is invertible and $\mathrm{rank}\left(\mathbb{E}\left[\mathbf{B}^{a}\left(k\right)\right]\right)=S.$
In this case, if the game value exists, then $\delta$ must satisfy
condition (\ref{eq: delta condition-suff-example-3}), and $\mathbf{R}^{a}$
must satisfy
\begin{align}
 & \left\Vert \mathbf{R}^{a}\right\Vert >\gamma_{4}\left(1-\left(1-\delta\right)\rho^{2}\left(\mathbf{A}\right)\right)^{-2}.\label{eq:Ra L2-norm-Necessary Condition}
\end{align}
where $\gamma_{4}$ is some positive constant scalar.
\end{cor}

\begin{IEEEproof}
Please see Appendix \ref{subsec:Proof-of-Coro: Example 1 Necessary}.
\end{IEEEproof}

\begin{remrk}
\textsl{(}\textit{Tightness of}\textsl{ $\delta$)} \emph{\label{Remark: Tightness}}Under
the system configurations in Example \ref{Exmp: E.g.1}, condition
(\ref{eq: =00005CrhomaxA_k}) in the proposed verifiable sufficient
condition is reduced to $\delta>1-\rho^{-2}\left(\mathbf{A}\right),$
which is both necessary and sufficient. This also delivers a physical
insight that for a more unstable plant, i.e., a larger $\rho\left(\mathbf{A}\right)$,
more communication resource is needed to guarantee the existence of
game value, i.e., the controllers should have larger transmission
probabilities to induce a larger $\delta$.
\end{remrk}

\begin{remrk}
\textsl{(}\textit{Tightness of}\textsl{ $\mathbf{R}^{a}$)} \emph{\label{Remark: Tightness-1}}Under
the system configurations in Example \ref{Exmp: E.g.1}, requirement
on $\mathbf{R}^{a}\left(\mathbf{T}^{*}\right)$ in the proposed verifiable
sufficient condition is reduced to condition (\ref{eq:Ra condition-suff-example-3}),
which implies that $\left\Vert \mathbf{R}^{a}\right\Vert $ \\=$\Omega\left(\left(1-\left(1-\delta\right)\rho^{2}\left(\mathbf{A}\right)\right)^{-2}\right)$.
Corollary \ref{Coro: Example 1 Necessary} demonstrates that, with
the additional requirements of invertible $\mathbf{A}$ and $\mathrm{rank}\left(\mathbb{E}\left[\mathbf{B}^{a}\left(k\right)\right]\right)=S$,
$\left\Vert \mathbf{R}^{a}\right\Vert =\Omega\left(\left(1-\left(1-\delta\right)\rho^{2}\left(\mathbf{A}\right)\right)^{-2}\right)$
is also necessary for the existence of the game value. Therefore,
the verifiable sufficient condition in Theorem \ref{Thm: Verifiable Suff Condition}
provides an $\mathbf{R}^{a}$ whose spectral norm is order-wise optimal.
\end{remrk}

We next consider a practical application example, where there is no
signal interference among the controllers.

\begin{example}
\emph{(}\textit{No Wireless Interference}\textit{\emph{ }}\emph{Among
Controllers)\label{Exmp: E.g.2}} The number of actuators, controllers
and attackers are the same, i.e., $\left|\mathcal{N}^{c}\right|=\left|\mathcal{N}^{a}\right|=M$.
The number of transmit antenna of each controller is equal to the
number of receive antenna of each actuator, i.e., $N_{t}^{c}=N_{r}=N$.
The plant state dimension is $S=M\cdot N$. Let $\delta_{j}^{c}\left(k\right)\in\{0,1\}$
and $\delta_{l}^{a}\left(k\right)\in\{0,1\}$ be the i.i.d. Bernoulli
random variables that represent whether the $j$-th controller and
the $l$-th attacker are active to transmit control action $\mathbf{u}_{j}^{c}\left(k\right)$
and malicious attack injection $\mathbf{u}_{l}^{a}\left(k\right)$
or not, respectively, with $\mathbb{E}\left[\delta_{j}^{c}\left(k\right)\right]=\delta_{j}^{c}$
and $\mathbb{E}\left[\delta_{l}^{a}\left(k\right)\right]=\delta_{l}^{a}$.
The control input gain matrix of the $i$-th actuator is $\mathbf{B}_{i}=\left[\begin{array}{ccc}
\mathbf{0}_{\left(i-1\right)N\times S} & \widetilde{\mathbf{B}}_{i}^{T} & \mathbf{0}_{\left(S-iN\right)\times S}\end{array}\right]^{T},\forall1\leq i\leq M$, where $\widetilde{\mathbf{B}}_{i}\in\mathbb{R}^{N\times N}$ is
full rank. There is no interference communication link between the
actuators and controllers, i.e., $\mathbf{H}_{j,i}^{c}\left(k\right)=\mathbf{0}_{N},\forall i\neq j,\forall1\leq i,j\leq M,\forall k\in\mathbb{Z}^{+}.$
Besides, for the direct communication link, the MIMO fading channel
matrix $\mathbf{H}_{i,i}^{c}\left(k\right)$ is full rank w. p. 1
,$\forall1\leq i\leq M,\forall k\in\mathbb{Z}^{+}.$ 
\end{example}

Note that in Example \ref{Exmp: E.g.2}, the equivalent control action
input gain matrix $\mathbf{B}^{c}\left(k\right)$ is given by
\begin{align}
 & \mathbf{B}^{c}\left(k\right)=\mathrm{diag}\left(\left[\delta_{i}^{c}\left(k\right)\widetilde{\mathbf{B}}_{i}\mathbf{H}_{i,i}^{c}\left(k\right)\right]_{1\leq i\leq M}\right)\in\mathbb{R}^{MN\times MN}.
\end{align}
The pair $\left(\mathbf{A},\mathbf{B}^{c}\left(k\right)\right)$ in
Example \ref{Exmp: E.g.2} is not controllable or stabilizable at
every time slot. This is because the $i$-th controller's random access
variable $\delta_{i}^{c}\left(k\right),\forall1\leq i\leq M,$ is
a Bernoulli random variable. Under the specific realization that $\delta_{i}^{c}\left(k\right)=0,\forall1\leq i\leq M,$
$\mathbf{B}^{c}\left(k\right)$ will become $\mathbf{B}^{c}\left(k\right)=\mathbf{0}$
and the pair $\left(\mathbf{A},\mathbf{0}\right)$ is clearly uncontrollable.
As a result, the existing infinite time-horizon LQ zero-sum difference
game results in {[}10{]}, {[}11{]}, {[}20{]}-{[}26{]}, which are strictly
built on controllability or stabilizability, are not applicable to
the above Example \ref{Exmp: E.g.2}, where $\left(\mathbf{A},\mathbf{B}^{c}\left(k\right)\right)$
can be uncontrollable for certain time slots. In other words, the
practical scenario considered in Example \ref{Exmp: E.g.2} falls
outside the application scope of the existing literature {[}10{]},
{[}11{]}, {[}20{]}-{[}26{]}, all of which fail to provide insight
or characterization of the existence and achievability of the game
value. Applying Theorem \ref{Thm: Verifiable Suff Condition} to Example
\ref{Exmp: E.g.2}, we arrive at the following sufficient condition
for existence of game value.

\begin{cor}
\textit{(Verifiable Sufficient Condition for Example }\textit{\emph{\ref{Exmp: E.g.2}}}\emph{)\label{Coro: Example 2 Suff}}
The game value of Example \ref{Exmp: E.g.2} exists if the following
two conditions are satisfied simultaneously:
\end{cor}

\begin{itemize}
\item $\rho\left(\mathbf{A}_{1}\right)<1$, where 
\begin{align*}
\mathbf{A}_{1}= & \left[\begin{array}{ccc}
\left(1-\delta_{1}^{c}\right)^{\frac{1}{2}}\mathbf{I}_{N}\\
 & \ddots\\
 &  & \left(1-\delta_{M}^{c}\right)^{\frac{1}{2}}\mathbf{I}_{N}
\end{array}\right]\mathbf{A}.
\end{align*}
\item $\mathbf{R}^{a}$ is chosen to satisfy (\ref{eq:Ra}) with $\mathbf{T}^{*}$
given by
\begin{align}
\mathbf{\mathbf{T}}^{*}= & \mathrm{vec}^{-1}\bigg(\left(\mathbf{I}-\mathbf{A}_{1}\otimes\mathbf{A}_{1}\right)^{-1}\\
 & \mathrm{vec}\left(\left\Vert \mathbf{A}\right\Vert ^{2}\mathbb{E}\left[\mathds{1}_{\left\{ r_{k}\neq0\right\} }\mathrm{Tr}\left(\left(\boldsymbol{\Sigma}_{k}\right)_{r_{k}}^{-1}\right)\right]\mathbf{I}+\mathbf{Q}\right)\bigg).\nonumber 
\end{align}
\end{itemize}

Condition $\rho\left(\mathbf{A}_{1}\right)<1$ in Corollary \ref{Coro: Example 2 Suff}
requires every controller has a relatively large transmission probability.
In the special case that $\mathbf{A}=\mathrm{diag}\left(\mathbf{A}_{1},\cdots,\mathbf{A}_{M}\right)$
with $\mathbf{A}_{i}\in\mathbb{R}^{N\times N},\forall1\leq i\leq M$,
each controller should have a sufficiently large transmission probability
to address the instability of its corresponding subsystem, i.e., $\delta_{i}^{c}$
should be larger than $\big(1-\rho^{-2}\left(\mathbf{A}\right)\big)$
for all $1\leq i\leq M$.

We finally consider an application example, where controllers and
actuators are equipped with an increased number of antennas.
\begin{example}
\emph{(}\textit{Equipment of Increased Number of Antennas}\emph{)\label{Exmp: E.g.3}}
Let $\left|\mathcal{N}^{c}\right|=\left|\mathcal{N}^{a}\right|=M$.
$N_{t}^{c}$ and $N_{r}$ are both increased to satisfy $N_{t}^{c}>M\cdot N_{r}\geq S$.
Let $\delta_{j}^{c}\left(k\right)\in\{0,1\}$ and $\delta_{l}^{a}\left(k\right)\in\{0,1\}$
be the i.i.d. Bernoulli random variables that represent whether the
$j$-th controller and the $l$-th attacker are active to transmit
control action $\mathbf{u}_{j}^{c}\left(k\right)$ and malicious attack
injection $\mathbf{u}_{l}^{a}\left(k\right)$ or not, respectively,
with $\mathbb{E}\left[\delta_{j}^{c}\left(k\right)\right]=\delta_{j}^{c}$
and $\mathbb{E}\left[\delta_{l}^{a}\left(k\right)\right]=\delta_{l}^{a}$.
The concatenated row block actuator input gain matrix $\mathbf{B}=\left[\mathbf{B}_{1},\mathbf{B}_{2},\cdots,\mathbf{B}_{M}\right]$
is full row rank. The concatenated column block MIMO channel matrices
of $i$-th controller $\mathbf{H}_{i}^{c}\left(k\right)=\left[\mathbf{H}_{i,1}^{c}\left(k\right);\mathbf{H}_{i,2}^{c}\left(k\right)\cdots;\mathbf{H}_{i,M}^{c}\left(k\right)\right]$
is full row rank w.p. 1, $\forall1\leq i\leq M,\forall k\in\mathbb{Z}^{+}.$ 
\end{example}

In Example \ref{Exmp: E.g.3}, the equivalent control action input
gain matrix $\mathbf{B}^{c}\left(k\right)$ is given by
\begin{align}
 & \mathbf{B}^{c}\left(k\right)=\left[\delta_{i}^{c}\left(k\right)\mathbf{B}\mathbf{H}_{i}^{c}\left(k\right)\right]_{1\leq i\leq M}\in\mathbb{R}^{S\times M\cdot N_{t}^{c}}.
\end{align}
The pair $\left(\mathbf{A},\mathbf{B}^{c}\left(k\right)\right)$ in
Example \ref{Exmp: E.g.3} can also be uncontrollable at certain time
slots under the possible realization of $\delta_{i}^{c}\left(k\right)=0,\forall1\leq i\leq M.$
Consequently, the existing results in {[}10{]}, {[}11{]}, {[}20{]}-{[}26{]},
which are based on controllability or stabilizability, also fail to
be applicable to the above Example \ref{Exmp: E.g.3}. In contrast,
we can directly apply our Theorems 2, 3, and 4, along with Lemma 6,
to Example \ref{Exmp: E.g.3} to analyze the existence of the game
value and obtain the following Corollary \ref{Coro: Example 3 Suff}.

\begin{cor}
\textit{(Verifiable Sufficient Condition for Example}\textit{\emph{
\ref{Exmp: E.g.3}}}\emph{)\label{Coro: Example 3 Suff}} The game
value of Example \ref{Exmp: E.g.3} exists if the following two conditions
are satisfied simultaneously:
\end{cor}

\begin{itemize}
\item $\prod_{i=1}^{M}\left(1-\delta_{i}^{c}\right)\rho^{2}\left(\mathbf{A}\right)<1.$ 
\item $\mathbf{R}^{a}$ is chosen to satisfy (\ref{eq:Ra}) with $\mathbf{T}^{*}$
given by
\begin{align}
\mathbf{\mathbf{T}}^{*}= & \mathrm{vec}^{-1}\Big(\Big(\mathbf{I}-\prod_{i=1}^{M}\left(1-\delta_{i}^{c}\right)\mathbf{A}\otimes\mathbf{A}\Big)^{-1}\\
 & \cdot\mathrm{vec}\left(\left\Vert \mathbf{A}\right\Vert ^{2}\mathbb{E}\left[\mathds{1}_{\left\{ r_{k}\neq0\right\} }\mathrm{Tr}\left(\left(\boldsymbol{\Sigma}_{k}\right)_{r_{k}}^{-1}\right)\right]\mathbf{I}+\mathbf{Q}\right)\Big).\nonumber 
\end{align}
\end{itemize}

In contrast to Example \ref{Exmp: E.g.2}, condition $\prod_{i=1}^{M}\left(1-\delta_{i}^{c}\right)\rho^{2}\left(\mathbf{A}\right)<1$
in Corollary \ref{Coro: Example 3 Suff} does not require all the
controllers to have a large transmission probability. It only requires
that at least one of the controllers has a large transmission probability
$\delta_{i}^{c}$. As a result, a larger number of antennas at the
controllers and actuators is more favorable for the existence of the
game value.

\section{Achievable Schemes}

The the standard LQ Nash equilibrium computations in {[}10{]}, {[}11{]},
{[}20{]} are based on the stabilizability of the closed-loop system
and are restricted to the class of strictly feedback stabilizing policies,
which is not applicable to our work. In this section, we discuss the
control and attack policies that are able to achieve the value of
the infinite-time horizon zero-sum game in a general class of policies
that may be possibly unstable, subject to intermittent controllability
or almost sure uncontrollability. 

We associate every $T_{0}\in\mathbb{Z}^{+}$ with two scalar functions:
\begin{align}
\alpha\left(T_{0}\right)= & \mathrm{max}\left\{ \left\Vert \mathbf{x}\left(0\right)\right\Vert _{2}^{2},\left\Vert \mathbf{\mathbf{W}}+\sum_{i\in\mathcal{N}^{c}}\mathbf{B}_{i}^{T}\mathbf{V}\mathbf{B}_{i}\right\Vert ^{2}\right\} \label{eq:alpha(T_0)}\\
 & \cdot\mathrm{Tr}\left(\sum_{i=0}^{T_{0}-1}\left(\mathbb{E}\left[\mathbf{\Xi}\left(k\right)\otimes\mathbf{\Xi}\left(k\right)\right]\right)^{i}\right),\nonumber \\
\beta\left(T_{0}\right)= & \mathrm{min}\Big\{\beta\in\mathbb{Z}^{+},\beta\geq T_{0}:\label{eq:beta(T_0)}\\
 & \left\Vert \mathbf{P}^{*}-f^{\beta-T_{0}}\left(\mathbf{Q}\right)\right\Vert \alpha\left(T_{0}\right)<1\Big\},\nonumber 
\end{align}
where $\mathbf{\Xi}\left(k\right)=\left(\mathbf{I}-\mathbf{B}\left(k\right)\Phi^{-1}\left(\mathbf{P}^{*};k\right)\mathbf{B}^{T}\left(k\right)\mathbf{P}^{*}\right)\mathbf{A}.$

The physical meanings of the two scalar functions $\alpha\left(T_{0}\right)$
and $\beta\left(T_{0}\right)$ can be interpreted as follows. Under
the control and attack inputs $\left\{ \left(\mathbf{u}^{c}\right)^{*}\left(\mathbf{\mathbf{P}^{*}};k\right),\left(\mathbf{u}^{a}\right)^{*}\left(\mathbf{\mathbf{P}^{*}};k\right)\right\} _{k\in\left[0,T_{0}-1\right]}$,
the closed-loop dynamics of the system state becomes, $\forall k\in\left[0,T_{0}-1\right],$
\begin{align}
 & \mathbf{x}\left(k+1\right)=\mathbf{\Xi}\left(k\right)\mathbf{x}\left(k\right)+\sum_{i\in\mathcal{N}^{c}}\mathbf{B}_{i}\mathbf{v}_{i}\left(k\right)+\mathbf{w}\left(k\right).
\end{align}
The Euclidean norm of the terminal state $\mathbf{x}\left(T_{0}\right)$
satisfies $\mathbb{E}\left[\left\Vert \mathbf{x}\left(T_{0}\right)\right\Vert _{2}^{2}\right]\leq\alpha\left(T_{0}\right)$.
As a result, the scalar $\alpha\left(T_{0}\right)$ in (\ref{eq:alpha(T_0)})
quantifies the maximum expected growth of the plant state's Euclidean
norm in the first $T_{0}$ time-slots. The scalar $\beta\left(T_{0}\right)$
in (\ref{eq:beta(T_0)}) characterizes the minimum positive integer
$\beta$ such that the weighted Euclidean norm of $\mathbf{x}\left(T_{0}\right)$
satisfies
\begin{align}
 & \mathbb{E}\left[\left\Vert \mathbf{x}\left(T_{0}\right)\right\Vert _{\left(\mathbf{\mathbf{P}^{*}}-f^{\beta-T_{0}}\left(\mathbf{Q}\right)\right)}^{2}\right]\leq\left\Vert \mathbf{P}^{*}-f^{\beta-T_{0}}\left(\mathbf{Q}\right)\right\Vert \alpha\left(T_{0}\right)<1.\label{eq: x-T_0 norm <1}
\end{align}

We next define a class of $\beta\left(T_{0}\right)$-length, $\forall T_{0}\in\mathbb{Z}^{+},$
admissible policy pairs
\begin{align}
\left(\mathcal{V}^{c}\times\mathcal{V}^{a}\right)_{\beta\left(T_{0}\right)}\triangleq & \Big\{\left\{ \mathbf{u}^{c}\left(k\right),\mathbf{u}^{a}\left(k\right)\right\} _{k\in\left[0,\beta\left(T_{0}\right)\right)}:\\
 & \mathbb{E}\left[\left\Vert \mathbf{x}\left(T_{0}\right)\right\Vert _{2}^{2}\right]\leq\alpha\left(T_{0}\right)\Big\}.\nonumber 
\end{align}

We now construct a control policy $\overline{\boldsymbol{\mu}}_{\beta\left(T_{0}\right)}^{c}=\left\{ \overline{\mathbf{u}}^{c}\left(k\right)\right\} _{k\in\left[0,\beta\left(T_{0}\right)\right)}$
and an attack policy $\overline{\boldsymbol{\mu}}_{\beta\left(T_{0}\right)}^{a}=\left\{ \overline{\mathbf{u}}^{a}\left(k\right)\right\} _{k\in\left[0,\beta\left(T_{0}\right)\right)}$
as follows 
\begin{align*}
 & \overline{\mathbf{u}}^{c}\left(k\right)=\begin{cases}
\left(\mathbf{u}^{c}\right)^{*}\left(\mathbf{P}^{*};k\right), & 0\leq k<T_{0};\\
\left(\mathbf{u}^{c}\right)^{*}\left(f^{\beta\left(T_{0}\right)-k-1}\left(\mathbf{Q}\right);k\right), & T_{0}\leq k<\beta\left(T_{0}\right).
\end{cases}\\
 & \overline{\mathbf{u}}^{a}\left(k\right)=\begin{cases}
\left(\mathbf{u}^{a}\right)^{*}\left(\mathbf{P}^{*};k\right), & 0\leq k<T_{0};\\
\left(\mathbf{u}^{a}\right)^{*}\left(f^{\beta\left(T_{0}\right)-k-1}\left(\mathbf{Q}\right);k\right), & T_{0}\leq k<\beta\left(T_{0}\right).
\end{cases}
\end{align*}

The control and attack policy pair $\left\{ \overline{\boldsymbol{\mu}}_{\beta\left(T_{0}\right)}^{c},\overline{\boldsymbol{\mu}}_{\beta\left(T_{0}\right)}^{a}\right\} $
adopts the steady state policy associated with the minimal solution
$\mathbf{P}^{*}$ to the MGARE (\ref{eq: general-GARE}) for time
slots that are equal or less than $T_{0}$, and another closed-form
time-varying policy for all time slots that are larger than $T_{0}$.
Note that $\left(\mathbf{u}^{c}\right)^{*}$ and $\left(\mathbf{u}^{a}\right)^{*}$
have the same form as the steady state controls and attacks defined
in {[}10{]}, {[}11{]}, {[}20{]}, and can be computed according to
equations (\ref{eq:def-u_c*}) and (\ref{eq:def-u_a*}) in Section
III-B, respectively. $\mathbf{P}^{*}$ is the minimal solution to
the MGARE and can be computed according to (\ref{eq: general-GARE})
in Section III-A. 

We have the following theorem on the achievability of the value of
the zero-sum game. 

\begin{thm}
\textit{(General Achievable Scheme)}\emph{ \label{Thm: achievability of NE}}Let
the class of admissible policy pairs be $\left(\mathcal{U}^{c}\right)^{\left|\mathcal{N}^{c}\right|}\times\left(\mathcal{U}^{a}\right)^{\left|\mathcal{N}^{a}\right|}=\left(\mathcal{V}^{c}\times\mathcal{V}^{a}\right)_{\beta\left(\infty\right)}$.
Assume $\mathcal{P}_{\mathbf{R}^{a}}\cap\mathcal{R}_{\mathbf{R}^{a}}\neq\textrm{Ø}$
such that the value of the infinite-time horizon zero-sum game exists.
Let $T_{0}\rightarrow\infty$, the policy pair $\left\{ \overline{\boldsymbol{\mu}}_{\beta\left(\infty\right)}^{c},\overline{\boldsymbol{\mu}}_{\beta\left(\infty\right)}^{a}\right\} $
constitutes a Nash equilibrium solution and achieves the value of
the zero-sum game, i.e., for all $\left\{ \boldsymbol{\mu}_{\beta\left(\infty\right)}^{c},\boldsymbol{\mu}_{\beta\left(\infty\right)}^{a}\right\} \in\left(\mathcal{V}^{c}\times\mathcal{V}^{a}\right)_{\beta\left(\infty\right)},$
we have $J\left(\overline{\boldsymbol{\mu}}_{\beta\left(\infty\right)}^{c},\boldsymbol{\mu}_{\beta\left(\infty\right)}^{a}\right)\leq J^{*}\leq J\left(\boldsymbol{\mu}_{\beta\left(\infty\right)}^{c},\overline{\boldsymbol{\mu}}_{\beta\left(\infty\right)}^{a}\right),$
where
\begin{align}
 & J^{*}=J\left(\overline{\boldsymbol{\mu}}_{\beta\left(\infty\right)}^{c},\overline{\boldsymbol{\mu}}_{\beta\left(\infty\right)}^{a}\right)=\mathrm{Tr}\left(\mathbf{P}^{*}\left(\mathbf{\mathbf{W}}+\sum_{i\in\mathcal{N}^{c}}\mathbf{B}_{i}^{T}\mathbf{V}\mathbf{B}_{i}\right)\right).
\end{align}
\end{thm}

\begin{IEEEproof}
Please see Appendix \ref{subsec:Proof-of-Thm: achievability of NE}.
\end{IEEEproof}

Note that for unstable policies in $\left(\mathcal{V}^{c}\times\mathcal{V}^{a}\right)_{\beta\left(T_{0}\right)}$,
both $\alpha\left(T_{0}\right)$ and $\beta\left(T_{0}\right)-T_{0}$
are infinite as $T_{0}\rightarrow\infty$. In this case, each of the
proposed policies will consist of two policy sequences, both of infinite
length, i.e., $\overline{\boldsymbol{\mu}}_{\beta\left(\infty\right)}^{c}=\Big\{\left(\mathbf{u}^{c}\right)^{*}\left(\mathbf{P}^{*}\right),\cdots$
$,\left(\mathbf{u}^{c}\right)^{*}\left(f^{\infty}\left(\mathbf{Q}\right)\right),\cdots,\left(\mathbf{u}^{c}\right)^{*}\left(f^{0}\left(\mathbf{Q}\right)\right)\Big\}$
and $\overline{\boldsymbol{\mu}}_{\beta\left(\infty\right)}^{a}=\Big\{\left(\mathbf{u}^{a}\right)^{*}\left(\mathbf{P}^{*}\right),\cdots$
$,\left(\mathbf{u}^{a}\right)^{*}\left(f^{\infty}\left(\mathbf{Q}\right)\right),\cdots,\left(\mathbf{u}^{a}\right)^{*}\left(f^{0}\left(\mathbf{Q}\right)\right)\Big\}$.
We cannot neglect the second infinite length sequence and just pick
the first infinite length sequence, otherwise it cannot achieve the
Nash equilibrium in $\left(\mathcal{V}^{c}\times\mathcal{V}^{a}\right)_{\beta\left(\infty\right)}$.
The same situation also occurs in the standard settings in {[}10{]},
{[}11{]}, {[}20{]}. By replacing the MGARE iteration $f\left(\cdot\right)$
to the standard GARE iteration $g\left(\cdot\right)$, the policies
pair $\left(\overline{\boldsymbol{\mu}}_{\beta\left(\infty\right)}^{c},\text{\ensuremath{\overline{\boldsymbol{\mu}}_{\beta\left(\infty\right)}^{a}}}\right)$
becomes $\overline{\boldsymbol{\mu}}_{\beta\left(\infty\right)}^{c}=\Big\{\left(\mathbf{u}^{c}\right)^{*}\left(g^{\infty}\left(\mathbf{Q}\right)\right),\cdots$
$,\left(\mathbf{u}^{c}\right)^{*}\left(g^{\infty}\left(\mathbf{Q}\right)\right),\cdots,\left(\mathbf{u}^{c}\right)^{*}\left(g^{0}\left(\mathbf{Q}\right)\right)\Big\}$
and $\overline{\boldsymbol{\mu}}_{\beta\left(\infty\right)}^{a}=\Big\{\left(\mathbf{u}^{a}\right)^{*}\left(g^{\infty}\left(\mathbf{Q}\right)\right),\cdots$
\\$,\left(\mathbf{u}^{a}\right)^{*}\left(g^{\infty}\left(\mathbf{Q}\right)\right),\cdots,\left(\mathbf{u}^{a}\right)^{*}\left(f^{0}\left(\mathbf{Q}\right)\right)\Big\}$,
which also achieves the Nash equilibrium in the general class of admissible
policies. However, if we only take the first infinite length sequence
to form a steady-state policies pair $\left(\overline{\boldsymbol{\mu}}_{s}^{c},\text{\ensuremath{\overline{\boldsymbol{\mu}}_{s}^{a}}}\right)$
as $\overline{\boldsymbol{\mu}}_{s}^{c}=\Big\{\left(\mathbf{u}^{c}\right)^{*}\left(g^{\infty}\left(\mathbf{Q}\right)\right),\cdots\Big\}$
and $\overline{\boldsymbol{\mu}}_{s}^{a}=\Big\{\left(\mathbf{u}^{a}\right)^{*}\left(g^{\infty}\left(\mathbf{Q}\right)\right),\cdots\Big\}$,
such a policies pair $\left(\overline{\boldsymbol{\mu}}_{s}^{c},\text{\ensuremath{\overline{\boldsymbol{\mu}}_{s}^{a}}}\right)$
is not a Nash equilibrium policies pair and can only achieve the upper
value of the game in the general class of admissible policies that
may not be stabilizing. 

A direct result of Theorem \ref{Thm: achievability of NE} is the
achievability of the Nash equilibrium in the class of strictly feedback
mean-square stabilizing policies.
\begin{cor}
\textit{(Mean-square Stabilizing Policies)}\emph{ \label{Coro: Strictly Stabilizing Class}}If
$\mathrm{\rho}\left(\mathbb{E}\left[\mathbf{\Xi}\left(k\right)\otimes\mathbf{\Xi}\left(k\right)\right]\right)<1$,
then the control and attack actions $\left\{ \left(\mathbf{u}^{c}\right)^{*}\left(\mathbf{P}^{*};k\right),\left(\mathbf{u}^{a}\right)^{*}\left(\mathbf{P}^{*};k\right)\right\} _{k\in\mathbb{Z}_{0}^{+}}$
constitutes a feedback Nash equilibrium in the class of strictly feedback
mean-square stabilizing policies $\mathcal{V}_{s}^{c}\times\mathcal{V}_{s}^{a}$,
where 
\begin{align}
\mathcal{V}_{s}^{c}\times\mathcal{V}_{s}^{a}\triangleq & \Big\{\left(\left\{ \boldsymbol{\mu}_{j}^{c}\right\} _{j\in\mathcal{N}^{c}},\text{\ensuremath{\left\{  \boldsymbol{\mu}_{l}^{a}\right\}  _{l\in\mathcal{N}^{a}}}}\right):\\
 & \limsup_{k\rightarrow\infty}\mathbb{E}\Big[\left\Vert \mathbf{x}\left(k\right)\right\Vert _{2}^{2}\Big]<\infty\Big\}.\nonumber 
\end{align}
\end{cor}

Note that if $\mathrm{\rho}\left(\mathbb{E}\left[\mathbf{\Xi}\left(k\right)\otimes\mathbf{\Xi}\left(k\right)\right]\right)<1$,
then $\mathrm{Tr}\left(\sum_{i=0}^{\infty}\left(\mathbb{E}\left[\mathbf{\Xi}\left(k\right)\otimes\mathbf{\Xi}\left(k\right)\right]\right)^{i}\right)$
is finite, which implies $\alpha\left(T_{0}\right)$ in (\ref{eq:alpha(T_0)})
is finite for all $T_{0}\in\mathbb{Z}^{+}$. This means that $\mathbb{E}\Big[\left\Vert \mathbf{x}\left(\infty\right)\right\Vert _{2}^{2}\Big]$
is finite, and hence $\left\{ \left(\mathbf{u}^{c}\right)^{*}\left(\mathbf{P}^{*};k\right),\left(\mathbf{u}^{a}\right)^{*}\left(\mathbf{P}^{*};k\right)\right\} _{k\in\mathbb{Z}_{0}^{+}}$
is strictly feedback mean-square stabilizing. Based on the definition
of $\beta\left(T_{0}\right)$ in (\ref{eq:beta(T_0)}), for the class
of strictly feedback mean-square stabilizing policies, $\beta\left(T_{0}\right)-T_{0}$
is finite. This is because if the constant matrix $\mathbf{Q}$ is
sufficiently close to $\mathbf{P}^{*}$ such that $\left\Vert \mathbf{P}^{*}-\mathbf{Q}\right\Vert <\frac{1}{\alpha\left(T_{0}\right)},$
then $\beta\left(T_{0}\right)-T_{0}=0$. Otherwise, if $\left\Vert \mathbf{P}^{*}-\mathbf{Q}\right\Vert \nless\frac{1}{\alpha\left(T_{0}\right)}$,
$\beta\left(T_{0}\right)-T_{0}=k_{0}$, where $k_{0}$ is the smallest
finite integer such $\left\Vert \mathbf{P}^{*}-f^{k_{0}}\left(\mathbf{Q}\right)\right\Vert <\frac{1}{\alpha\left(T_{0}\right)}.$
As a result, when we let $T_{0}\rightarrow\infty$, there are only
a finite number of terms $\left(\mathbf{u}^{c}\right)^{*}\left(f^{\beta\left(T_{0}\right)-k-1}\left(\mathbf{Q}\right);k\right)$
and $\left(\mathbf{u}^{a}\right)^{*}\left(f^{\beta\left(T_{0}\right)-k-1}\left(\mathbf{Q}\right);k\right)$,
$\forall T_{0}<k<\beta\left(T_{0}\right)$, in $\overline{\mathbf{u}}^{c}\left(k\right)$
and $\overline{\mathbf{u}}^{a}\left(k\right)$, respectively, which
have negligible impact on the Nash equilibrium of the infinite horizon
problem. It follows that $\left\{ \left(\mathbf{u}^{c}\right)^{*}\left(\mathbf{P}^{*};k\right),\left(\mathbf{u}^{a}\right)^{*}\left(\mathbf{P}^{*};k\right)\right\} _{k\in\mathbb{Z}_{0}^{+}}$
achieves the Nash equilibrium.

\begin{remrk}
\textit{(Backward Compatibility Discussion)} \emph{\label{Remark: Compatibility}}For
linear time-invariant systems without randomness, existing literature
shows that if the steady state policy associated with the minimal
solution to the GARE $\left\{ \left(\mathbf{u}^{c}\right)^{*}\left(\mathbf{P}^{*};k\right),\left(\mathbf{u}^{a}\right)^{*}\left(\mathbf{P}^{*};k\right)\right\} _{k\in\mathbb{Z}_{0}^{+}}$
is strictly feedback stabilizing, then $\left\{ \left(\mathbf{u}^{c}\right)^{*}\left(\mathbf{P}^{*};k\right),\left(\mathbf{u}^{a}\right)^{*}\left(\mathbf{P}^{*};k\right)\right\} _{k\in\mathbb{Z}_{0}^{+}}$
forms a Nash equilibrium pair in the restricted class of feedback
stabilizing policies {[}10{]}, {[}11{]}, {[}20{]}, {[}48{]}, {[}49{]}.
Note that for deterministic systems, condition $\mathrm{\rho}\left(\mathbb{E}\left[\mathbf{\Xi}\left(k\right)\otimes\mathbf{\Xi}\left(k\right)\right]\right)<1$
in Corollary \ref{Coro: Strictly Stabilizing Class} is reduced to
$\mathrm{\rho}\left(\mathbf{\Xi}\right)<1$ with $\mathbf{\Xi}=\left(\mathbf{I}-\mathbf{B}\Phi^{-1}\left(\mathbf{P}^{*}\right)\mathbf{B}^{T}\mathbf{P}^{*}\right)\mathbf{A}$,
which is equivalent to $\left\{ \left(\mathbf{u}^{c}\right)^{*}\left(\mathbf{P}^{*};k\right),\left(\mathbf{u}^{a}\right)^{*}\left(\mathbf{P}^{*};k\right)\right\} _{k\in\mathbb{Z}_{0}^{+}}$
being strictly feedback stabilizing. Therefore, Corollary \ref{Coro: Strictly Stabilizing Class}
not only encompasses the existing result {[}10{]}, {[}11{]}, {[}20{]},
{[}48{]}, {[}49{]}, but also further extends it to the linear time-variant
stochastic systems with mean-square stability.
\end{remrk}

\section{Numerical Results}

In this section, we verify our developed theoretical results via numerical
simulations

\subsection{Tightness Verifications \label{subsec:Tightness-Verification}}

In this subsection, we demonstrate the tightness of our proposed verifiable
sufficient condition numerically. Specifically, we consider a WNCS
with $\left|\mathcal{N}^{c}\right|=3$ actuators, $\left|\mathcal{N}^{c}\right|=3$
controllers, and $\left|\mathcal{N}^{a}\right|=3$ attackers. The
total numbers of transmit antennas of the controllers and the attackers
are $N_{t}^{c}=3$ and $N_{t}^{a}=2$, respectively. The number of
receive antennas of the actuator is $N_{r}=2$. The plant state dimension
is $S=6$. We randomly pick a full rank symmetric state transition
matrix $\mathbf{A}$ given by {\small{}
\begin{align*}
 & \mathbf{A}=\left[\begin{array}{cccccc}
0.2750 & 0.2745 & 0.2466 & 0.2724 & 0.2516 & 0.2975\\
0.2745 & 0.2862 & 0.2535 & 0.2793 & 0.2450 & 0.2957\\
0.2466 & 0.2535 & 0.2272 & 0.2489 & 0.2188 & 0.2655\\
0.2724 & 0.2793 & 0.2489 & 0.2836 & 0.2495 & 0.2946\\
0.2516 & 0.2450 & 0.2188 & 0.2495 & 0.2412 & 0.2742\\
0.2975 & 0.2957 & 0.2655 & 0.2946 & 0.2742 & 0.3245
\end{array}\right].
\end{align*}
}The spectral radius of $\mathbf{A}$ is $\rho\left(\mathbf{A}\right)=1.6016.$
The control input gain matrix of the actuators are randomly generated
as {\small{}
\begin{align*}
 & \mathbf{B}_{1}=\left[\begin{array}{cccccc}
0.5125 & 0.5750 & 0.0875 & 0.5750 & 0.4000 & 0.0625\\
0.1750 & 0.3500 & 0.600 & 0.6125 & 0.1000 & 0.6125
\end{array}\right]^{T},\\
 & \mathbf{B}_{2}=\left[\begin{array}{cccccc}
0.6000 & 0.3125 & 0.5125 & 0.1000 & 0.2750 & 0.5750\\
0.5000 & 0.6000 & 0.4125 & 0.0250 & 0.5375 & 0.5875
\end{array}\right]^{T},\\
 & \mathbf{B}_{3}=\left[\begin{array}{cccccc}
0.4250 & 0.4750 & 0.4750 & 0.2510 & 0.4125 & 0.1125\\
0.4500 & 0.0251 & 0.1750 & 0.0375 & 0.0625 & 0.5250
\end{array}\right]^{T}.
\end{align*}
The weight matrices of the plant state and the control signals are
chosen as $\mathbf{Q}=\mathbf{I}_{6}$ and $\mathbf{R}_{j}^{c}=\mathbf{I}_{6}$,
$\forall1\leq j\leq3$, respectively. }The wireless MIMO fading channel
between the $j$-th controller and the $i$-th actuator is $\delta\left(k\right)\mathbf{H}_{j,i}^{c}\left(k\right)\in\mathbb{R}^{2\times3}$,
where each controller's random access variable $\delta\left(k\right)$
is Bernoulli distributed with constant mean $\mathbb{E}\left[\delta\left(k\right)\right]=\delta\in\left(0,1\right]$
and is i.i.d. over the time-slots, and each element of $\mathbf{H}_{j,i}^{c}\left(k\right)$
is i.i.d. Gaussian distributed with zero mean and unit variance. The
wireless MIMO fading channel between the $l$-th attacker and the
$i$-th actuator is $\mathbf{H}_{l,i}^{a}\left(k\right)\in\mathbb{R}^{2\times2}$,
where each element of is $\mathbf{H}_{l,i}^{a}\left(k\right)$ is
Gaussian distributed with $\mathbb{E}\left[\mathbf{H}_{l,i}^{a}\left(k\right)\right]=\mathrm{cov}\left(\mathbf{H}_{l,i}^{a}\left(k\right)\right)=\mathbf{I}_{2}$.
These numerical configurations satisfy all the assumptions in Example
\ref{Exmp: E.g.1} and Corollary \ref{Coro: Example 1 Necessary}.

\begin{figure}
\begin{centering}
\subfloat[\label{fig: sim-figure-a}Transitions of the bounds of $\mathrm{Tr}\left(\mathbf{P}^{*}\right)$
to instability.]{\includegraphics[clip,width=0.6\columnwidth]{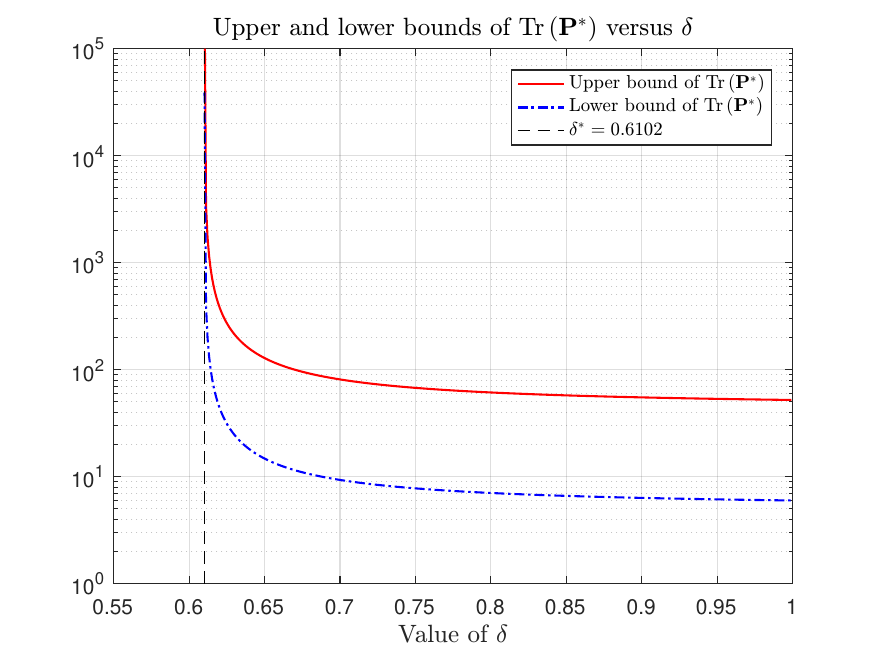}}\ \ \ \subfloat[\label{fig: sim-figure-b}The log-log plot of the lower bounds of
$\left\Vert \mathbf{R}^{a}\right\Vert $ versus $\left(1-\left(1-\delta\right)\rho^{2}\left(\mathbf{A}\right)\right)^{-2}$.]{\includegraphics[clip,width=0.54\columnwidth]{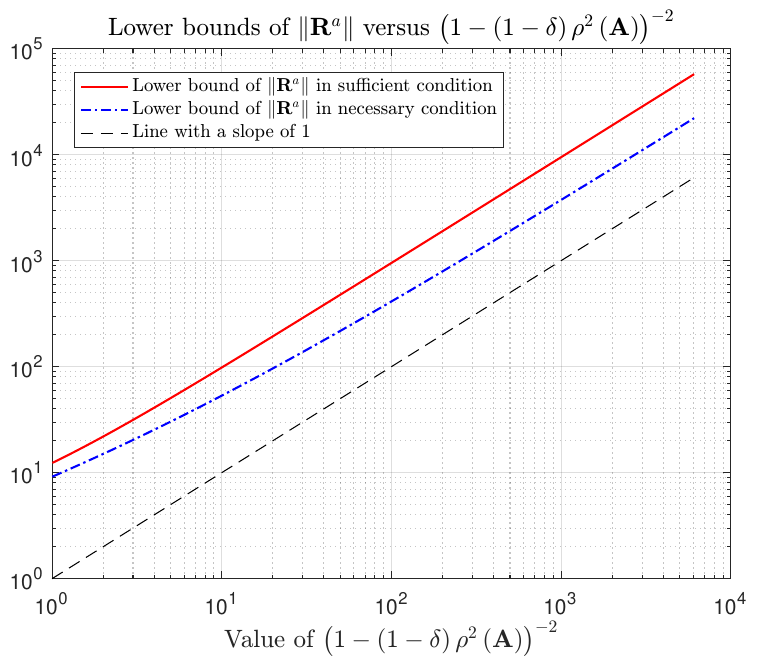}}
\par\end{centering}
\caption{\label{fig: sim-figure}Illustrations of the tightness of the proposed
verifiable sufficient condition.}
\end{figure}

We first demonstrate the tightness of the controllers random access
probability $\delta$ in Fig. \ref{fig: sim-figure-a}. The lower
bound of $\mathrm{Tr}\left(\mathbf{P}^{*}\right)$ in Fig. \ref{fig: sim-figure-a}
is obtained based on inequality (\ref{eq:P* DARE}) in the proof of
Corollary \ref{Coro: Example 1 Necessary}. The upper bound of $\mathrm{Tr}\left(\mathbf{P}^{*}\right)$
in Fig. \ref{fig: sim-figure-a} is calculated based on equality (\ref{eq:T* constructed})
in the proof of Corollary \ref{Coro: Example 1 Suff}. The critical
value of $\delta$ is computed as $\delta^{*}=1-\rho^{-2}\left(\mathbf{A}\right)=1-1.6016^{-2}=0.6102$.
The transition clearly appears in Fig. \ref{fig: sim-figure-a}, where
we see that the value of both upper and lower bound tends to infinity
as $\delta$ approaches $\delta^{*}=0.6102$. As a result, we can
conclude $\delta>\delta^{*}=0.6102$ is both sufficient and necessary
for the boundedness of $\mathrm{Tr}\left(\mathbf{P}^{*}\right)$,
which justifies the tightness of the controllers random access probability
$\delta$ obtained by our proposed verifiable sufficient condition.

The sufficient requirement and necessary requirement on the lower
bound of $\left\Vert \mathbf{R}^{a}\right\Vert $ are plotted against
$\left(1-\left(1-\delta\right)\rho^{2}\left(\mathbf{A}\right)\right)^{-2}$
on a log-log scale in Fig. \ref{fig: sim-figure-b}. Specifically,
the lower bound of $\left\Vert \mathbf{R}^{a}\right\Vert $ in the
proposed sufficient condition is obtained based on equality (\ref{eq:T* constructed})
in the proof of Corollary \ref{Coro: Example 1 Suff} and inequality
(\ref{eq:Ra}) in Theorem \ref{Thm: g(T)<T}. The necessary requirement
on the lower bound of $\left\Vert \mathbf{R}^{a}\right\Vert $ is
obtained based on inequality (\ref{eq:P* DARE}) and (\ref{eq:Ra > Theta-tilde})
in the proof of Corollary \ref{Coro: Example 1 Necessary}. We can
see from Fig. \ref{fig: sim-figure-b} that, as $\left(1-\left(1-\delta\right)\rho^{2}\left(\mathbf{A}\right)\right)^{-2}$
approaches infinity, the slopes of the resulting plots approach 1.
This demonstrates that $\lim_{\delta\rightarrow\delta^{*}}\mathrm{log}\left(\left\Vert \mathbf{R}^{a}\right\Vert \right)>\lim_{\delta\rightarrow\delta^{*}}\mathrm{log}\left(\left(1-\left(1-\delta\right)\rho^{2}\left(\mathbf{A}\right)\right)^{-2}\right)+c$,
where $c$ is some constant representing the $y$-intercept, is both
sufficient and necessary. As a result, our verifiable sufficient condition
provides an $\mathbf{R}^{a}$ with $\left\Vert \mathbf{R}^{a}\right\Vert =\Omega\left(\left(1-\left(1-\delta\right)\rho^{2}\left(\mathbf{A}\right)\right)^{-2}\right),$
which is order-wise optimal

\subsection{Extended Practical Application Scenarios}

We further consider another two typical practical application scenarios
as outlined in Examples \ref{Exmp: E.g.2} and \ref{Exmp: E.g.3}.
Specifically, Example \ref{Exmp: E.g.2} encompasses scenarios where
there is no wireless signal interference between all wireless controllers.
While Example \ref{Exmp: E.g.3} focuses on scenarios where the controllers
and actuators have been equipped with an increased number of transmit
and receive antennas, respectively. The key numerical simulation parameters
for Example \ref{Exmp: E.g.2} is configured as $\left|\mathcal{N}^{c}\right|=\left|\mathcal{N}^{a}\right|=3,$
$N_{t}^{c}=N_{r}=2,$ $S=6$, $\widetilde{\mathbf{B}}_{1}=\left[\begin{array}{cc}
0.5125 & 0.175\\
0.5750 & 0.3500
\end{array}\right],$$\widetilde{\mathbf{B}}_{2}=\left[\begin{array}{cc}
0.5125 & 0.4125\\
0.1000 & 0.0250
\end{array}\right]$, $\widetilde{\mathbf{B}}_{3}=\left[\begin{array}{cc}
0.4125 & 0.0625\\
0.1125 & 0.5250
\end{array}\right]$. The key numerical simulation parameters for Example \ref{Exmp: E.g.3}
is configured as $\left|\mathcal{N}^{c}\right|=\left|\mathcal{N}^{a}\right|=3,$
$N_{t}^{c}=6,$ $N_{r}=2,$ $S=6.$ The configurations of $\mathbf{B}_{1}$,
$\mathbf{B}_{2}$ and $\mathbf{B}_{3}$ remain the same as in Section
\ref{subsec:Tightness-Verification}. In Examples \ref{Exmp: E.g.2}
and \ref{Exmp: E.g.3}, the configurations of the state transition
matrix $\mathbf{A},$ wireless MIMO fading channels $\mathbf{H}_{j,i}^{c}\left(k\right)$
and $\mathbf{H}_{l,i}^{a}\left(k\right)$ remain the same as in Section
\ref{subsec:Tightness-Verification}, and we let $\delta_{l}^{a}=0.5,\forall1\leq l\leq3$
and $\delta_{j}^{c}=\delta,1\leq j\leq3$. The parameter $\varsigma$
is defined as $\text{\ensuremath{\varsigma}}=1-\rho\Big(\mathbb{E}\big[\left(\mathbf{U}_{k}^{T}\mathrm{diag}\left(\left[\boldsymbol{0}_{r_{k}},\mathbf{I}_{S-r_{k}}\right]\right)\mathbf{U}_{k}\mathbf{A}\right)\otimes\big(\mathbf{U}_{k}^{T}\mathrm{diag}\left(\left[\boldsymbol{0}_{r_{k}},\mathbf{I}_{S-r_{k}}\right]\right)\mathbf{U}_{k}\mathbf{A}\big)\big]\Big),$
which represents a margin of the sufficient condition in Theorem \ref{Thm: Verifiable Suff Condition}. 

\begin{figure}
\begin{centering}
\subfloat[\label{fig: sim-figure-a-1} Upper bound of $\mathrm{Tr}\left(\mathbf{P}^{*}\right)$
versus controllers' activation probability.]{\includegraphics[width=0.6\columnwidth]{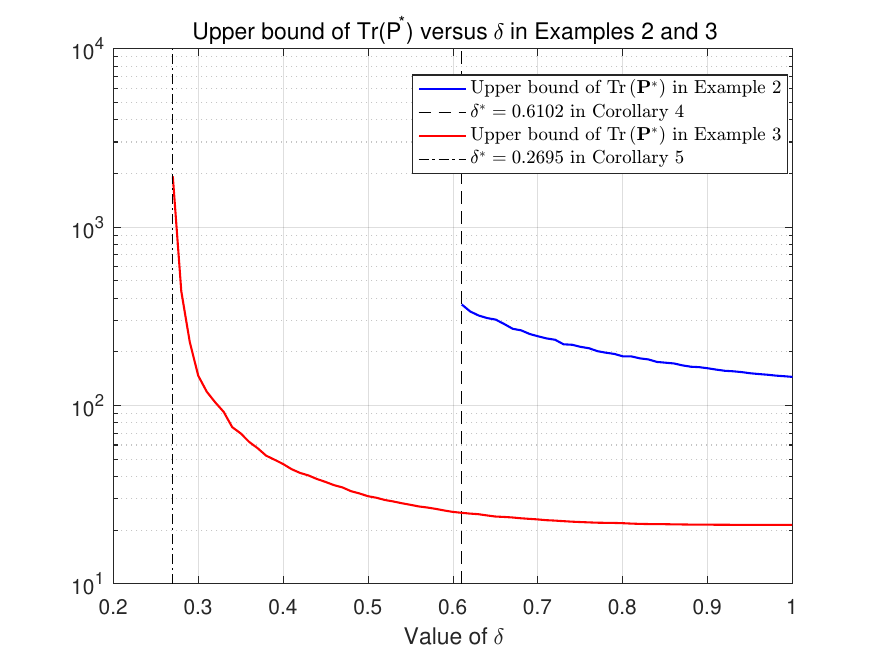}}\ \ \ \subfloat[\label{fig: sim-figure-b-1}The log-log plot of the lower bound of
$\left\Vert \mathbf{R}^{a}\right\Vert $ versus $\varsigma^{-1}$.]{\includegraphics[width=0.6\columnwidth]{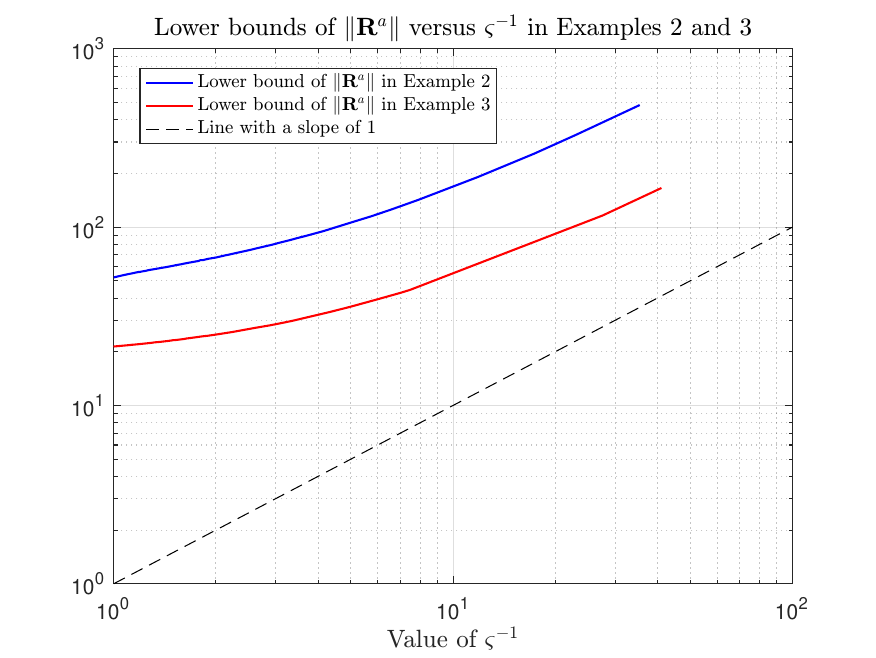}}
\par\end{centering}
\caption{\label{fig: sim-figure-1}Illustrations of the bounds of $\mathrm{Tr}\left(\mathbf{P}^{*}\right)$
and $\left\Vert \mathbf{R}^{a}\right\Vert $ in Examples \ref{Exmp: E.g.2}
and \ref{Exmp: E.g.3}.}
\end{figure}

Figure \ref{fig: sim-figure-a-1} delivers various practical application
insights. Based on Corollaries 4 and 5, to guarantee the existence
of the game value, the controllers' activation probability should
be sufficiently large such that $\delta>\text{\ensuremath{\delta}}^{*},$
where $\text{\ensuremath{\delta}}^{*}$ equals to $0.6012$ and 0.2695
in Examples 2 and 3, respectively. This fact has been demonstrated
in Figure \ref{fig: sim-figure-a-1}, where the finite upper bound
of $\mathrm{Tr}\left(\mathbf{P}^{*}\right)$ indicates the solvability
of the MGARE. Note that $\text{\ensuremath{\delta}}^{*}$ in Example
3 is significantly smaller than that in Example 2, this means that
equipping more transmit antennas at the controllers and more receive
antennas at the actuators is more favorable because it can ensure
the existence of the game value even if the controllers' activation
rate is small. 

The variable $\varsigma$ in Figure \ref{fig: sim-figure-b-1} is
a key parameter that characterizes the margin of the solvability of
the MGARE. The value of $\varsigma$ is small when the controllers'
activation probability $\delta$ is small, where any increase in the
instability of the state transition matrix $\mathbf{A}$ may render
the MGARE unsolvable, ultimately leading to the non-existence of the
game value. Figure \ref{fig: sim-figure-b-1} illustrates that the
logarithm of the lower bound of $\left\Vert \mathbf{R}^{a}\right\Vert $
grows approximately linearly with $\mathrm{log}\left(\varsigma^{-1}\right)$.
This indicates that the magnitude of $\left\Vert \mathbf{R}^{a}\right\Vert $
should be large when the controllers' activation probability $\delta$
is small. To maintain a same penalty in the total cost of the zero-sum
game, a larger $\left\Vert \mathbf{R}^{a}\right\Vert $ will promote
a smaller power of the attack signals, i.e., $\left\Vert \mathbf{u}^{a}\left(k\right)\right\Vert _{2}$,
and this in turn enhances the robustness of the system when controllers
are not active frequently. Note that under the same $\mathrm{log}\left(\varsigma^{-1}\right)$,
$\left\Vert \mathbf{R}^{a}\right\Vert $ in Example 3 is relatively
smaller than that in Example 2, this means that equipping more antennas
at the controllers and the actuators enhances the closed-loop system's
robustness against the attack signals.

\section{Conclusion}

In this paper, we have investigated the LQ zero-sum game for wireless
networked control systems with multi-controllers and multi-attackers,
where the wireless MIMO fading channels and random accesses of the
controllers and attackers will introduce multiplicative randomness
to the closed-loop system. Firstly, we have proposed a general sufficient
and necessary condition, which is difficult to verify. We have next
proposed a novel PSD kernel decomposition induced by multiplicative
randomness. This has facilitated the derivation of a tight closed-form
verifiable sufficient condition for the existence of the game value.
Finally, under the existence condition, we have constructed a saddle-point
policy that is able to achieve the value of the zero-sum game. 

\appendix

\subsection{\label{subsec:Proof-of-Lemma-J_k-representations}Proof of Lemma
\ref{Lemma: J_K representations}}

We first define the per-stage conditional cost $c\left(k\right)$
as
\begin{align}
 & c\left(k\right)=\mathbb{E}\left[\left.\left\Vert \mathbf{x}\left(k+1\right)\right\Vert _{\mathbf{Q}}^{2}+\left\Vert \mathbf{u}^{c}\left(k\right)\right\Vert _{\mathbf{R}^{c}}^{2}-\left\Vert \mathbf{u}^{a}\left(k\right)\right\Vert _{\mathbf{R}^{a}}^{2}\right|\mathcal{I}\left(k\right)\right],\forall0\leq k\leq K-1.\label{eq:per-stage-cost-1}
\end{align}
The finite time horizon cost $J_{K}$ can thus be represented as $J_{K}=\frac{1}{K}\sum_{k=0}^{K-1}\mathbb{E}\left[c\left(k\right)\right]$.

Note that under the necessary and sufficient condition (\ref{eq:Ra-Jk-Condi})
in Lemma \ref{Lemma: J_K representations}, we have, $\forall0\leq k\leq K-1$,
\begin{align}
 & \Phi_{3}\left(f^{K-k-1}\left(\mathbf{Q}\right)\right)\succ\mathbf{0},\ \Phi\left(f^{K-k-1}\left(\mathbf{Q}\right);k\right)\big/\left(-\Phi_{3}\left(f^{K-k-1}\left(\mathbf{Q}\right)\right)\right)\succ\mathbf{0},
\end{align}
where the matrix functions $\Phi\left(\cdot;k\right)$ and $\Phi_{3}\left(\cdot\right)$
are related to the MGARE and are formally defined in Section III-B.

Therefore, substitute the system state dynamics into (\ref{eq:per-stage-cost-1})
and apply the technique of completion of squares w.r.t. $\mathbf{u}_{\Delta}^{c}\left(\mathbf{Q};K-1\right)$,
the average terminal cost $\mathbb{E}\left[c\left(K-1\right)\right]$
can be represented as
\begin{align}
 & \mathbb{E}\left[c\left(K-1\right)\right]\label{eq:c(K-1)-upper-1}\\
 & =\mathbb{E}\Big[\mathbb{E}\Big[\left.\left\Vert \mathbf{x}\left(K\right)\right\Vert _{\mathbf{Q}}^{2}+\left\Vert \mathbf{u}^{c}\left(K-1\right)\right\Vert _{\mathbf{R}^{c}}^{2}-\left\Vert \mathbf{u}^{a}\left(K-1\right)\right\Vert _{\mathbf{R}^{a}}^{2}\right|\mathcal{I}\left(K-1\right)\Big]\Big]\nonumber \\
 & =\mathbb{E}\Big[\left\Vert \mathbf{x}\left(K-1\right)\right\Vert _{\left(f\left(\mathbf{Q}\right)-\mathbf{Q}\right)}^{2}\Big]+\mathrm{Tr}\big(\mathbf{Q}\big(\mathbf{\mathbf{W}}+\sum_{i\in\mathcal{N}^{c}}\mathbf{B}_{i}^{T}\mathbf{V}\mathbf{B}_{i}\big)\big)\nonumber \\
 & +\mathbb{E}\Big[\left\Vert \mathbf{u}_{\Delta}^{c}\left(\mathbf{Q};K-1\right)\right\Vert _{\Phi\left(\mathbf{Q};K-1\right)\big/\left(-\Phi_{3}\left(\mathbf{Q}\right)\right)}^{2}-\big\Vert\mathbf{u}_{\Delta}^{a}\left(\mathbf{Q};K-1\right)+\Phi_{3}^{-1}\left(\mathbf{Q}\right)\Phi_{2}\left(\mathbf{Q};K-1\right)\mathbf{u}_{\Delta}^{c}\left(\mathbf{Q};K-1\right)\big\Vert_{\Phi_{3}\left(\mathbf{Q}\right)}^{2}\Big].\nonumber 
\end{align}

Next consider $\mathbb{E}\left[c\left(K-1\right)+c\left(K-2\right)\right]$,
we have
\begin{align}
 & \mathbb{E}\left[c\left(K-1\right)+c\left(K-2\right)\right]\label{eq:C(K-1)+C(K-2)-1}\\
 & =\mathbf{\mathbb{E}}\Big[\mathbb{E}\Big[\left.\left\Vert \mathbf{x}\left(K-1\right)\right\Vert _{f\left(\mathbf{Q}\right)}^{2}+\left\Vert \mathbf{u}^{c}\left(K-2\right)\right\Vert _{\mathbf{R}^{c}}^{2}-\left\Vert \mathbf{u}^{a}\left(K-2\right)\right\Vert _{\mathbf{R}^{a}}^{2}\right|\mathcal{I}\left(K-2\right)\Big]\Big]+\mathrm{Tr}\big(\mathbf{Q}\big(\mathbf{\mathbf{W}}+\sum_{i\in\mathcal{N}^{c}}\mathbf{B}_{i}^{T}\mathbf{V}\mathbf{B}_{i}\big)\big)\nonumber \\
 & +\mathbb{E}\left[\left\Vert \mathbf{u}_{\Delta}^{c}\left(\mathbf{Q};K-1\right)\right\Vert _{\Phi\left(\mathbf{Q};K-1\right)\big/\left(-\Phi_{3}\left(\mathbf{Q}\right)\right)}^{2}-\left\Vert \mathbf{u}_{\Delta}^{a}\left(\mathbf{Q};K-1\right)+\Phi_{3}^{-1}\left(\mathbf{Q}\right)\Phi_{2}\left(\mathbf{Q};K-1\right)\mathbf{u}_{\Delta}^{c}\left(\mathbf{Q};K-1\right)\right\Vert _{\Phi_{3}\left(\mathbf{Q}\right)}^{2}\right].\nonumber 
\end{align}
Via completing the squares w.r.t. $\mathbf{u}_{\Delta}^{c}\left(f\left(\mathbf{Q}\right);K-2\right)$,
the expectation term $\mathbb{E}\left[\cdot\right]$ in the second
line of (\ref{eq:C(K-1)+C(K-2)-1}) can be represented as 
\begin{align}
 & \mathbf{\mathbb{E}}\Big[\mathbb{E}\Big[\left.\left\Vert \mathbf{x}\left(K-1\right)\right\Vert _{f\left(\mathbf{Q}\right)}^{2}+\left\Vert \mathbf{u}^{c}\left(K-2\right)\right\Vert _{\mathbf{R}^{c}}^{2}-\left\Vert \mathbf{u}^{a}\left(K-2\right)\right\Vert _{\mathbf{R}^{a}}^{2}\right|\mathcal{I}\left(K-2\right)\Big]\Big]\label{eq:C(K-1)+C(K-2)-2-1}\\
 & =\mathbf{\mathbb{E}}\left[\left\Vert \mathbf{x}\left(K-2\right)\right\Vert _{\left(f^{2}\left(\mathbf{Q}\right)-\mathbf{Q}\right)}^{2}\right]+\mathrm{Tr}\big(\mathbf{Q}\big(\mathbf{\mathbf{W}}+\sum_{i\in\mathcal{N}^{c}}\mathbf{B}_{i}^{T}\mathbf{V}\mathbf{B}_{i}\big)\big)\nonumber \\
 & +\mathbb{E}\Big[\left\Vert \mathbf{u}_{\Delta}^{c}\left(f\left(\mathbf{Q}\right);K-2\right)\right\Vert _{\Phi\left(f\left(\mathbf{Q}\right);K-2\right)\big/\left(-\Phi_{3}\left(f\left(\mathbf{Q}\right)\right)\right)}^{2}\nonumber \\
 & -\left\Vert \mathbf{u}_{\Delta}^{a}(f\left(\mathbf{Q}\right);K-2)+\Phi_{3}^{-1}\left(f\left(\mathbf{Q}\right)\right)\Phi_{2}\left(f\left(\mathbf{Q}\right);K-2\right)\mathbf{u}_{\Delta}^{c}\left(f\left(\mathbf{Q}\right);K-2\right)\right\Vert _{\Phi_{3}\left(f\left(\mathbf{Q}\right)\right)}^{2}\Big].\nonumber 
\end{align}
Combining (\ref{eq:C(K-1)+C(K-2)-1}) and (\ref{eq:C(K-1)+C(K-2)-2-1}),
we have 
\begin{align}
 & \mathbb{E}\left[c\left(K-1\right)+c\left(K-2\right)\right]=\mathbb{E}\left[\left\Vert \mathbf{x}\left(K-2\right)\right\Vert _{\left(f^{2}\left(\mathbf{Q}\right)-\mathbf{Q}\right)}^{2}\right]+\mathrm{Tr}\big(\mathbf{Q}\big(\mathbf{\mathbf{W}}+\sum_{i\in\mathcal{N}^{c}}\mathbf{B}_{i}^{T}\mathbf{V}\mathbf{B}_{i}\big)\big)\label{eq:C(K-1)+C(K-2)-final}\\
 & +\mathbb{E}\left[\left\Vert \mathbf{u}_{\Delta}^{c}\left(\mathbf{Q};K-1\right)\right\Vert _{\Phi\left(\mathbf{Q};K-1\right)\big/\left(-\Phi_{3}\left(\mathbf{Q}\right)\right)}^{2}+\left\Vert \mathbf{u}_{\Delta}^{c}\left(f\left(\mathbf{Q}\right);K-2\right)\right\Vert _{\Phi\left(f\left(\mathbf{Q}\right);K-2\right)\big/\left(-\Phi_{3}\left(f\left(\mathbf{Q}\right)\right)\right)}^{2}\right]\nonumber \\
 & -\mathbb{E}\bigg[\left\Vert \mathbf{u}_{\Delta}^{a}\left(\mathbf{Q};K-1\right)+\Phi_{3}^{-1}\left(\mathbf{Q}\right)\Phi_{2}\left(\mathbf{Q};K-1\right)\mathbf{u}_{\Delta}^{c}\left(\mathbf{Q};K-1\right)\right\Vert _{\Phi_{3}\left(\mathbf{Q}\right)}^{2}\nonumber \\
 & +\left\Vert \mathbf{u}_{\Delta}^{a}\left(f\left(\mathbf{Q}\right);K-2\right)+\Phi_{3}^{-1}\left(f\left(\mathbf{Q}\right)\right)\Phi_{2}\left(f\left(\mathbf{Q}\right);K-2\right)\mathbf{u}_{\Delta}^{c}\left(f\left(\mathbf{Q}\right);K-2\right)\right\Vert _{\Phi_{3}\left(f\left(\mathbf{Q}\right)\right)}^{2}\bigg].\nonumber 
\end{align}

We can next proceed to consider the terms $\mathbb{E}\left[\sum_{i=1}^{k}c\left(K-i\right)\right],3\leq k\leq K-1,$
and repeat the steps from (\ref{eq:C(K-1)+C(K-2)-1}) to (\ref{eq:C(K-1)+C(K-2)-final}).
By mathematical induction, we can obtain $\overline{J}_{K}$ in (\ref{eq:J_K-Upper}).

Using similar approaches to complete the squares w.r.t. $\mathbf{u}_{\Delta}^{a}$
(instead of w.r.t. $\mathbf{u}_{\Delta}^{c}$) and note that under
the necessary and sufficient condition (\ref{eq:Ra-Jk-Condi}) in
Lemma \ref{Lemma: J_K representations}, we have
\begin{align}
 & \Phi\left(f^{K-k-1}\left(\mathbf{Q}\right);k\right)\big/\left(\Phi_{1}\left(f^{K-k-1}\left(\mathbf{Q}\right);k\right)\right)\prec\mathbf{0},\forall0\leq k\leq K-1,
\end{align}
where the matrix functions $\Phi\left(\cdot;k\right)$ and $\Phi_{1}\left(\cdot\right)$
are MGARE related and formally defined in Section III-B. We can then
obtain another form of $J_{K}$ as $\underline{J}_{K}$ in (\ref{eq:J_K-Lower})
follow the the same steps in the above derivation of $\overline{J}_{K}$.

\subsection{\label{subsec:Proof-of-Thm General Suff and Necessary-1}Proof of
Theorem \ref{Thm: General case- suff and necessary} }

\textit{Sufficiency}: Based on Lemma \ref{Lemma: P*=00003Df^infty(Q)},
we know that if condition (\ref{eq: general-GARE condition}) is satisfied,
then condition (\ref{eq:Ra-Jk-Condi}) in Lemma \ref{Lemma: J_K representations}
is satisfied for all $K\in\mathbb{Z}^{+}.$ This means that under
condition (\ref{eq: general-GARE condition}), the representation
of $J_{K}$ as $\overline{J}_{K}$ or $\underline{J}_{K}$ holds true
for all $K\in\mathbb{Z}^{+}.$ 

As a result, under condition (\ref{eq: general-GARE condition}),
we have, $\forall K\in\mathbb{Z}^{+},$
\begin{align}
 & \inf_{\left\{ \boldsymbol{\mu}_{l}^{c}\right\} _{j\in\mathcal{N}^{c}}\in\left(\mathcal{U}^{c}\right)^{\left|\mathcal{N}^{c}\right|}}\sup_{\left\{ \boldsymbol{\mu}_{l}^{a}\right\} _{l\in\mathcal{N}^{a}}\in\left(\mathcal{U}^{a}\right)^{\left|\mathcal{N}^{a}\right|}}J_{K}=\label{eq: min max J-K}\\
 & \min_{\left\{ \mathbf{u}_{\Delta}^{c}\left(f^{K-k-1}\left(\mathbf{Q}\right);k\right)\right\} _{k\in\left[0,K-1\right]}}\max_{\left\{ \left(\mathbf{u}_{\Delta}^{a}\right)\left(f^{K-k-1}\left(\mathbf{Q}\right);k\right)\right\} _{k\in\left[0,K-1\right]}}\overline{J}_{K}\nonumber \\
 & \overset{\left(\ref{eq: min max J-K}.a\right)}{=}\frac{1}{K}\left\Vert \mathbf{x}\left(0\right)\right\Vert _{\left(f^{K}\left(\mathbf{Q}\right)-\mathbf{Q}\right)}^{2}\nonumber \\
 & +\mathrm{Tr}\left(\frac{1}{K}\sum_{k=0}^{K-1}\left(f^{k}\left(\mathbf{Q}\right)\right)\left(\mathbf{\mathbf{W}}+\sum_{i\in\mathcal{N}^{c}}\mathbf{B}_{i}^{T}\mathbf{V}\mathbf{B}_{i}\right)\right).\nonumber 
\end{align}

The equality ($\ref{eq: min max J-K}.a$) holds because of the following
reasons. Since $\overline{J}_{K}$ in (\ref{eq:J_K-Upper}) is a summation
of nonnegative quadratic terms, we have 
\begin{align}
 & \max_{\left\{ \left(\mathbf{u}_{\Delta}^{a}\right)\left(f^{K-k-1}\left(\mathbf{Q}\right);k\right)\right\} _{k\in\left[0,K-1\right]}}\overline{J}_{K}\overset{\left(\ref{eq:max-J_K-Upper}.a\right)}{=}\frac{1}{K}\left\Vert \mathbf{x}\left(0\right)\right\Vert _{\left(f^{K}\left(\mathbf{Q}\right)-\mathbf{Q}\right)}^{2}\label{eq:max-J_K-Upper}\\
 & +\mathrm{Tr}\left(\frac{1}{K}\sum_{k=0}^{K-1}\left(f^{k}\left(\mathbf{Q}\right)\right)\left(\mathbf{\mathbf{W}}+\sum_{i\in\mathcal{N}^{c}}\mathbf{B}_{i}^{T}\mathbf{V}\mathbf{B}_{i}\right)\right)+\frac{1}{K}\sum_{k=0}^{K-1}\nonumber \\
 & \mathbb{E}\left[\left\Vert \mathbf{u}_{\Delta}^{c}\left(f^{K-k-1}\left(\mathbf{Q}\right);k\right)\right\Vert _{\Phi\left(f^{K-k-1}\left(\mathbf{Q}\right);k\right)\big/\left(-\Phi_{3}\left(f^{K-k-1}\left(\mathbf{Q}\right)\right)\right)}^{2}\right],\nonumber 
\end{align}
where ($\ref{eq:max-J_K-Upper}.a$) holds by choosing, $\forall0\leq k\leq K-1,$
\begin{align}
 & \left(\mathbf{u}_{\Delta}^{a}\right)^{*}\left(f^{K-k-1}\left(\mathbf{Q}\right);k\right)=-\Phi_{3}^{-1}\left(f^{K-k-1}\left(\mathbf{Q}\right)\right)\label{eq:U^a_delta* Solution}\\
 & \cdot\Phi_{2}\left(f^{K-k-1}\left(\mathbf{Q}\right);k\right)\mathbf{u}_{\Delta}^{c}\left(f^{K-k-1}\left(\mathbf{Q}\right);k\right),\nonumber 
\end{align}
or equivalently (based on the relationship between $\mathbf{u}_{\Delta}^{a}$
and $\mathbf{u}^{a}$ in (\ref{eq:delta-u-a})), $\forall0\leq k\leq K-1,$
\begin{align}
 & \left(\mathbf{u}^{a}\left(k\right)\right)^{*}=-\Phi_{3}^{-1}\left(f^{K-k-1}\left(\mathbf{Q}\right)\right)\label{eq:U^a*-Solution}\\
 & \cdot\Phi_{2}\left(f^{K-k-1}\left(\mathbf{Q}\right);k\right)\mathbf{u}_{\Delta}^{c}\left(f^{K-k-1}\left(\mathbf{Q}\right);k\right)\nonumber \\
 & -\left(\Phi^{-1}\left(\mathbf{P};k\right)\mathbf{B}^{T}\left(k\right)\mathbf{P}\mathbf{A}\mathbf{x}\left(k\right)\right)_{\left(N_{t}^{c}\left|\mathcal{N}^{c}\right|+1\right):\left(N_{t}^{c}\left|\mathcal{N}^{c}\right|+N_{t}^{a}\left|\mathcal{N}^{a}\right|\right)}.\nonumber 
\end{align}
Then taking the minimization w.r.t. $\mathbf{u}_{\Delta}^{c}$ on
both sides of (\ref{eq:max-J_K-Upper}), we can obtain equality ($\ref{eq: min max J-K}.a$),
and the corresponding minimizing solution $\left(\mathbf{u}_{\Delta}^{c}\right)^{*}$
is, $\forall0\leq k\leq K-1,$ 
\begin{align}
 & \left(\mathbf{u}_{\Delta}^{c}\right)^{*}\left(f^{K-k-1}\left(\mathbf{Q}\right);k\right)=\mathbf{0},\label{eq:U-delat-c^*}
\end{align}
or equivalently (based on the relationship between $\mathbf{u}_{\Delta}^{c}$
and $\mathbf{u}^{c}$ in (\ref{eq:delta-u-c})), $\forall0\leq k\leq K-1,$
\begin{align}
\left(\mathbf{u}^{c}\left(k\right)\right)^{*}=- & \left(\Phi^{-1}\left(\mathbf{P};k\right)\mathbf{B}^{T}\left(k\right)\mathbf{P}\mathbf{A}\mathbf{x}\left(k\right)\right)_{1:N_{t}^{c}\left|\mathcal{N}^{c}\right|}.\label{eq:u-c^*}
\end{align}
Substitute (\ref{eq:U-delat-c^*}) back to (\ref{eq:U^a*-Solution}),
we have,$,\forall0\leq k\leq K-1,$ 
\begin{align}
 & \left(\mathbf{u}^{a}\left(k\right)\right)^{*}=\label{eq:u-a^*}\\
 & -\left(\Phi^{-1}\left(\mathbf{P};k\right)\mathbf{B}^{T}\left(k\right)\mathbf{P}\mathbf{A}\mathbf{x}\left(k\right)\right)_{\left(N_{t}^{c}\left|\mathcal{N}^{c}\right|+1\right):\left(N_{t}^{c}\left|\mathcal{N}^{c}\right|+N_{t}^{a}\left|\mathcal{N}^{a}\right|\right)}.\nonumber 
\end{align}
Therefore, we conclude that (\ref{eq:u-c^*}) and (\ref{eq:u-a^*})
achieve the upper value of $J_{K}$.

Similarly, under condition (\ref{eq: general-GARE condition}), we
also have, $\forall K\in\mathbb{Z}^{+},$
\begin{align}
 & \sup_{\left\{ \boldsymbol{\mu}_{l}^{a}\right\} _{l\in\mathcal{N}^{a}}\in\left(\mathcal{U}^{a}\right)^{\left|\mathcal{N}^{a}\right|}}\inf_{\left\{ \boldsymbol{\mu}_{l}^{c}\right\} _{j\in\mathcal{N}^{c}}\in\left(\mathcal{U}^{c}\right)^{\left|\mathcal{N}^{c}\right|}}J_{K}\label{eq: max min J-K}\\
 & =\max_{\left\{ \mathbf{u}_{\Delta}^{a}\left(f^{K-k-1}\left(\mathbf{Q}\right);k\right)\right\} _{k\in\left[0,K-1\right]}}\min_{\left\{ \mathbf{u}_{\Delta}^{c}\left(f^{K-k-1}\left(\mathbf{Q}\right);k\right)\right\} _{k\in\left[0,K-1\right]}}\underline{J}_{K}\nonumber \\
 & \overset{\left(\ref{eq: max min J-K}.a\right)}{=}\frac{1}{K}\left\Vert \mathbf{x}\left(0\right)\right\Vert _{\left(f^{K}\left(\mathbf{Q}\right)-\mathbf{Q}\right)}^{2}\nonumber \\
 & +\mathrm{Tr}\left(\frac{1}{K}\sum_{k=0}^{K-1}\left(f^{k}\left(\mathbf{Q}\right)\right)\left(\mathbf{\mathbf{W}}+\sum_{i\in\mathcal{N}^{c}}\mathbf{B}_{i}^{T}\mathbf{V}\mathbf{B}_{i}\right)\right),
\end{align}
where the equality ($\ref{eq: max min J-K}.a$) holds by choosing
$\mathbf{u}^{c}\left(k\right)=\left(\mathbf{u}^{c}\left(k\right)\right)^{*}$
in (\ref{eq:u-c^*}) and $\mathbf{u}^{a}\left(k\right)=\left(\mathbf{u}^{a}\left(k\right)\right)^{*}$
in (\ref{eq:u-a^*}).

Note that condition (\ref{eq: general-GARE condition}) is sufficient
to guarantee 
\begin{align}
 & \lim_{K\rightarrow\infty}\frac{1}{K}\sum_{k=0}^{K-1}f^{k}\left(\mathbf{Q}\right)=f^{\infty}\left(\mathbf{Q}\right),\\
 & \lim_{K\rightarrow\infty}\frac{1}{K}\left\Vert \mathbf{x}\left(0\right)\right\Vert _{\left(f^{K}\left(\mathbf{Q}\right)-\mathbf{Q}\right)}^{2}=0.
\end{align}

Combining (\ref{eq: min max J-K}), (\ref{eq: max min J-K}) and letting
$K\rightarrow\infty$, we conclude that, under condition (\ref{eq: general-GARE condition}),
the zero-sum game has equal upper value and lower value, i.e., $\overline{J}=\underline{J}=J^{*}=\mathrm{Tr}\left(f^{\infty}\left(\mathbf{Q}\right)\left(\mathbf{\mathbf{W}}+\sum_{i\in\mathcal{N}^{c}}\mathbf{B}_{i}^{T}\mathbf{V}\mathbf{B}_{i}\right)\right).$
This completes the proof of sufficiency.

\textit{Necessity}: Let us first consider the finite time horizon
case with horizon length being $K$. For any given $K\in\mathbb{Z}^{+},$
we construct a $K$-length sequence of indicator functions $\left\{ \mathds{1}_{\mathbf{R}^{a}}^{K}\left(k\right)\right\} _{k\in\left[0,K-1\right]},$
where 
\begin{align}
 & \mathds{1}_{\mathbf{R}^{a}}^{K}\left(k\right)=\mathds{1}_{\left\{ \mathbf{R}^{a}\succ\mathbb{E}\left[\left(\mathbf{B}^{a}\left(k\right)\right)^{T}f^{k}\left(\mathbf{Q}\right)\mathbf{B}^{a}\left(k\right)\right]\right\} }.
\end{align}
We next prove by contradiction that if the value of the $K$-length
time horizon game exists, then $\mathds{1}_{\mathbf{R}^{a}}^{K}\left(k\right)\neq0,$
$\forall k\in\left[0,K-1\right],\forall K\in\mathbb{Z}^{+}$.

For some $K\in\mathbb{Z}^{+}$, we suppose the value of the $K$-length
time horizon game exists, and there also exists at least a $k\in\left[0,K-1\right]$
such that $\mathds{1}_{\mathbf{R}^{a}}^{K}\left(k\right)=0$. Let
\begin{align}
 & \tilde{k}=\min\left\{ k:\mathds{1}_{\mathbf{R}^{a}}^{K}\left(k\right)=0,k\in\left[0,K-1\right]\right\} .
\end{align}
According to Lemma \ref{Lemma: J_K representations}, we have 
\begin{align}
 & \sum_{k=0}^{\tilde{k}}\mathbb{E}\left[c\left(K-k-1\right)\right]=\mathbb{E}\left[\left\Vert \mathbf{x}(K-\tilde{k})\right\Vert _{\left(f^{\tilde{k}}\left(\mathbf{Q}\right)-\mathbf{Q}\right)}^{2}\right]\\
 & +\sum_{k=0}^{\tilde{k}}\mathrm{Tr}\left(f^{k}\left(\mathbf{Q}\right)\left(\mathbf{\mathbf{W}}+\sum_{i\in\mathcal{N}^{c}}\mathbf{B}_{i}^{T}\mathbf{V}\mathbf{B}_{i}\right)\right)\nonumber \\
 & +\sum_{k=0}^{\tilde{k}}\mathbb{E}\left[\vphantom{\left[\left(\mathbf{u}_{\Delta}^{a}\left(f^{k}\left(\mathbf{Q}\right);K-k-1\right)\right)^{T}\right]}\right.\left(\mathbf{u}_{\Delta}^{a}\left(f^{k}\left(\mathbf{Q}\right);K-k-1\right)\right)^{T}\nonumber \\
 & \cdot\left(\Phi\left(f^{k}\left(\mathbf{Q}\right);K-k-1\right)\big/\Phi_{1}\left(f^{k}\left(\mathbf{Q}\right);K-k-1\right)\right)\nonumber \\
 & \cdot\mathbf{u}_{\Delta}^{a}\left(f^{k}\left(\mathbf{Q}\right);K-k-1\right)\left.\vphantom{\left[\left(\mathbf{u}_{\Delta}^{a}\left(f^{k}\left(\mathbf{Q}\right);K-k-1\right)\right)^{T}\right]}\right]\nonumber \\
 & +\sum_{k=0}^{\tilde{k}}\mathbb{E}\left[\vphantom{\left[\left\Vert \mathbf{x}\left(K-1\right)\right\Vert _{\left(f\left(\mathbf{Q}\right)-\mathbf{Q}\right)}^{2}\right]}\right.\left\Vert \vphantom{\left\Vert \mathbf{u}_{\Delta}^{c}\left(\mathbf{Q};K-1\right)\right\Vert }\right.\mathbf{u}_{\Delta}^{c}\left(f^{k}\left(\mathbf{Q}\right);K-k-1\right)\nonumber \\
 & +\Phi_{1}^{-1}\left(f^{k}\left(\mathbf{Q}\right);K-k-1\right)\Phi_{2}\left(f^{k}\left(\mathbf{Q}\right);K-k-1\right)\nonumber \\
 & \cdot\mathbf{u}_{\Delta}^{a}\left(f^{k}\left(\mathbf{Q}\right);K-k-1\right)\left.\vphantom{\left\Vert \mathbf{u}_{\Delta}^{c}\left(\mathbf{Q};K-1\right)\right\Vert }\right\Vert _{\Phi_{1}\left(f^{k}\left(\mathbf{Q}\right);K-k-1\right)}^{2}\left.\vphantom{\vphantom{\left[\left\Vert \mathbf{x}\left(K-1\right)\right\Vert _{\left(f\left(\mathbf{Q}\right)-\mathbf{Q}\right)}^{2}\right]}}\right].\nonumber 
\end{align}

The Hessian matrix of $\sum_{k=0}^{\tilde{k}}\mathbb{E}\left[c\left(K-k-1\right)\right]$
at $\mathbf{u}_{\Delta}^{a}\left(f^{k}\left(\mathbf{Q}\right);K-k-1\right)$
is 
\begin{align}
 & \frac{\partial^{2}\sum_{k=0}^{\tilde{k}}\mathbb{E}\left[c\left(K-k-1\right)\right]}{\partial^{2}\mathbf{u}_{\Delta}^{a}\left(f^{\tilde{k}}\left(\mathbf{Q}\right);K-\tilde{k}-1\right)}\label{eq:Hessian}\\
 & =\Phi\left(f^{\tilde{k}}\left(\mathbf{Q}\right);K-\tilde{k}-1\right)\big/\Phi_{1}\left(f^{\tilde{k}}\left(\mathbf{Q}\right);K-\tilde{k}-1\right)\nonumber \\
 & +\Phi_{2}\left(f^{\tilde{k}}\left(\mathbf{Q}\right);K-\tilde{k}-1\right)\Phi_{1}^{-1}\left(f^{\tilde{k}}\left(\mathbf{Q}\right);K-\tilde{k}-1\right)\nonumber \\
 & \cdot\Phi_{2}\left(f^{\tilde{k}}\left(\mathbf{Q}\right);K-\tilde{k}-1\right)=-\mathbf{R}^{a}\nonumber \\
 & +\mathbb{E}\left[\Big(\mathbf{B}^{a}\left(K-\tilde{k}-1\right)\Big)^{T}f^{\tilde{k}}\left(\mathbf{Q}\right)\mathbf{B}^{a}\Big(K-\tilde{k}-1\Big)\right]\nprec\mathbf{0}.\nonumber 
\end{align}

Since the Hessian matrix in (\ref{eq:Hessian}) is indefinite, under
any given $\tilde{k}$-length control action sequence $\left\{ \mathbf{u}_{\Delta}^{c}\left(f^{k}\left(\mathbf{Q}\right);K-k-1\right)\right\} _{k\in\left[0,\tilde{k}\right]}$
with $0<\tilde{k}<K$ (or equivalently $\left\{ \mathbf{u}^{c}\left(K-k-1\right)\right\} _{k\in\left[0,\tilde{k}\right]}$),
we have
\begin{align}
 & \max_{\left\{ \mathbf{u}^{a}\left(K-k-1\right)\right\} _{k\in\left[0,\tilde{k}\right]}}\sum_{k=0}^{\tilde{k}}\mathbb{E}\left[c\left(K-k-1\right)\right]=\infty.
\end{align}
 This means that the upper value of the $K$-length time horizon game
does not exist, which contradicts with our assumption that the value
of the game exists.

As a result, we can conclude that, for any $K\in\mathbb{Z}^{+},$
if value of the $K$-length time horizon game exists, then 
\begin{align}
 & \mathds{1}_{\mathbf{R}^{a}}^{K}\left(k\right)\neq0,\forall k\in\left[0,K-1\right].\label{eq:condi-R-K-ra}
\end{align}

Condition (\ref{eq:condi-R-K-ra}) coincides with condition (\ref{eq:Ra-Jk-Condi})
in Lemma \ref{Lemma: J_K representations}. Therefore, based on Lemma
\ref{Lemma: J_K representations} and letting $K\rightarrow\infty,$
it is clear that if value of the infinite time horizon game exists,
$\lim_{K\rightarrow\infty}\frac{1}{K}\sum_{k=0}^{K-1}\left(f^{k}\left(\mathbf{Q}\right)\right)$
must be bounded. Further note that, under condition (\ref{eq:condi-R-K-ra}),
$f^{k}\left(\mathbf{Q}\right)$ is monotonically increasing as $k$
increases, $\forall k\in\mathbb{Z}^{+}$. We can conclude that $f^{\infty}\left(\mathbf{Q}\right)$
exists and is bounded. 

The existence of $f^{\infty}\left(\mathbf{Q}\right)$ implies that
$f^{\infty}\left(\mathbf{Q}\right)=f\left(f^{\infty}\left(\mathbf{Q}\right)\right)$,
which means that $f^{\infty}\left(\mathbf{Q}\right)\in\mathcal{P}_{\mathbf{R}^{a}}$.
Let $k\rightarrow\infty$, the existence of $f^{\infty}\left(\mathbf{Q}\right)$
together with condition (\ref{eq:condi-R-K-ra}) also imply that $f^{\infty}\left(\mathbf{Q}\right)\in\mathcal{R}_{\mathbf{R}^{a}}$.
Therefore, we conclude that if value of the game exists, then $f^{\infty}\left(\mathbf{Q}\right)\in\mathcal{P}_{\mathbf{R}^{a}}\cap\mathcal{R}_{\mathbf{R}^{a}}\neq\textrm{Ø}.$
This completes the proof of necessity.

\subsection{\label{subsec:Proof-of-Corollary: Refined Suff=000026Necessary}Proof
of Corollary \ref{Corollary: Refined Suff=000026Necessary}}

\textsl{Sufficiency}: Based on Lemma \ref{Lemma: Monotonicity f(.)},
we know that if $f^{k}\left(\mathbf{Q}\right)\in\mathcal{R}_{\mathbf{R}^{a}},\forall k\in\mathbb{Z}_{0}^{+},$
then $\left\{ f^{k}\left(\mathbf{Q}\right)\right\} _{k\in\mathbb{Z}_{0}^{+}}$is
a monotonically increasing sequence. Since $f^{k}\left(\mathbf{Q}\right)$
is bounded, $\forall k\in\mathbb{Z}_{0}^{+},$ we obtain that $\left\{ f^{k}\left(\mathbf{Q}\right)\right\} _{k\in\mathbb{Z}_{0}^{+}}$
is convergent. Therefore, $f^{\infty}\left(\mathbf{Q}\right)$ exists
and $f^{\infty}\left(\mathbf{Q}\right)\in\mathcal{R}_{\mathbf{R}^{a}},$
which means that $f^{\infty}\left(\mathbf{Q}\right)\in\mathcal{P}_{\mathbf{R}^{a}}\cap\mathcal{R}_{\mathbf{R}^{a}}$
and the general sufficient and necessary condition (\ref{eq: general-GARE condition})
is satisfied.

\textsl{Necessity}: Based on Lemma \ref{Lemma: P*=00003Df^infty(Q)},
if $\mathcal{P}_{\mathbf{R}^{a}}\cap\mathcal{R}_{\mathbf{R}^{a}}\neq\textrm{Ø}$,
then $\mathbf{Q}\prec f^{\infty}\left(\mathbf{Q}\right)$ and $\mathbf{Q}\in\mathcal{R}_{\mathbf{R}^{a}}.$
Continuously applying the operator $f$ to $\mathbf{Q}\prec f^{\infty}\left(\mathbf{Q}\right)$
and using the monotonicity in Lemma \ref{Lemma: Monotonicity f(.)},
we have, $\forall k\in\mathbb{Z}_{0}^{+},$
\begin{align}
 & \mathbf{Q}\prec f\left(\mathbf{Q}\right)\prec\cdots\prec f^{k}\left(\mathbf{Q}\right)\prec f^{k+1}\left(\mathbf{Q}\right)\prec\label{eq:increasing f(lambda, P)-2}\\
 & \cdots\prec f\left(f^{\infty}\left(\mathbf{Q}\right)\right)=f^{\infty}\left(\mathbf{Q}\right).\nonumber 
\end{align}
Therefore, $f^{k}\left(\mathbf{Q}\right)$ is bounded, $\forall k\in\mathbb{Z}_{0}^{+}.$
Since $f^{\infty}\left(\mathbf{Q}\right)\in\mathcal{R}_{\mathbf{R}^{a}},$
(\ref{eq:increasing f(lambda, P)-2}) implies $f^{k}\left(\mathbf{Q}\right)\in\mathcal{R}_{\mathbf{R}^{a}},\forall k\in\mathbb{Z}_{0}^{+}.$
This completes the proof of Corollary \ref{Corollary: Refined Suff=000026Necessary}.

\subsection{\label{subsec:Proof-of-Thm: g(T)<T}Proof of Theorem \ref{Thm: g(T)<T}}

For any $\mathbf{T}\in\mathbb{S}_{+}^{S}$ that satisfies condition
(\ref{eq: g(T)<T}), we can construct a matrix pair $\left(\widetilde{\mathbf{P}},\mathbf{R}^{a}\right)$,
where $\widetilde{\mathbf{P}}\in\mathbb{S}_{+}^{S}$ and $\mathbf{R}^{a}\succ\mathrm{Cov}\left(\mathbf{B}^{a}\left(k\right)\mathbf{\widetilde{\mathbf{P}}}^{\frac{1}{2}}\right)$,
such that 
\begin{align}
 & \widetilde{\mathbf{P}}^{-1}=\mathbb{E}\left[\mathbf{B}^{a}\left(k\right)\right]\left(\mathbf{R}^{a}-\mathrm{cov}\left(\mathbf{B}^{a}\left(k\right)\widetilde{\mathbf{P}}^{\frac{1}{2}}\right)\right)^{-1}\label{eq:inv P-tilde}\\
 & \cdot\mathbb{E}\left[\left(\mathbf{B}^{a}\left(k\right)\right)^{T}\right]+\mathbf{T}^{-1}.\nonumber 
\end{align}
By the definition of $\mathcal{R}_{\mathbf{R}^{a}}$ in (\ref{eq: Definition R^a-set}),
it is clear that $\widetilde{\mathbf{P}}$ defined in (\ref{eq:inv P-tilde})
satisfies $\widetilde{\mathbf{P}}\in\mathcal{R}_{\mathbf{R}^{a}}.$
We next analyze the requirement on $\mathbf{R}^{a}$ such that $\widetilde{\mathbf{P}}$
constructed in (\ref{eq:inv P-tilde}) also satisfies $f\left(\widetilde{\mathbf{P}}\right)\prec\mathbf{\widetilde{\mathbf{P}}}.$

We now defined a time varying matrix function $g\left(\mathbf{T};k\right)$
as follows
\begin{align}
 & g\left(\mathbf{T};k\right)\triangleq\mathbf{A}^{T}\left(\mathbf{B}^{c}\left(k\right)\left(\mathbf{R}^{c}\right)^{-1}\left(\mathbf{B}^{c}\left(k\right)\right)^{T}+\mathbf{T}^{-1}\right)^{-1}\mathbf{A},
\end{align}
which will be later frequently used for compact presentations of the
proofs. It is important to note that $g\left(\mathbf{T};k\right)$
is a random matrix function with randomness induced by $\mathbf{B}^{c}\left(k\right)$.

Substitute $\widetilde{\mathbf{P}}^{-1}$ in (\ref{eq:inv P-tilde})
into the MGARE (\ref{eq: general-GARE}), we obtain 
\begin{align}
 & f\left(\widetilde{\mathbf{P}}\right)=\mathbb{E}\left[g\left(\mathbf{T};k\right)\right]+\mathbf{Q}.\label{eq:f(P-tilde)}
\end{align}
Combining (\ref{eq:inv P-tilde}) and (\ref{eq:f(P-tilde)}), $f\big(\widetilde{\mathbf{P}}\big)\prec\mathbf{\widetilde{\mathbf{P}}}$
is equivalent to 
\begin{align}
 & \left(\mathbb{E}\left[g\left(\mathbf{T};k\right)\right]+\mathbf{Q}\right)^{-1}-\mathbf{T}^{-1}\label{eq:f(P)<P equivalent}\\
 & \succ\mathbb{E}\left[\mathbf{B}^{a}\left(k\right)\right]\left(\mathbf{R}^{a}-\mathrm{cov}\left(\mathbf{B}^{a}\left(k\right)\widetilde{\mathbf{P}}^{\frac{1}{2}}\right)\mathbb{E}\right)^{-1}\left[\left(\mathbf{B}^{a}\left(k\right)\right)^{T}\right].\nonumber 
\end{align}
The left hand side of (\ref{eq:f(P)<P equivalent}) is positive definite
because of condition (\ref{eq: g(T)<T}). Condition (\ref{eq:f(P)<P equivalent})
is further equivalent to
\begin{align}
 & \mathbf{R}^{a}\succ\mathbb{E}\left[\left(\mathbf{B}^{a}\left(k\right)\right)^{T}\right]\left(\left(\mathbb{E}\left[g\left(\mathbf{T};k\right)\right]+\mathbf{Q}\right)^{-1}-\mathbf{T}^{-1}\right)^{-1}\label{eq:Ra-lower bound}\\
 & \cdot\mathbb{E}\left[\mathbf{B}^{a}\left(k\right)\right]+\mathrm{cov}\left(\mathbf{B}^{a}\left(k\right)\widetilde{\mathbf{P}}^{\frac{1}{2}}\right).\nonumber 
\end{align}

Let $\sigma_{\mathbf{B}^{a}}^{2}$ denote the maximum value of the
covariance of each element in matrix $\mathbf{B}^{a}\left(k\right)$,
i.e., $\sigma_{\mathbf{B}^{a}}^{2}=\max\left\{ \mathrm{Cov}\left(\left(\mathbf{B}^{a}\left(k\right)\right)_{i,j}\right),1\leq i\leq S,1\leq j\leq\left(N_{t}^{a}\left|\mathcal{N}^{a}\right|\right)\right\} .$
By the construction of $\widetilde{\mathbf{P}}$ in (\ref{eq:inv P-tilde}),
we have $\widetilde{\mathbf{P}}\prec\mathbf{T}$ and
\begin{align}
 & \mathrm{cov}\left(\mathbf{B}^{a}\left(k\right)\widetilde{\mathbf{P}}^{\frac{1}{2}}\right)\prec\sigma_{\mathbf{B}^{a}}^{2}\mathrm{Tr}\left(\widetilde{\mathbf{P}}\right)\mathbf{I}\prec\sigma_{\mathbf{B}^{a}}^{2}\mathrm{Tr}\left(\mathbf{T}\right)\mathbf{I}.
\end{align}
As a result, the choice of $\mathbf{R}^{a}$ in (\ref{eq:Ra}) is
sufficient to guarantee that inequality (\ref{eq:f(P)<P equivalent})
and (\ref{eq:Ra-lower bound}) are both satisfied, which in turn guarantees
$f\left(\widetilde{\mathbf{P}}\right)\prec\mathbf{\widetilde{\mathbf{P}}}$is
satisfied. Therefore, we conclude that $\left(\widetilde{\mathbf{P}},\mathbf{R}^{a}\right)$
given by (\ref{eq:P-tilde}) and (\ref{eq:Ra}) satisfies condition
(\ref{eq:condition P-tilde}) in Lemma \ref{Lemma: P-tilde bound}.

\subsection{\label{subsec:Proof-of-Lemma: Suff-Condition-Existence of T}Proof
of Lemma \ref{Lemma: Suff-Condition-Existence of T}}

We prove Lemma \ref{Lemma: Suff-Condition-Existence of T} by applying
the proposed closed-form PSD kernel decomposition in Theorem \ref{Thm: Ker-Decomposition}
to the matrix inequalities in Theorem \ref{Thm: g(T)<T}. 

Note that the L.H.S. of (\ref{eq: g(T)<T}) can be represented as
\begin{align}
 & \mathbb{E}\left[g\left(\mathbf{T};k\right)\right]+\mathbf{Q}=\mathbf{A}^{T}\mathbb{E}\bigg[\mathbf{T}-\mathbf{T}\mathbf{B}^{c}\left(k\right)\label{eq:expanded g(T)}\\
 & \cdot\left(\mathbf{R}^{c}+\left(\mathbf{B}^{c}\left(k\right)\right)^{T}\mathbf{T}\mathbf{B}^{c}\left(k\right)\right)^{-1}\left(\mathbf{B}^{c}\left(k\right)\right)^{T}\mathbf{T}\bigg]\mathbf{A}+\mathbf{Q}.\nonumber 
\end{align}
Substitute $\mathbf{T}=\mathbf{T}_{k}^{\mathrm{ker}}+\mathbf{T}_{k}^{\mathrm{\mathbf{0}}}$
into (\ref{eq:expanded g(T)}) and note that $\mathbf{T}_{k}^{\mathrm{\mathbf{0}}}\mathbf{B}^{c}\left(k\right)=\mathbf{0}$,
we obtain 
\begin{align}
 & \mathbb{E}\left[g\left(\mathbf{T};k\right)\right]+\mathbf{Q}\label{eq:Expanded substitute}\\
 & \preceq\mathrm{Pr}\left(r_{k}=0\right)\mathbf{A}^{T}\mathbf{T}\mathbf{A}+\mathbb{E}\left[\mathds{1}_{\left\{ r_{k}\neq0,\left(\boldsymbol{\Sigma}_{k}\right)_{r_{k}}\preceq\xi\mathbf{I}_{r_{k}}\right\} }\right]\nonumber \\
 & \cdot\mathbf{A}^{T}\mathbf{T}\mathbf{A}+\mathbf{A}^{T}\mathbb{E}\left[\vphantom{\left(\mathbf{R}^{c}+\left(\mathbf{B}^{c}\left(k\right)\right)^{T}\mathbf{T}_{k}^{\mathrm{ker}}\mathbf{B}^{c}\left(k\right)\right)^{-1}}\right.\mathds{1}_{\left\{ r_{k}\neq0,\left(\boldsymbol{\Sigma}_{k}\right)_{r_{k}}\succ\xi\mathbf{I}_{r_{k}}\right\} }\bigg(\mathbf{T}_{k}^{\mathrm{ker}}-\mathbf{T}_{k}^{\mathrm{ker}}\nonumber \\
 & \cdot\mathbf{B}^{c}\left(k\right)\left(\mathbf{R}^{c}+\left(\mathbf{B}^{c}\left(k\right)\right)^{T}\mathbf{T}_{k}^{\mathrm{ker}}\mathbf{B}^{c}\left(k\right)\right)^{-1}\left(\mathbf{B}^{c}\left(k\right)\right)^{T}\nonumber \\
 & \cdot\mathbf{T}_{k}^{\mathrm{ker}}\bigg)\left.\vphantom{\vphantom{\left(\mathbf{R}^{c}+\left(\mathbf{B}^{c}\left(k\right)\right)^{T}\mathbf{T}_{k}^{\mathrm{ker}}\mathbf{B}^{c}\left(k\right)\right)^{-1}}}\right]\mathbf{A}+\mathbf{A}^{T}\mathbb{E}\left[\mathds{1}_{\left\{ r_{k}\neq0,\left(\boldsymbol{\Sigma}_{k}\right)_{r_{k}}\succ\xi\mathbf{I}_{r_{k}}\right\} }\mathbf{T}_{k}^{\mathrm{\mathbf{0}}}\right]\mathbf{A}+\mathbf{Q},\nonumber 
\end{align}
where $\xi>0$ is a tunable constant. 

We next analyze each term in (\ref{eq:Expanded substitute}). Denote
$\mathbf{T}_{r_{k}}=\left(\mathbf{U}_{k}\mathbf{T}\mathbf{U}_{k}^{T}\right)_{r_{k}}$
and $\overline{\mathbf{T}}_{r_{k}}=\left(\boldsymbol{\Sigma}_{k}\right)_{r_{k}}^{\frac{1}{2}}\left(\mathbf{U}_{k}\mathbf{T}\mathbf{U}_{k}^{T}\right)_{r_{k}}\left(\boldsymbol{\Sigma}_{k}\right)_{r_{k}}^{\frac{1}{2}}.$
Substitute (\ref{eq:T-ker-closed-form}) into (\ref{eq:Expanded substitute}),
the $\mathbf{T}_{k}^{\mathrm{ker}}$ dependent terms in (\ref{eq:Expanded substitute})
thus can be represented as
\begin{align}
 & \mathbf{T}_{k}^{\mathrm{ker}}-\mathbf{T}_{k}^{\mathrm{ker}}\mathbf{B}^{c}\left(k\right)\left(\mathbf{R}^{c}+\left(\mathbf{B}^{c}\left(k\right)\right)^{T}\mathbf{T}_{k}^{\mathrm{ker}}\mathbf{B}^{c}\left(k\right)\right)^{-1}\label{eq:T-ker dependence}\\
 & \cdot\left(\mathbf{B}^{c}\left(k\right)\right)^{T}\mathbf{T}_{k}^{\mathrm{ker}}=\mathbf{U}_{k}^{T}\mathbf{\Xi}_{k}\mathbf{U}_{k},\nonumber 
\end{align}
where
\begin{align}
 & \mathbf{\Xi}_{k}=\left[\begin{array}{cc}
\left(\left(\boldsymbol{\Sigma}_{k}\right)_{r_{k}}+\mathbf{T}_{r_{k}}^{-1}\right)^{-1} & \left(\mathbf{I}+\overline{\mathbf{T}}_{r_{k}}^{-1}\right)^{-1}\mathbf{T}_{r_{k}}\mathbf{L}_{k}\\
\mathbf{L}_{k}^{T}\mathbf{T}_{r_{k}}\left(\mathbf{I}+\overline{\mathbf{T}}_{r_{k}}^{-1}\right)^{-1} & \widetilde{\mathbf{\Xi}}_{k}
\end{array}\right],\\
 & \widetilde{\mathbf{\Xi}}_{k}=\mathbf{L}_{k}^{T}\mathbf{T}_{r_{k}}\mathbf{L}_{k}-\mathbf{L}_{k}^{T}\mathbf{T}_{r_{k}}\left(\mathbf{I}+\overline{\mathbf{T}}_{r_{k}}^{-1}\right)^{-1}\mathbf{T}_{r_{k}}\mathbf{L}_{k}.
\end{align}

Substitute (\ref{eq:T-0-closed-form}) into (\ref{eq:Expanded substitute}),
the $\mathbf{T}_{k}^{\mathrm{\mathbf{0}}}$ dependent terms in (\ref{eq:Expanded substitute})
can be represented as
\begin{align}
 & \mathbf{A}^{T}\mathbf{T}_{k}^{\mathrm{\mathbf{0}}}\mathbf{A}=\mathbf{A}^{T}\widetilde{\mathbf{U}}_{k}^{T}\widetilde{\mathbf{U}}_{k}\mathbf{T}\widetilde{\mathbf{U}}_{k}^{T}\widetilde{\mathbf{U}}_{k}\mathbf{A}\label{eq:T-0 dependence}\\
 & -\mathbf{A}^{T}\mathbf{U}_{k}^{T}\left[\begin{array}{cc}
\mathbf{0}_{r_{k}} & \mathbf{0}_{r_{k}\times\left(S-r_{k}\right)}\\
\mathbf{0}_{\left(S-r_{k}\right)\times r_{k}} & \mathbf{L}_{k}^{T}\mathbf{T}_{r_{k}}\mathbf{L}_{k}
\end{array}\right]\mathbf{U}_{k}\mathbf{A},\nonumber 
\end{align}
where we define
\begin{align}
 & \widetilde{\mathbf{U}}_{k}\triangleq\mathrm{diag}\left(\left[\boldsymbol{0}_{r_{k}},\mathbf{I}_{S-r_{k}}\right]\right)\mathbf{U}_{k}
\end{align}
for conciseness.

Now substitute (\ref{eq:T-ker dependence}) and (\ref{eq:T-0 dependence})
into (\ref{eq:Expanded substitute}), we can obtain an upper bound
of the L.H.S. of (\ref{eq: g(T)<T}) as follows

\begin{align}
 & \mathbb{E}\left[g\left(\mathbf{T};k\right)\right]+\mathbf{Q}\preceq\mathbb{E}\left[\mathds{1}_{\left\{ r_{k}\neq0,\left(\boldsymbol{\Sigma}_{k}\right)_{r_{k}}\preceq\xi\mathbf{I}_{r_{k}}\right\} }\right]\mathbf{A}^{T}\mathbf{T}\mathbf{A}\label{eq: g(T) substitute <T}\\
 & +\mathbb{E}\left[\mathbf{A}^{T}\widetilde{\mathbf{U}}_{k}^{T}\widetilde{\mathbf{U}}_{k}\mathbf{T}^{*}\widetilde{\mathbf{U}}_{k}^{T}\widetilde{\mathbf{U}}_{k}\mathbf{A}\right]\nonumber \\
 & +\mathbf{A}^{T}\mathbb{E}\bigg[\mathds{1}_{\left\{ r_{k}\neq0,\left(\boldsymbol{\Sigma}_{k}\right)_{r_{k}}\succ\xi\mathbf{I}_{r_{k}}\right\} }\mathbf{U}_{k}^{T}\nonumber \\
 & \cdot\bigg(\mathbf{\Xi}_{k}-\left[\begin{array}{cc}
\mathbf{0}_{r_{k}} & \mathbf{0}_{r_{k}\times\left(S-r_{k}\right)}\\
\mathbf{0}_{\left(S-r_{k}\right)\times r_{k}} & \mathbf{L}_{k}^{T}\mathbf{T}_{r_{k}}\mathbf{L}_{k}
\end{array}\right]\bigg)\mathbf{U}_{k}\bigg]\mathbf{A}+\mathbf{Q}.\nonumber 
\end{align}

Note that 
\begin{align}
 & \widetilde{\mathbf{\Xi}}_{k}-\mathbf{L}_{k}^{T}\mathbf{T}_{r_{k}}\mathbf{L}_{k}=-\mathbf{L}_{k}^{T}\mathbf{T}_{r_{k}}\left(\mathbf{I}+\overline{\mathbf{T}}_{r_{k}}^{-1}\right)^{-1}\mathbf{T}_{r_{k}}\mathbf{L}_{k}\prec\boldsymbol{0},
\end{align}
we have
\begin{align}
 & \mathds{1}_{\left\{ r_{k}\neq0,\left(\boldsymbol{\Sigma}_{k}\right)_{r_{k}}\succ\xi\mathbf{I}_{r_{k}}\right\} }\bigg(\mathbf{\Xi}_{k}-\left[\begin{array}{cc}
\mathbf{0}_{r_{k}} & \mathbf{0}_{r_{k}\times\left(S-r_{k}\right)}\\
\mathbf{0}_{\left(S-r_{k}\right)\times r_{k}} & \mathbf{L}_{k}^{T}\mathbf{T}_{r_{k}}\mathbf{L}_{k}
\end{array}\right]\bigg)\label{eq:upper bound-T_ker}\\
 & \prec\mathds{1}_{\left\{ r_{k}\neq0,\left(\boldsymbol{\Sigma}_{k}\right)_{r_{k}}\succ\xi\mathbf{I}_{r_{k}}\right\} }\mathrm{Tr}\left(\left(\boldsymbol{\Sigma}_{k}\right)_{r_{k}}^{-1}\right)\mathbf{I}.\nonumber 
\end{align}

Combining (\ref{eq: g(T) substitute <T}) and (\ref{eq:upper bound-T_ker}),
we can conclude that if we can find a $\mathbf{T}^{*}\in\mathbb{S}_{+}^{S}$
that satisfies the following equality
\begin{align}
 & \text{\ensuremath{\mathrm{Pr}}}\left(r_{k}\neq0,\left(\boldsymbol{\Sigma}_{k}\right)_{r_{k}}\preceq\xi\mathbf{I}_{r_{k}}\right)\mathbf{A}^{T}\mathbf{T}^{*}\mathbf{A}\label{eq:T Lya Eq}\\
 & +\mathbb{E}\left[\mathbf{A}^{T}\widetilde{\mathbf{U}}_{k}^{T}\widetilde{\mathbf{U}}_{k}\mathbf{T}^{*}\widetilde{\mathbf{U}}_{k}^{T}\widetilde{\mathbf{U}}_{k}\mathbf{A}\right]\nonumber \\
 & +\left\Vert \mathbf{A}\right\Vert ^{2}\left\Vert \mathbb{E}\left[\mathds{1}_{\left\{ r_{k}\neq0,\left(\boldsymbol{\Sigma}_{k}\right)_{r_{k}}\succ\xi\mathbf{I}_{r_{k}}\right\} }\mathrm{Tr}\left(\left(\boldsymbol{\Sigma}_{k}\right)_{r_{k}}^{-1}\right)\right]\mathbf{I}+\mathbf{Q}\right\Vert \mathbf{I}\nonumber \\
 & =\mathbf{T}^{*},\nonumber 
\end{align}
then $\mathbf{T}=\mathbf{T}^{*}$ will automatically satisfy condition
(\ref{eq: g(T) substitute <T}). 

Vectorizing both sides of (\ref{eq:T Lya Eq}), we obtain 
\begin{align}
 & \mathrm{vec\left(\mathbf{T}^{*}\right)}=\widetilde{\mathbf{A}}\mathrm{vec\left(\mathbf{T}^{*}\right)}+\left\Vert \mathbf{A}\right\Vert ^{2}\label{eq:T Lya Eq-1}\\
 & \cdot\mathrm{vec}\left(\left\Vert \mathbb{E}\left[\mathds{1}_{\left\{ r_{k}\neq0,\left(\boldsymbol{\Sigma}_{k}\right)_{r_{k}}\succ\xi\mathbf{I}_{r_{k}}\right\} }\mathrm{Tr}\left(\left(\boldsymbol{\Sigma}_{k}\right)_{r_{k}}^{-1}\right)\right]\mathbf{I}+\mathbf{Q}\right\Vert \mathbf{I}\right),\nonumber 
\end{align}
where $\widetilde{\mathbf{A}}$ is given by 
\begin{align}
\widetilde{\mathbf{A}}= & \text{\ensuremath{\mathrm{Pr}}}\left(r_{k}\neq0,\left(\boldsymbol{\Sigma}_{k}\right)_{r_{k}}\preceq\xi\mathbf{I}_{r_{k}}\right)\mathbf{A}\otimes\mathbf{A}\label{eq:GLyaEq-A-tilde}\\
 & +\mathbb{E}\left[\left(\widetilde{\mathbf{U}}_{k}^{T}\widetilde{\mathbf{U}}_{k}\mathbf{A}\right)\otimes\left(\widetilde{\mathbf{U}}_{k}^{T}\widetilde{\mathbf{U}}_{k}\mathbf{A}\right)\right].\nonumber 
\end{align}

Equation (\ref{eq:T Lya Eq-1}) is a generalized generalized Lyapunov
equation, and the existence of a positive definite solution $\mathbf{T}^{*}$
that satisfies (\ref{eq:T Lya Eq}) is closely related to the solvability
of the corresponding system of linear equations in (\ref{eq:T Lya Eq-1}).
Based on the standard Lyapunov stability theory, we can conclude that
the generalized Lyapunov equation (\ref{eq:T Lya Eq}) has a positive
definite solution $\mathbf{T}^{*}$ if $\rho\left(\widetilde{\mathbf{A}}\right)<1,$Besides,
under the condition $\rho\left(\widetilde{\mathbf{A}}\right)<1,$
the solution $\mathbf{T}^{*}$ to equality (\ref{eq:T Lya Eq}) is
given by (\ref{eq:T*-Closed-form}). 

We next show that condition (\ref{eq: =00005CrhomaxA_k}) is sufficient
to guarantee $\rho\left(\widetilde{\mathbf{A}}\right)<1$. Note that
for every $\varepsilon>0$, there exists an operator norm $\left\Vert \cdot\right\Vert _{\varepsilon}$
such that 
\begin{align}
 & \left\Vert \mathbb{E}\left[\left(\widetilde{\mathbf{U}}_{k}^{T}\widetilde{\mathbf{U}}_{k}\mathbf{A}\right)\otimes\left(\widetilde{\mathbf{U}}_{k}^{T}\widetilde{\mathbf{U}}_{k}\mathbf{A}\right)\right]\right\Vert _{\varepsilon}\label{eq:epsilon-norm}\\
 & =\text{\ensuremath{\rho}\ensuremath{\left(\mathbb{E}\left[\left(\widetilde{\mathbf{U}}_{k}^{T}\widetilde{\mathbf{U}}_{k}\mathbf{A}\right)\otimes\left(\widetilde{\mathbf{U}}_{k}^{T}\widetilde{\mathbf{U}}_{k}\mathbf{A}\right)\right]\right)}}+\varepsilon.\nonumber 
\end{align}
As a result, under condition (\ref{eq: =00005CrhomaxA_k}), we are
able to choose the positive constant $\varepsilon$ in (\ref{eq:epsilon-norm})
to be sufficiently small such that 
\begin{align}
 & \left\Vert \mathbb{E}\left[\left(\widetilde{\mathbf{U}}_{k}^{T}\widetilde{\mathbf{U}}_{k}\mathbf{A}\right)\otimes\left(\widetilde{\mathbf{U}}_{k}^{T}\widetilde{\mathbf{U}}_{k}\mathbf{A}\right)\right]\right\Vert _{\varepsilon}<1.\label{eq:epsilon-matirx-norm}
\end{align}

Applying Gelfand's formula to $\widetilde{\mathbf{A}}$ in (\ref{eq:GLyaEq-A-tilde})
using the $\left\Vert \cdot\right\Vert _{\varepsilon}$ in (\ref{eq:epsilon-norm}),
we obtain
\begin{align}
 & \rho\left(\widetilde{\mathbf{A}}\right)=\lim_{k\rightarrow\infty}\left\Vert \widetilde{\mathbf{A}}^{k}\right\Vert _{\varepsilon}^{\frac{1}{k}}<\left\Vert \widetilde{\mathbf{A}}\right\Vert _{\varepsilon}\label{eq:eq:GLyaEq-rho}\\
 & <\mathrm{Pr}\left(r_{k}\neq0,\left(\boldsymbol{\Sigma}_{k}\right)_{r_{k}}\preceq\xi\mathbf{I}_{r_{k}}\right)\left\Vert \mathbf{A}\otimes\mathbf{A}\right\Vert _{\varepsilon}\nonumber \\
 & +\left\Vert \mathbb{E}\left[\left(\widetilde{\mathbf{U}}_{k}^{T}\widetilde{\mathbf{U}}_{k}\mathbf{A}\right)\otimes\left(\widetilde{\mathbf{U}}_{k}^{T}\widetilde{\mathbf{U}}_{k}\mathbf{A}\right)\right]\right\Vert _{\varepsilon}.\nonumber 
\end{align}

Note that $\lim_{\xi\rightarrow0}\mathrm{Pr}\left(r_{k}\neq0,\left(\boldsymbol{\Sigma}_{k}\right)_{r_{k}}\preceq\xi\mathbf{I}_{r_{k}}\right)=0$
and the positive constant $\xi$ in (\ref{eq:eq:GLyaEq-rho}) is tunable.
Under condition (\ref{eq: =00005CrhomaxA_k}), we can choose a sufficiently
small $\xi>0$ such that 
\begin{align}
 & \mathrm{Pr}\left(r_{k}\neq0,\left(\boldsymbol{\Sigma}_{k}\right)_{r_{k}}\preceq\xi\mathbf{I}_{r_{k}}\right)<\left\Vert \mathbf{A}\otimes\mathbf{A}\right\Vert _{\varepsilon}^{-1}\label{eq: choice of kesi-1}\\
 & \cdot\left(1-\left\Vert \mathbb{E}\left[\left(\widetilde{\mathbf{U}}_{k}^{T}\widetilde{\mathbf{U}}_{k}\mathbf{A}\right)\otimes\left(\widetilde{\mathbf{U}}_{k}^{T}\widetilde{\mathbf{U}}_{k}\mathbf{A}\right)\right]\right\Vert _{\varepsilon}\right).\nonumber 
\end{align}

Substitute (\ref{eq: choice of kesi-1}) into (\ref{eq:eq:GLyaEq-rho}),
it follows that $\rho\left(\widetilde{\mathbf{A}}\right)<1$. This
completes the proof of Lemma \ref{Lemma: Suff-Condition-Existence of T}.

\subsection{\label{subsec:Proof-of-Coro: Example 1 Suff}Proof of Corollary \ref{Coro: Example 1 Suff}}

In the system setup of Example \ref{Exmp: E.g.1}, we have $\widetilde{\mathbf{U}}_{k}^{T}\widetilde{\mathbf{U}}_{k}\mathbf{A}=\mathbf{A}$
$w.p.\ \left(1-\delta\right)$ and $\widetilde{\mathbf{U}}_{k}^{T}\widetilde{\mathbf{U}}_{k}\mathbf{A}=\mathbf{0}$
$w.p.\ \delta.$ As a result, we obtain
\begin{align}
 & \rho\left(\mathbb{E}\left[\left(\widetilde{\mathbf{U}}_{k}^{T}\widetilde{\mathbf{U}}_{k}\mathbf{A}\right)\otimes\left(\widetilde{\mathbf{U}}_{k}^{T}\widetilde{\mathbf{U}}_{k}\mathbf{A}\right)\right]\right)\\
 & =\rho\left(\left(1-\delta\right)\mathbf{A}\otimes\mathbf{A}\right)=\left(1-\delta\right)\rho^{2}\left(\mathbf{A}\right).\nonumber 
\end{align}
As a result, condition (\ref{eq: =00005CrhomaxA_k}) is reduced to
condition (\ref{eq: delta condition-suff-example-3}), which is a
part of sufficient condition because of Theorem \ref{Thm: Verifiable Suff Condition}. 

Since the minimum non-zero singular value of $\mathbf{B}^{c}\left(k\right)$
in Example \ref{Exmp: E.g.1} is bounded away from 0 w. p. 1., we
can choose $\xi$ in Lemma \ref{Lemma: Suff-Condition-Existence of T}
to be sufficiently small such that $\mathrm{Pr}\left(\left(\boldsymbol{\Sigma}_{k}\right)_{r_{k}}\succ\xi\mathbf{I}_{r_{k}}\right)=1$.
As a result, condition (\ref{eq: g(T)<T}) in Theorem \ref{Thm: g(T)<T}
is reduced to find a positive definite $\mathbf{T}$ such that 
\begin{align}
 & \left(1-\delta\right)\mathbf{A}\mathbf{T}\mathbf{A}+\left\Vert \mathbf{A}\right\Vert ^{2}\mathbb{E}\left[\mathds{1}_{\left\{ r_{k}\neq0\right\} }\mathrm{Tr}\left(\left(\boldsymbol{\Sigma}_{k}\right)_{r_{k}}^{-1}\right)\right]\mathbf{I}+\mathbf{Q}\prec\mathbf{T}.\label{eq:Example 1- suff}
\end{align}
We can construct a $\mathbf{T}^{*}$ given by
\begin{align}
\mathbf{T}^{*}= & \left(1+\gamma_{2}\right)\left(1-\left(1-\delta\right)\rho^{2}\left(\mathbf{A}\right)\right)^{-1}\label{eq:T* constructed}\\
 & \cdot\left\Vert \mathbb{E}\left[\mathds{1}_{\left\{ r_{k}\neq0\right\} }\mathrm{Tr}\left(\left(\boldsymbol{\Sigma}_{k}\right)_{r_{k}}^{-1}\right)\right]\mathbf{I}+\mathbf{Q}\right\Vert \mathbf{I},\nonumber 
\end{align}
where $\gamma_{2}>0$ is some constant. It is clear that $\mathbf{T}^{*}$
in (\ref{eq:T* constructed}) satisfies equality (\ref{eq:Example 1- suff}). 

We next analyze the second condition in Theorem \ref{Thm: Verifiable Suff Condition}
based on the $\mathbf{T}^{*}$ in (\ref{eq:T* constructed}). Note
that 
\begin{align}
 & \left(\mathbb{E}\left[g\left(\mathbf{\mathbf{T}^{*}};k\right)\right]+\mathbf{Q}\right)^{-1}-\left(\mathbf{\mathbf{T}^{*}}\right)^{-1}\label{eq:E=00005B=00005D-T*}\\
 & \succ\left\Vert \mathbb{E}\left[\mathds{1}_{\left\{ r_{k}\neq0\right\} }\mathrm{Tr}\left(\left(\boldsymbol{\Sigma}_{k}\right)_{r_{k}}^{-1}\right)\right]\mathbf{I}+\mathbf{Q}\right\Vert ^{-1}\nonumber \\
 & \cdot\left(1-\left(1-\delta\right)\rho^{2}\left(\mathbf{A}\right)\right)\left[\left(1+\gamma_{2}\rho^{2}\left(\mathbf{A}\right)\right)^{-1}-\left(1+\gamma_{2}\right)^{-1}\right]\mathbf{I}\nonumber \\
 & =\gamma_{3}\left(1-\left(1-\delta\right)\rho^{2}\left(\mathbf{A}\right)\right)^{2}\mathbf{I},\nonumber 
\end{align}
where
\begin{align}
\gamma_{3}= & \gamma_{2}\left(\gamma_{2}+1\right)^{-1}\left(1+\gamma_{2}\rho^{2}\left(\mathbf{A}\right)\right)^{-1}\\
 & \cdot\left\Vert \mathbb{E}\left[\mathds{1}_{\left\{ r_{k}\neq0\right\} }\mathrm{Tr}\left(\left(\boldsymbol{\Sigma}_{k}\right)_{r_{k}}^{-1}\right)\right]\mathbf{I}+\mathbf{Q}\right\Vert ^{-1}.\nonumber 
\end{align}

Substitute (\ref{eq:E=00005B=00005D-T*}) back into (\ref{eq:Ra}),
we obtain the sufficient requirement on $\mathbf{R}^{a}$ as
\begin{align}
 & \mathbf{R}^{a}\succ\gamma_{3}^{-1}\left(1-\left(1-\delta\right)\rho^{2}\left(\mathbf{A}\right)\right)^{-2}\mathbb{E}\left[\left(\mathbf{B}^{a}\left(k\right)\right)^{T}\right]\label{eq:requirement Ra}\\
 & \cdot\mathbb{E}\left[\left(\mathbf{B}^{a}\left(k\right)\right)\right]+N_{t}^{a}\left|\mathcal{N}^{a}\right|\sigma_{\mathbf{B}^{a}}^{2}\left(1+\gamma_{2}\right)\left(1-\left(1-\delta\right)\rho^{2}\left(\mathbf{A}\right)\right)^{-1}\nonumber \\
 & \cdot\left\Vert \mathbb{E}\left[\mathds{1}_{\left\{ r_{k}\neq0\right\} }\mathrm{Tr}\left(\left(\boldsymbol{\Sigma}_{k}\right)_{r_{k}}^{-1}\right)\right]\mathbf{I}+\mathbf{Q}\right\Vert \mathbf{I}.\nonumber 
\end{align}

Therefore, to satisfy requirement \ref{eq:requirement Ra}, it suffices
to choose\\ $\mathbf{R}^{a}\succ\gamma_{1}\left(1-\left(1-\delta\right)\rho^{2}\left(\mathbf{A}\right)\right)^{-2}\mathbf{I},$
where $\gamma_{1}$ satisfies
\begin{align}
\gamma_{1}\geq & \gamma_{3}^{-1}\left\Vert \mathbb{E}\left[\left(\mathbf{B}^{a}\left(k\right)\right)\right]\right\Vert ^{2}+N_{t}^{a}\left|\mathcal{N}^{a}\right|\sigma_{\mathbf{B}^{a}}^{2}\left(1+\gamma_{2}\right)\\
 & \cdot\left\Vert \mathbb{E}\left[\mathds{1}_{\left\{ r_{k}\neq0\right\} }\mathrm{Tr}\left(\left(\boldsymbol{\Sigma}_{k}\right)_{r_{k}}^{-1}\right)\right]\mathbf{I}+\mathbf{Q}\right\Vert .\nonumber 
\end{align}
This completes the proof of Corollary \ref{Coro: Example 1 Suff}.

\subsection{\label{subsec:Proof-of-Coro: Example 1 Necessary}Proof of Corollary
\ref{Coro: Example 1 Necessary} }

We first analyze the necessary requirement on $\delta$. In the system
setup of Example \ref{Exmp: E.g.1}, we know from Theorem \ref{Thm: General case- suff and necessary}
that if $\mathcal{P}_{\mathbf{R}^{a}}\cap\mathcal{R}_{\mathbf{R}^{a}}\neq\textrm{Ø},$
then
\begin{align}
 & \mathbf{P}^{*}=f\left(\mathbf{P}^{*}\right)\succ\left(1-\delta\right)\mathbf{A}\mathbf{P}^{*}\mathbf{A}+\label{eq:P* DARE}\\
 & \mathbb{E}\left[\mathbf{A}^{T}\left(\mathds{1}_{\left\{ r_{k}\neq0\right\} }\mathbf{B}^{c}\left(k\right)\left(\mathbf{R}^{c}\right)^{-1}\left(\mathbf{B}^{c}\left(k\right)\right)^{T}+\left(\mathbf{P}^{*}\right)^{-1}\right)^{-1}\mathbf{A}\right]\nonumber \\
 & +\mathbf{Q}\succ\left(1-\delta\right)\mathbf{A}\mathbf{P}^{*}\mathbf{A}+\mathbf{Q}.\nonumber 
\end{align}
This implies $\left(1-\delta\right)\rho^{2}\left(\mathbf{A}\right)<1$,
which is exactly condition (\ref{eq: delta condition-suff-example-3}). 

We now analyze the necessary requirement on $\mathbf{R}^{a}$ under
condition (\ref{eq: delta condition-suff-example-3}). Let the eigenvalue
decomposition of $\mathbf{A}$ be $\mathbf{A}=\mathbf{V}\mathbf{\Lambda}_{\mathbf{A}}\mathbf{V}^{-1},$
where $\mathbf{V}\in\mathbb{R}^{S\times S}$ is an orthogonal matrix
and $\mathbf{\Lambda}_{\mathbf{A}}$ is a diagonal matrix with the
diagonal elements in descending order. Based on (\ref{eq:P* DARE}),
we have
\begin{align}
 & \mathbf{P}^{*}\succ\sigma_{\mathrm{min}}\left(\mathbf{V}^{-1}\mathbf{Q}\mathbf{V}\right)\mathbf{V}\left(\mathbf{I}-\left(1-\delta\right)\mathbf{\Lambda}_{\mathbf{A}}^{2}\right)^{-1}\mathbf{V}.\label{eq:P*>diag-1}
\end{align}

Since $\mathcal{R}_{\mathbf{R}^{a}}\neq\textrm{Ø}$, there exists
a sufficiently small $\mathbf{\Theta}\in\mathbb{S}_{+}^{S}$ such
that 
\begin{align}
 & \left(\mathbf{P}^{*}\right)^{-1}-\mathbb{E}\left[\mathbf{B}^{a}\left(k\right)\right]\left(\mathbf{R}^{a}-\mathrm{Cov}\left(\mathbf{B}^{a}\left(k\right)\left(\mathbf{P}^{*}\right)^{\frac{1}{2}}\right)\right)^{-1}\label{eq:Theta P*}\\
 & \cdot\mathbb{E}\left[\left(\mathbf{B}^{a}\left(k\right)\right)^{T}\right]=\mathbf{\Theta}.\nonumber 
\end{align}
As a result, we have 
\begin{align}
 & f\left(\mathbf{P}^{*}\right)=\left(1-\delta\right)\mathbf{A}\mathbf{\Theta}^{-1}\mathbf{A}\\
 & +\mathbb{E}\left[\mathbf{A}^{T}\left(\mathds{1}_{\left\{ \mathbf{B}^{c}\left(k\right)\neq\mathbf{0}\right\} }\mathbf{B}^{c}\left(k\right)\left(\mathbf{R}^{c}\right)^{-1}\left(\mathbf{B}^{c}\left(k\right)\right)^{T}+\mathbf{\Theta}\right)^{-1}\mathbf{A}\right]\nonumber \\
 & +\mathbf{Q}=\bigg(\mathbf{\Theta}+\mathbb{E}\left[\mathbf{B}^{a}\left(k\right)\right]\left(\mathbf{R}^{a}-\mathrm{Cov}\left(\mathbf{B}^{a}\left(k\right)\left(\mathbf{P}^{*}\right)^{\frac{1}{2}}\right)\right)^{-1}\nonumber \\
 & \cdot\mathbb{E}\left[\left(\mathbf{B}^{a}\left(k\right)\right)^{T}\right]\bigg)^{-1}.\nonumber 
\end{align}

Note that $\mathbf{A}$ is invertible, it follows that
\begin{align}
 & \left(1-\delta\right)^{-1}\mathbf{A}^{-1}\mathbf{\Theta}\mathbf{A}^{-1}-\mathbf{\Theta}\\
 & \succ\mathbb{E}\left[\mathbf{B}^{a}\left(k\right)\right]\left(\mathbf{R}^{a}-\mathrm{Cov}\left(\mathbf{B}^{a}\left(k\right)\left(\mathbf{P}^{*}\right)^{\frac{1}{2}}\right)\right)^{-1}\mathbb{E}\left[\left(\mathbf{B}^{a}\left(k\right)\right)^{T}\right]\nonumber 
\end{align}
This is equivalent to 
\begin{align}
\mathbf{R}^{a}\succ & \mathbb{E}\left[\mathbf{B}^{a}\left(k\right)^{T}\right]\mathbf{V}\widetilde{\mathbf{\Theta}}^{-1}\mathbf{V}^{-1}\mathbb{E}\left[\mathbf{B}^{a}\left(k\right)\right]\label{eq:Ra > Theta-tilde}\\
 & +\mathrm{Cov}\left(\mathbf{B}^{a}\left(k\right)\left(\mathbf{P}^{*}\right)^{\frac{1}{2}}\right),\nonumber 
\end{align}
where $\widetilde{\mathbf{\Theta}}=\left(1-\delta\right)^{-1}\mathbf{\Lambda}_{\mathbf{A}}^{-1}\mathbf{V}^{-1}\mathbf{\Theta}\mathbf{V}\mathbf{\Lambda}_{\mathbf{A}}^{-1}-\mathbf{V}^{-1}\mathbf{\Theta}\mathbf{V}.$
We have the following relationship on the $i$-th diagonal elements
of $\widetilde{\mathbf{\Theta}},\forall1\leq i\leq S,$
\begin{align}
 & \widetilde{\mathbf{\Theta}}_{\left(i,i\right)}=\left(\left(1-\delta\right)^{-1}\sigma_{i}^{-2}\left(\mathbf{A}\right)-1\right)\left(\mathbf{V}^{-1}\mathbf{\Theta}\mathbf{V}\right)_{\left(i,i\right)}.\label{eq:Theta-Tilde}
\end{align}

From (\ref{eq:Theta P*}), we know that the diagonal elements satisfy
$\left(\mathbf{\mathbf{V}^{-1}\mathbf{\Theta}\mathbf{V}}\right)_{\left(i,i\right)}<\left(\mathbf{V}^{-1}\left(\mathbf{P}^{*}\right)^{-1}\mathbf{V}\right)_{\left(i,i\right)},\forall1\leq i\leq S.$
From (\ref{eq:P*>diag-1}), we have 
\begin{align}
 & \left(\mathbf{V}^{-1}\left(\mathbf{P}^{*}\right)^{-1}\mathbf{V}\right)_{\left(i,i\right)}\label{eq:P*>diag}\\
 & <\sigma_{\mathrm{min}}^{-1}\left(\mathbf{V}^{-1}\mathbf{Q}\mathbf{V}\right)\left(1-\left(1-\delta\right)\sigma_{i}^{2}\left(\mathbf{A}\right)\right)^{-1}.\nonumber 
\end{align}
Substitute (\ref{eq:P*>diag}) into (\ref{eq:Theta-Tilde}), we obtain
\begin{align}
 & \widetilde{\mathbf{\Theta}}{}_{\left(i,i\right)}\label{eq:Theta-tilde-ii}\\
 & <\sigma_{\mathrm{min}}^{-1}\left(\mathbf{V}^{-1}\mathbf{Q}\mathbf{V}\right)\sigma_{i}^{-2}\left(\mathbf{A}\right)\left(1-\left(1-\delta\right)\sigma_{i}^{2}\left(\mathbf{A}\right)\right)^{2}.\nonumber 
\end{align}
Since $\widetilde{\mathbf{\Theta}}$ is symmetric and positive definite,
we have 
\begin{align}
\left\Vert \widetilde{\mathbf{\Theta}}^{-1}\right\Vert  & \geq\mathrm{max}\left(\Big(\widetilde{\mathbf{\Theta}}^{-1}\Big)_{\left(i,i\right)}\right)\geq\mathrm{max}\left(\left(\widetilde{\mathbf{\Theta}}{}_{\left(i,i\right)}\right)^{-1}\right)\label{eq:Theta-tilde-L2-norm}\\
 & \geq\sigma_{\mathrm{min}}\left(\mathbf{V}^{-1}\mathbf{Q}\mathbf{V}\right)\rho^{2}\left(\mathbf{A}\right)\left(1-\left(1-\delta\right)\rho^{2}\left(\mathbf{A}\right)\right)^{-2}.\nonumber 
\end{align}

Applying Ostrwoski's theorem {[}50{]} to (\ref{eq:Ra > Theta-tilde})
and note that $\mathbf{V}^{-1}=\mathbf{V}^{T}$, we obtain 
\begin{align}
 & \left\Vert \mathbf{R}^{a}\right\Vert >\gamma_{5}\left\Vert \widetilde{\mathbf{\Theta}}^{-1}\right\Vert \\
 & \geq\gamma_{5}\sigma_{\mathrm{min}}\left(\mathbf{V}^{-1}\mathbf{Q}\mathbf{V}\right)\rho^{2}\left(\mathbf{A}\right)\left(1-\left(1-\delta\right)\rho^{2}\left(\mathbf{A}\right)\right)^{-2},\nonumber 
\end{align}
where $\gamma_{5}$ is a positive constant that lies in a closed interval
satisfying
\begin{align}
\gamma_{5}\in\Big[ & \sigma_{S}\left(\mathbf{V}^{T}\mathbb{E}\left[\mathbf{B}^{a}\left(k\right)\right]\mathbb{E}\left[\left(\mathbf{B}^{a}\left(k\right)\right)^{T}\right]\mathbf{V}\right)\\
 & ,\sigma_{1}\left(\mathbf{V}^{T}\mathbb{E}\left[\mathbf{B}^{a}\left(k\right)\right]\mathbb{E}\left[\left(\mathbf{B}^{a}\left(k\right)\right)^{T}\right]\mathbf{V}\right)\Big].\nonumber 
\end{align}

Let $\gamma_{4}=\gamma_{5}\sigma_{\mathrm{min}}\left(\mathbf{V}^{-1}\mathbf{Q}\mathbf{V}\right)\rho^{2}\left(\mathbf{A}\right)$,
we arrive at condition (\ref{eq:Ra L2-norm-Necessary Condition}).
This completes the proof Corollary \ref{Coro: Example 1 Necessary}.

\subsection{\label{subsec:Proof-of-Thm: achievability of NE}Proof of Theorem
\ref{Thm: achievability of NE}}

We prove Theorem \ref{Thm: achievability of NE} based on the completed
squares form of $J_{K}$ in Lemma \ref{Lemma: J_K representations}. 

We first consider a $\beta\left(T_{0}\right)$-length finite-time
horizon problem. We complete the squares w.r.t. $\mathbf{u}_{\Delta}^{c}\left(\mathbf{\mathbf{P}^{*}};k\right)$
in the time slots $k\in\left[0,T_{0}\right]$ and complete the squares
w.r.t. $\mathbf{u}_{\Delta}^{c}\left(f^{\beta\left(T_{0}\right)-k-1}\left(\mathbf{Q}\right);k\right)$
in the time slots $k\in\left[T_{0}+1,\beta\left(T_{0}\right)-1\right]$,
$\forall T_{0}\in\mathbb{Z}^{+}$. We have
\begin{align}
 & \inf_{\boldsymbol{\mu}_{\beta\left(T_{0}\right)}^{c}}\sup_{\boldsymbol{\mu}_{\beta\left(T_{0}\right)}^{c}}J_{\beta\left(T_{0}\right)}\left(\boldsymbol{\mu}_{\beta\left(T_{0}\right)}^{c},\text{\ensuremath{\boldsymbol{\mu}_{\beta\left(T_{0}\right)}^{a}}}\right)=\min_{\mathbf{u}_{\Delta}^{c}}\max_{\mathbf{u}_{\Delta}^{c}}\overline{J}_{\beta\left(T_{0}\right)}\label{eq:min max J(beta t-0)}\\
 & \overset{\left(\ref{eq:min max J(beta t-0)}.a\right)}{=}\frac{1}{\beta\left(T_{0}\right)}\Bigg\{\mathrm{Tr}\Bigg(\Bigg(\sum_{k=0}^{\beta\left(T_{0}\right)-T_{0}-1}f^{k}\left(\mathbf{Q}\right)+\sum_{k=0}^{T_{0}-1}\mathbf{P}^{*}\Bigg)\nonumber \\
 & \cdot\left(\mathbf{\mathbf{W}}+\sum_{i\in\mathcal{N}^{c}}\mathbf{B}_{i}^{T}\mathbf{V}\mathbf{B}_{i}\right)\Bigg)+\left\Vert \mathbf{x}\left(0\right)\right\Vert _{\left(\mathbf{\mathbf{P}^{*}}-\mathbf{Q}\right)}^{2}\nonumber \\
 & -\mathbb{E}\left[\left\Vert \mathbf{x}\left(T_{0}\right)\right\Vert _{\left(\mathbf{\mathbf{P}^{*}}-f^{\beta\left(T_{0}\right)-T_{0}}\left(\mathbf{Q}\right)\right)}^{2}\right]\Bigg\},\nonumber 
\end{align}
where the equality ($\ref{eq:min max J(beta t-0)}.a$) holds by choosing
$\mathbf{u}^{c}\left(k\right)=\overline{\mathbf{u}}^{c}\left(k\right)$
and $\mathbf{u}^{a}\left(k\right)=\overline{\mathbf{u}}^{a}\left(k\right)$
in Theorem (\ref{Thm: achievability of NE}). 

Recall that the definition of $\alpha\left(T_{0}\right)$ in (\ref{eq:alpha(T_0)})
implies that $\mathbb{E}\left[\left\Vert \mathbf{x}\left(T_{0}\right)\right\Vert _{2}^{2}\right]<\alpha\left(T_{0}\right)$.
Moreover, the definition of $\beta\left(T_{0}\right)$ in (\ref{eq:beta(T_0)})
implies that $\mathbb{E}\left[\left\Vert \mathbf{x}\left(T_{0}\right)\right\Vert _{\left(\mathbf{\mathbf{P}^{*}}-f^{\beta\left(T_{0}\right)-T_{0}}\left(\mathbf{Q}\right)\right)}^{2}\right]<1.$
Therefore, let $T_{0}\rightarrow\infty$, we have, for all $\left\{ \boldsymbol{\mu}_{\beta\left(\infty\right)}^{c},\boldsymbol{\mu}_{\beta\left(\infty\right)}^{a}\right\} \in\left(\mathcal{V}^{c}\times\mathcal{V}^{a}\right)_{\beta\left(\infty\right)},$
$\lim_{T_{0}\rightarrow\infty}\frac{1}{\beta\left(T_{0}\right)}\left\Vert \mathbf{x}\left(0\right)\right\Vert _{\left(\mathbf{\mathbf{P}^{*}}-\mathbf{Q}\right)}^{2}$
$=\lim_{T_{0}\rightarrow\infty}\frac{1}{\beta\left(T_{0}\right)}\left\Vert \mathbf{x}\left(T_{0}\right)\right\Vert _{\left(\mathbf{\mathbf{P}^{*}}-f^{\beta\left(T_{0}\right)-T_{0}}\left(\mathbf{Q}\right)\right)}^{2}$
$=0$ and 
\begin{align}
 & \lim_{T_{0}\rightarrow\infty}\frac{1}{\beta\left(T_{0}\right)}\mathrm{Tr}\Bigg(\Bigg(\sum_{k=0}^{\beta\left(T_{0}\right)-T_{0}-1}f^{k}\left(\mathbf{Q}\right)+\sum_{k=0}^{T_{0}-1}\mathbf{P}^{*}\Bigg)\\
 & \cdot\left(\mathbf{\mathbf{W}}+\sum_{i\in\mathcal{N}^{c}}\mathbf{B}_{i}^{T}\mathbf{V}\mathbf{B}_{i}\right)\Bigg)=\mathbf{P}^{*}\left(\mathbf{\mathbf{W}}+\sum_{i\in\mathcal{N}^{c}}\mathbf{B}_{i}^{T}\mathbf{V}\mathbf{B}_{i}\right)\nonumber \\
 & =J^{*}.\nonumber 
\end{align}
This proves that $\left\{ \overline{\boldsymbol{\mu}}_{\beta\left(\infty\right)}^{c},\overline{\boldsymbol{\mu}}_{\beta\left(\infty\right)}^{a}\right\} $
achieves the upper value of the zero-sum game in the class of admissible
policy pairs $\left(\mathcal{V}^{c}\times\mathcal{V}^{a}\right)_{\beta\left(\infty\right)}$.

Similarly, based on the $\underline{J}_{K}$ form in Lemma \ref{Lemma: J_K representations},
we still consider a $\beta\left(T_{0}\right)$-length finite-time
horizon problem, where we complete the squares w.r.t. $\mathbf{u}_{\Delta}^{a}\left(\mathbf{\mathbf{P}^{*}};k\right)$
in the time slots $k\in\left[0,T_{0}\right]$ and complete the squares
w.r.t. $\mathbf{u}_{\Delta}^{a}\left(f^{\beta\left(T_{0}\right)-k-1}\left(\mathbf{Q}\right);k\right)$
in the time slots $k\in\left[T_{0}+1,\beta\left(T_{0}\right)-1\right]$,
$\forall T_{0}\in\mathbb{Z}^{+}$. By letting $T_{0}\rightarrow\infty$,
using similar approaches, we can verify that $\left\{ \overline{\boldsymbol{\mu}}_{\beta\left(\infty\right)}^{c},\overline{\boldsymbol{\mu}}_{\beta\left(\infty\right)}^{a}\right\} $
also achieves the lower value of the zero-sum game in the class of
admissible policy pairs $\left(\mathcal{V}^{c}\times\mathcal{V}^{a}\right)_{\beta\left(\infty\right)}$.

\end{document}